\documentclass[12pt,twoside]{article}

\usepackage{amsgen,amsmath,amstext,amsbsy,amsopn,amsfonts,amssymb}

\usepackage{memoire}

\def\poi#1#2#3#4#5#6#7{\def\un{#5#6#7}\def\deux{#6#7}
\def\trois{#2#4} \def\cinq{#3#4#5}
\ifx\un\empty {#1}_{#2}{#3}{#1}_{#4} \else
\ifx\deux\empty {#5}(#1_{#2}){#3}{#5}(#1_{#4}) \else
\ifx\trois\empty {#5}_{#6}(#1){#3}{#5}_{#7}(#1) \else
{#5_{#6}}(#1_{#2}){#3}{#5_{#7}}(#1_{#4}) \fi \fi \fi}
\def\rond{\raisebox{.3mm}{\scriptsize$\circ$}}

\def\tens{\raisebox{.3mm}{\scriptsize$\otimes$}}

\def\dv#1#2{\langle {#1},{#2}\rangle}
    
\def\tk#1#2{{#2}\otimes _{#1}}

\def\ec#1#2#3#4#5{\def\un{#3#4#5}\def\deux{#3#5}\def\trois{#3}
\def\four{#2#4#5}\def\five{#2#5}\def\six{#2}\def\seven{#3#4}
\def\eight{#2#4} \def\nine{#2#3#4}
\ifx\nine\empty {\rm #1}_{#5} \else
\ifx\un\empty {\rm #1}({\goth #2}) \else
\ifx\deux\empty {\rm #1}({\goth #2}_{#4}) \else
\ifx\trois\empty {\rm #1}_{#5}({\goth #2}_{#4}) \else
\ifx\four\empty {\rm #1}(#3) \else
\ifx\five\empty {\rm #1}(#3_{#4}) \else
\ifx\six\empty {\rm #1}_{#5}(#3_{#4}) \else
\ifx\seven\empty {\rm #1}_{#5} ({\goth#2})\else                                               
\ifx\eight\empty {\rm #1}_{#5}({#3})                                               
\fi \fi \fi \fi \fi \fi \fi \fi \fi}
\def\hec#1#2#3#4#5{\def\un{#3#4#5}\def\deux{#3#5}\def\trois{#3}
\def\four{#2#4#5}\def\five{#2#5}\def\six{#2}\def\seven{#3#4}
\def\eight{#2#4} \def\nine{#2#3#4}
\ifx\nine\empty \hat{{\rm #1}}_{#5} \else
\ifx\un\empty \hat{{\rm #1}}({\goth #2}) \else
\ifx\deux\empty \hat{{\rm #1}}({\goth #2}_{#4}) \else
\ifx\trois\empty \hat{{\rm #1}}_{#5}({\goth #2}_{#4}) \else
\ifx\four\empty \hat{{\rm #1}}(#3) \else
\ifx\five\empty \hat{{\rm #1}}(#3_{#4}) \else
\ifx\six\empty \hat{{\rm #1}}_{#5}(#3_{#4}) \else 
\ifx\seven\empty \hat{{\rm #1}}_{#5} ({\goth#2})  \else                                             
\ifx\eight\empty \hat{{\rm #1}}_{#5}({#3}) 
\fi \fi \fi \fi \fi \fi \fi \fi \fi}
\def\e#1#2{\ec {#1}#2{}{}{}}
\def\es#1#2{\ec {#1}{}{#2}{}{}}

\def\ad{{\rm ad}\hskip .1em}

\def\j#1#2{\def\deux{#2} \ifx\deux\empty {\rm rk}\hskip .125em{{\goth #1}} \else {\rm rk}\hskip .125em{{\goth #1}_{#2}} \fi}
\def\k#1#2{\def\deux{#2} \ifx\deux\empty {\rm b}_{{\goth #1}} \else {\rm b}_{{\goth #1}_{#2}} \fi}

\def\an#1#2{\def\deux{#2} \ifx\deux\empty {\cal O}_{#1} \else {\cal O}_{#1,#2} \fi }
\def\han#1#2{\def\deux{#2} \ifx\deux\empty {\hat{{\cal O}}}_{#1} \else {\hat{{\cal O}}}_{#1,#2} \fi }

\def\r#1#2#3{\relax \def\un{#3} \relax \ifx\un\empty #1_{{\goth #2}} \else #1_{{\goth #2},{\goth #3}} \fi}
\def\u{\rm u}

\def\dim{{\rm dim}\hskip .125em}

\def\dd{{\rm d}}
\def\rk{{\rm rk}\hskip .125em}
\def\ad{{\rm ad}\hskip .1em}

\def\det{{\rm det}\hskip .125em}

\def\E{{\rm e}}
\def\r{{\rm r}}

\def\n{{\rm n}}
\def\s{{\rm s}}

\def\sy#1#2{{\rm S}^{#1}(#2)}

\def \ex #1#2{\mbox{\large$\wedge$}^{#1}(#2)}

\def\sl2{{\goth s}{\goth l}_{2}({\Bbb C})}

\def\A #1{{\rm Aut}({\goth #1})}

\renewcommand{\version}{{\rm Version r\'evis\'ee le~}\date} 
\renewcommand{\title}{\centerline{Complexe canonique de deuxi\`eme esp\`ece,}
\bigskip 
\centerline{vari\'et\'e commutante et bic\^one nilpotent 
d'une alg\`ebre de Lie r\'eductive.}}
\renewcommand{\author}{Jean-Yves CHARBONNEL}
\renewcommand{\institution}{Universit\'e Paris 7 - CNRS}
\renewcommand{\lastname}{J-Y. Charbonnel}

\renewcommand{\entries}{%
\nextref{Bo} {N. Bourbaki} {Groupes et Alg\`ebres de Lie.
Chapitres 4, 5 et 6} {Masson. (1981), Paris, New-York,
Barcelone, Milan, Mexico, Rio de Janeiro}

\nextref{Bo1} {N. Bourbaki} {Groupes et Alg\`ebres de Lie.
Chapitres 7 et 8} {Diffusion CCLS. (1975), Paris}

\nextref{Ca}{R. W. Carter}{Finite Groups of Lie Type. Conjugacy
classes and complex characters} {Pure and Applied Mathematics. A
Wiley-Interscience Series of Texts, Monographs \& Tracts
(1985), Chichester, New-York, Brisbane, Toronto, Singapore}

\nextref{Ch}{J-Y. Charbonnel}{Complexe canonique d'une alg\`ebre de Lie 
r\'eductive}{Preprint (2004)}

\nextref{De}{D. Deturck, H. Goldschmidt, J. Talvacchia}{Connections with 
prescribed curvature and Yang-Mills Currents: the semi-simple case}{Annales 
Scientifiques de l'\'Ecole Normale Sup\'erieure {\bf t. 24} (1991), p. 57--112}

\nextref{Dix}{J. Dixmier} {Champs de vecteurs adjoints sur les
groupes et alg\`ebres de Lie semi-simples} {Journal f\"ur die reine
und angewandte Mathematik {\bf Band. 309} (1979), p. 183--190}

\nextref{Gi}{W. L. Gan and V. Ginzburg}{Almost-commuting variety, 
${\cal D}$-modules, and Cherednik Algebras}{arXiv:math.RT/0409262}

\nextref{Ha}{R. Hartshorne} {Algebraic Geometry} {Graduate
Texts in Mathematics {\bf n$^{\circ}$ 52} (1977), Springer-Verlag, Berlin
Heidelberg New York}

\nextref{Ka}{R. Hotta and M. Kashiwara} {The invariant holonomic system
on a semi-simple Lie algebra} {Inventiones Mathematicae 
{\bf n$^{\circ}$ 75} (1984), p. 327--358}

\nextref{Ko}{B. Kostant} {The principal three-dimensional subgroup and
the Betti numbers of a complex simple Lie Group} {American Journal of
Mathematics {\bf n$^{\circ}$ 81} (1959), p. 973--1032}

\nextref{Ko1}{B. Kostant} {Lie group representations on polynomial rings} {American 
Journal of Mathematics {\bf n$^{circ}$ 85} (1963), p. 327--404}

\nextref{Le}{T. Levasseur and J. T. Stafford} {The Kernel of an
Homomorphism of Harish-Chandra}{Annales Scientifiques de l'\'Ecole
Normale Sup\'erieure {\bf t. 29} (1996), p. 385--397}

\nextref{Mat}{H. Matsumura} {Commutative ring theory} {Cambridge
studies in advanced mathematics {\bf n$^{\circ}$ 8} (1986), Cambridge
University Press, Cambridge, London, New York, New Rochelle, Melbourne,
Sydney}

\nextref{Na}{M. Nagata}{Local Rings}{Interscience Tracts in Pure and Applied
Mathematics {\bf n$^{\circ}$ 13} (1962), Interscience Publishers, New-York,
London, Sydney}

\nextref{Pa}{D. I. Panyushev} {The jacobian modules of a representation of a 
Lie algebra and geometry of commuting varieties} {Compositio Mathematica
{\bf n$^{\circ}$ 94} (1994), p. 181--199}

\nextref{Pr}{A. Premet} {Nilpotent commuting varieties of
reductive Lie algebras} {arXiv:math.RT/0302204}

\nextref{Ric} {R. W. Richardson} {Commuting varieties of
semisimple Lie algebras and algebraic groups} {Compositio
Mathematica {\bf vol. 38} (1979), p. 311--322}

\nextref{Sh} {I. R. Shafarevich} {Basic Algebraic Geometry 2} { 
(1994) Springer-Verlag, Berlin, Heidelberg, New York, London, Paris,
Tokyo, Hong-Kong, Barcelona, Budapest} 

\nextref{V} {F. D. Veldkamp} {The center of the universal
enveloping algebra of a Lie algebra in characteristic $p$} {Annales
Scientifiques de L'\'Ecole Normale Sup\'erieure {\bf n$^{\circ}$ 5}
(1972) , p. 217--240}

}

\renewcommand{\address}{Jean-Yves CHARBONNEL

Universit\'e Paris 7 - CNRS 

Institut de Math\'ematiques de Jussieu 

Th\'eorie des groupes

Case 7012

2 Place Jussieu

75251 Paris Cedex 05, France.

\makebox[1in][l]{jyc@math.jussieu.fr}

}

\begin{document}

\firstpage

\begin{abstract}
Let ${\goth g}$ be a finite dimensional complex reductive Lie algebra and
$\dv ..$ an invariant non degenerated bilinear form on 
${\goth g}\times {\goth g}$ which extends the Killing form of 
$[{\goth g},{\goth g}]$. We define a subcomplex $E_{\bullet}({\goth g})$ of 
the canonical complex $C_{\bullet}({\goth g})$ of ${\goth g}$. The space of 
$C_{\bullet}({\goth g})$ is the algebra 
$\tk {{\Bbb C}}{\e Sg}\tk {{\Bbb C}}{\e Sg}\ex {}{{\goth g}}$ where 
$\e Sg$ and $\ex {}{{\goth g}}$ are the symmetric and exterior algebras of 
${\goth g}$ and its differential is the
$\tk {{\Bbb C}}{\e Sg}\e Sg$-derivation which associates to the
element $v$ of ${\goth g}$ the function 
$(x,y)\mapsto \dv v{[x,y]}$ on ${\goth g}\times {\goth g}$. 
There exists a well defined sub-$\tk {{\Bbb C}}{\e Sg}\e Sg$-module 
$B_{{\goth g}}$ of $\tk {{\Bbb C}}{\e Sg}\tk {{\Bbb C}}{\e Sg}{\goth g}$
which is free of rank equal to the dimension $\k g{}$ of the borel
subalgebras of ${\goth g}$. Moreover, $B_{{\goth g}}$ is contained in the
space of cycles of $C_{\bullet}({\goth g})$. The complex
$E_{\bullet}({\goth g})$ is the ideal of $C_{\bullet}({\goth g})$ generated
by $\ex {\k g{}}{B_{{\goth g}}}$. We denote by ${\cal N}_{{\goth g}}$ the set of elements
in ${\goth g}\times {\goth g}$ whose components generate a subsbspace contained in the 
nilpotent cone of ${\goth g}$ and we say that ${\goth g}$ has property $({\bf N})$ if
the codimension of ${\cal N}_{{\goth g}}$ in ${\goth g}\times {\goth g}$ is strictly
bigger than $\k g{}-\j g{}$. Let $I_{{\goth g}}$ be the ideal in 
$\tk {{\Bbb C}}{\e Sg}\e Sg$ generated by the functions 
$(x,y) \mapsto \dv v{[x,y]}$ where $v$ is in ${\goth g}$ and
$\sqrt {I_{{\goth g}}}$ its radical. The main result is the theorem:
\vskip .4em
{\it Let us suppose that for any semi-simple element in ${\goth g}$, the simple factors
of its centralizer in ${\goth g}$ have the property $({\bf N})$. Then the complex 
$E_{\bullet}({\goth g})$ has no homology in degree different 
from $\k g{}$ and $\sqrt {I_{{\goth g}}} \ex {\k g{}}{B_{{\goth g}}}$ is 
contained in the space of boundaries of degree $\k g{}$ of 
$E_{\bullet}({\goth g})$. In particular, $I_{{\goth g}}$ is a prime ideal 
whose set of zeros in ${\goth g}\times {\goth g}$ is the commuting variety of 
${\goth g}$.}

\end{abstract}

\catcode`\@=11
\catcode`\;=\active\def;{\relax\ifhmode\ifdim\lastskip>\z@
\unskip\fi\kern.2em\fi\string;}
\catcode`\:=\active\def:{\relax\ifhmode\ifdim\lastskip>\z@\unsskip\fi
\penalty\@M\ \fi\string:}
\catcode`\!=\active\def!{\relax\ifhmode\ifdim\lastskip>\z@
\unskip\fi\kern.2em\fi\string!}
\catcode`\?=\active\def?{\relax\ifhmode\ifdim\lastskip>\z@
\unskip\fi\kern.2em\fi\string?}\frenchspacing\catcode`\@=12

\section{Introduction.} 
Soient ${\goth g}$ une alg\`ebre de Lie r\'eductive complexe de dimension
finie, $\e Sg$ et $\ex {}{{\goth g}}$ les alg\`ebres sym\'etrique et
ext\'erieure de ${\goth g}$. La repr\'esentation coadjointe de ${\goth g}$
s'identifie \`a sa repr\'esentation adjointe au moyen d'une forme
bilin\'eaire sym\'etrique invariante, non d\'eg\'en\'er\'ee, $\dv ..$ qui 
prolonge la forme de Killing de $[{\goth g},{\goth g}]$. On appelle vari\'et\'e
commutante de ${\goth g}$ la sous-vari\'et\'e ${\goth C}_{{\goth g}}$ des
points $(x,y)$ de ${\goth g}\times {\goth g}$ tels que $[x,y]$ soit nul. La
vari\'et\'e ${\goth C}_{{\goth g}}$ est alors la vari\'et\'e des z\'eros de
l'id\'eal $I_{{\goth g}}$ de $\tk {{\Bbb C}}{\e Sg}\e Sg$ engendr\'e par les
fonctions $(x,y) \mapsto \dv {v}{[x,y]}$ o\`u $v$ est dans ${\goth g}$.
D'apr\`es \cite{Ric}, la vari\'et\'e ${\goth C}_{{\goth g}}$ est
irr\'eductible. La question de savoir si $I_{{\goth g}}$ est un id\'eal
premier s'\'etait pos\'ee depuis plusieurs ann\'ees. Elle est la principale
motivation de ce m\'emoire. D'apr\`es le r\'esultat de R. W. Richardson, il
s'agit en fait de savoir si l'id\'eal $I_{{\goth g}}$ est radiciel. Une
r\'eponse partielle avait \'et\'e donn\'ee par J. Dixmier dans \cite{Dix}. Il
y d\'emontre qu'un champ de vecteurs polynomial sur ${\goth g}$, tangent
aux orbites adjointes, est du type
$x\mapsto [x,\varphi (x)]$ o\`u $\varphi $ est une application polynomiale
de ${\goth g}$ dans ${\goth g}$. Cela revient \`a dire que l'\'egalit\'e de
$I_{{\goth g}}$ et de son radical $\sqrt {I_{{\goth g}}}$ est satisfaite en
degr\'e $1$. Tr\`es r\'ecemment, W. L. Gan et V. Ginzburg ont montr\'e dans
\cite{Gi} que le sous-espace des \'el\'ements ${\goth g}$-invariants de
$I_{{\goth g}}$ est un id\'eal radiciel de l'alg\`ebre des \'el\'ements
${\goth g}$-invariants de $\e Sg$ pour ${\goth g}$ simple de type ${\rm A}$. En 
s'inspirant de la m\'ethode utilis\'ee par J. Dixmier dans \cite{Dix}, on 
introduit dans \cite{Ch} le complexe canonique $C_{\bullet}({\goth g})$ de 
l'alg\`ebre de Lie ${\goth g}$ et on montre qu'il n'a pas d'homologie en 
degr\'e strictement sup\'erieur au rang $\j g{}$ de ${\goth g}$. L'espace 
sous-jacent \`a $C_{\bullet}({\goth g})$ est l'alg\`ebre 
$\tk {{\Bbb C}}{\e Sg}\tk {{\Bbb C}}{\e Sg}\ex {}{{\goth g}}$
et sa diff\'erentielle est la d\'erivation 
$\tk {{\Bbb C}}{\e Sg}\e Sg$-lin\'eaire qui \`a l'\'el\'ement $v$ de 
${\goth g}$ associe la fonction $(x,y)\mapsto \dv v{[x,y]}$ sur ${\goth g}$. 

On appelle {\it sous-module caract\'eristique} pour ${\goth g}$ le sous-module
$B_{{\goth g}}$ des \'el\'ements $\varphi $ de 
$\tk {{\Bbb C}}{\e Sg}\tk {{\Bbb C}}{\e Sg}{\goth g}$ qui ont la propri\'et\'e
suivante: pour tout $(x,y)$ dans un ouvert non vide de 
${\goth g}\times {\goth g}$, $\varphi (x,y)$ appartient \`a la somme des 
centralisateurs dans ${\goth g}$ des \'el\'ements non nuls du sous-espace engendr\'e 
par $x$ et $y$. Le module $B_{{\goth g}}$ a trois propri\'et\'es remarquables:
\begin{list}{}{}
\item 1) $B_{{\goth g}}$ est un module libre de rang $\k g{}$ \'egal \`a la
dimension des sous-alg\`ebres de Borel de ${\goth g}$,
\item 2) $B_{{\goth g}}$ est contenu dans l'espace des cycles du complexe
$C_{\bullet}({\goth g})$,
\item 3) pour toute sous-alg\`ebre parabolique ${\goth p}$ de ${\goth g}$ et
pour tout $\varphi $ dans $B_{{\goth g}}$, ${\goth p}$ contient $\varphi (x,y)$
pour tout $(x,y)$ dans ${\goth p}\times {\goth p}$.
\end{list}
Ce module avait \'et\'e \'etudi\'e en \cite{De} mais introduit de fa\c con diff\'erente.
En particulier, on doit \`a M. Ra\"{\i}s le fait que $B_{{\goth g}}$ est libre de rang 
$\k g{}$. Il r\'esulte de la propri\'et\'e (2) que l'id\'eal $E_{\bullet}({\goth g})$
de $C_{\bullet}({\goth g})$, engendr\'e par $\ex {\k g{}}{B_{{\goth g}}}$,
est un sous-complexe de $C_{\bullet}({\goth g})$. On l'appelle
{\it complexe canonique de deuxi\`eme esp\`ece} de l'alg\`ebre de Lie 
${\goth g}$. On note $X_{{\goth g}}$ l'ensemble des \'el\'ements $(x,y)$ de
${\goth g}\times {\goth g}$ tels que l'image de $B_{{\goth g}}$ par l'application
$\varphi \mapsto \varphi (x,y)$ soit de dimension strictement inf\'erieure \`a
$\k g{}$. 

\begin{Defd}
Soit ${\cal N}_{{\goth g}}$ l'ensemble des \'el\'ements $(x,y)$ de 
${\goth g}\times {\goth g}$ tels que le sous-espace engendr\'e par $x$ et $y$ soit
contenu dans le c\^one nilpotent de ${\goth g}$. On dira que l'alg\`ebre de Lie 
${\goth g}$ a la propri\'et\'e $({\bf N})$ si la codimension de ${\cal N}_{{\goth g}}$
dans ${\goth g}\times {\goth g}$ est strictement sup\'erieure \`a $\k g{}-\j g{}$.
\end{Defd}

Le sous-ensemble ${\cal N}_{{\goth g}}$ de ${\goth g}\times {\goth g}$ est la vari\'et\'e
des z\'eros communs \`a \sloppy \hbox{$\k g{}+\j g{}$} fonctions polynomiales sur 
${\goth g}\times {\goth g}$; donc la codimension de chacune de ses composantes 
irr\'eductibles est inf\'erieure \`a $\k g{}+\j g{}$. Le r\'esultat principal de ce 
m\'emoire est le th\'eor\`eme:

\begin{Thd}
On suppose que les facteurs simples du centralisateur dans ${\goth g}$ de tout \'el\'ement
semi-simple de ${\goth g}$ ont la propri\'et\'e $({\bf N})$. Soient $J_{{\goth g}}$ le 
radical de $I_{{\goth g}}$ et $\overline{E}_{\bullet}({\goth g})$ le sous-complexe de 
$E_{\bullet}({\goth g})$ tel que $E_{j}({\goth g})$ soit \'egal \`a
$\overline{E}_{j}({\goth g})$ pour $j$ diff\'erent de $\k g{}$ et tel que
$\overline{E}_{\k g{}}({\goth g})$ soit \'egal \`a 
$J_{{\goth g}}\ex {\k g{}}{B_{{\goth g}}}$.

{\rm i)} Le complexe $\overline{E}_{\bullet}({\goth g})$ est acyclique.

{\rm ii)} L'id\'eal $I_{{\goth g}}$ est premier.

{\rm iii)} Le complexe $E_{\bullet}({\goth g})$ n'a pas d'homologie en degr\'e
diff\'erent de $\k g{}$ et son homologie est isomorphe \`a l'alg\`ebre
des fonctions r\'eguli\`eres sur la vari\'et\'e commutante de ${\goth g}$.
\end{Thd}

Puisque la vari\'et\'e commutante ${\goth C}_{{\goth g}}$ de ${\goth g}$ est
irr\'eductible, les assertions (ii) et (iii) r\'esultent de l'assertion (i) car
$I_{{\goth g}}\ex {\k g{}}{B_{{\goth g}}}$ est l'espace des bords de degr\'e
$\k g{}$ de $E_{\bullet}({\goth g})$ et de 
$\overline{E}_{\bullet}({\goth g})$. En utilisant les r\'esultats de 
cohomologie \`a support qui sont rappel\'es dans \cite{Ch}, il r\'esulte
du th\'eor\`eme suivant:

\begin{Thd}
On suppose que les facteurs simples de ${\goth g}$ ont la propri\'et\'e $({\bf N})$. 
Alors pour tout entier naturel $k$, la dimension projective de \sloppy \hbox{ 
$\ex k{{\goth g}}\wedge \ex {\k g{}}{B_{{\goth g}}}$} est inf\'erieure \`a
$k$.
\end{Thd}
 
\noindent que le complexe $E_{\bullet}({\goth g})$ n'a pas d'homologie en 
degr\'e diff\'erent de $\k g{}$. On pr\'ecise \`a ce propos que cette 
argumentation g\'en\'eralise l'argumentation de J. Dixmier dans \cite{Dix}. 
Un raisonnement par r\'ecurrence sur la dimension de ${\goth g}$ montre que 
la projection sur ${\goth g}$ du support dans ${\goth g}\times {\goth g}$ du 
quotient $J_{{\goth g}}/I_{{\goth g}}$ ne contient pas d'\'el\'ement 
semi-simple non central. Il en r\'esulte que la dimension de ce support est 
inf\'erieure \`a $\dim {\goth g}$; donc par les arguments ci-dessus, 
$J_{{\goth g}}$ est \'egal \`a $I_{{\goth g}}$. Le deuxi\`eme th\'eor\`eme
est une cons\'equence simple du cas des alg\`ebres de Lie simples. Pour ${\goth g}$
simple, on introduit un complexe de cohomologie 
$D_{k}^{\bullet}({\goth g},B_{{\goth g}})$ dont le terme de plus haut degr\'e est le
module $\ex k{{\goth g}}\wedge \ex {\k g{}}{B_{{\goth g}}}$. Un raisonnement par 
r\'ecurrence sur $k$ montre que l'exactitude de ces complexes pour 
$k=0,\ldots,\k g{}-\j g{}$ donne le th\'eor\`eme ci-dessus. Il est facile de voir
que $X_{{\goth g}}$ contient le support de la cohomologie de ces complexes. En outre, 
d'apr\`es la propri\'et\'e (3) pour le module $B_{{\goth g}}$, le support de la 
cohomologie de ces complexes ne rencontre pas une partie $\Sigma _{{\goth g}}$ de
${\goth g}\times {\goth g}$ d\'efinie au moyen des sous-alg\`ebres paraboliques de
${\goth g}$ dont les facteurs r\'eductifs ont une alg\`ebre de Lie d\'eriv\'ee simple. 

\begin{Defd}
Pour ${\goth g}$ simple, on dira que ${\goth g}$ a la propri\'et\'e $({\bf D})$ 
si toute partie ferm\'ee irr\'eductible de $X_{{\goth g}}$, invariante pour les actions 
de $G$ et de ${\rm GL}_{2}({\Bbb C})$, de codimension inf\'erieure \`a 
$\dim {\goth g}-4$, qui ne rencontre pas $\Sigma _{{\goth g}}$, est contenue dans la 
r\'eunion de ${\cal N}_{{\goth g}}$ et de la sous-vari\'et\'e ${\cal X}_{{\goth g}}$ 
des \'el\'ements $(x,y)$ de ${\goth g}\times {\goth g}$ tels que $x$ et $y$ appartiennent
\`a une m\^eme sous-alg\`ebre de Borel de ${\goth g}$.
\end{Defd}
 
Par d\'efinition, la vari\'et\'e d'incidence de ${\goth g}$ est la sous-vari\'et\'e des
\'el\'ements $(u,x,y)$ o\`u $u$ est une sous-alg\`ebre de Borel de ${\goth g}$
et o\`u $(x,y)$ est un \'el\'ement de $u\times u$. La vari\'et\'e
d'incidence est introduite dans \cite{Ka}. On rappelle que dans \cite{Ka}
est donn\'ee une d\'emonstration simple de l'irr\'eductibilit\'e de la 
vari\'et\'e commutante de ${\goth g}$. Par des arguments de cohomologie \`a support,
l'exactitude des complexes $D_{k}^{\bullet}({\goth g})$ r\'esulte des propri\'et\'es
$({\bf D})$ et $({\bf N})$ pour ${\goth g}$. En remarquant que $X_{{\goth g}}$
ne contient pas d'\'el\'ement $(x,y)$ de ${\cal X}_{{\goth g}}$ tel que $x$ et $y$ soient
des \'el\'ements r\'eguliers de ${\goth g}$, respectivement nilpotent et 
semi-simple, on montre la propri\'et\'e $({\bf D})$.

Dans ce m\'emoire, le corps de base est le corps des nombres complexes,
les espaces vectoriels et les alg\`ebres de Lie consid\'er\'es sont de
dimension finie. On d\'esigne par ${\goth g}$ une alg\`ebre de Lie
r\'eductive, par $G$ son groupe adjoint, par $\j g{}$ son rang, par
$\k g{}$ la dimension des sous-alg\`ebres de Borel de ${\goth g}$ et par
$\dv ..$ une forme bilin\'eaire sym\'etrique, non d\'eg\'en\'er\'ee 
${\goth g}$-invariante sur ${\goth g}$ qui prolonge la forme de Killing de
${\goth g}$. On identifie ${\goth g}$ \`a son dual au moyen de $\dv ..$.
Si $V$ est un espace vectoriel, $V^{*}$ d\'esigne son dual.
Si ${\goth a}$ est une alg\`ebre de Lie et si ${\goth m}$ est un
sous-espace de ${\goth g}$, on note ${\goth m}(x)$ le sous-espace des
\'el\'ements de ${\goth m}$ qui centralisent l'\'el\'ement $x$ de
${\goth a}$ et ${\goth m}(x')$ le sous-espace des \'el\'ements de 
${\goth m}$ qui stabilisent la forme lin\'eaire $x'$ sur ${\goth a}$ pour
l'action coadjointe de ${\goth a}$ dans ${\goth a}^{*}$. On utilise la
topologie de Zariski sur les vari\'et\'es alg\'ebriques consid\'er\'ees. 
Si $X$ est une vari\'et\'e alg\'ebrique, $\an X{}$ d\'esigne le faisceau
structural de $X$ et $\an Xx$ l'anneau local au point $x$ de $X$. Sauf
mention contraire, un point de $X$ est un point ferm\'e. On appelle 
{\it grand ouvert} de $X$, tout ouvert dont le compl\'ementaire est de
codimension sup\'erieure \`a $2$. Conform\'ement \`a l'usage, pour tout ouvert
$Y$ de $X$ et pour tout faisceau ${\cal F}$ sur $X$, $\Gamma (Y,{\cal F})$ 
d\'esigne l'espace des sections de ${\cal F}$ au dessus de $Y$. La suite de ce
m\'emoire se divise en $10$ sections:
\begin{list}{}{}
\item 2) dimension projective et cohomologie, 
\item 3) sur les id\'eaux,
\item 4) exemples de complexes,
\item 5) sous-module caract\'eristique,
\item 6) complexes canoniques d'une alg\`ebre de Lie,
\item 7) au voisinage d'un \'el\'ement semi-simple,
\item 8) sur certaines sous-vari\'et\'es de ${\goth g}\times {\goth g}$,
\item 9) d\'efinition de la propri\'et\'e $({\bf D})$,
\item 10) th\'eor\`eme d'exactitude et id\'eaux premiers,
\item 11) la propri\'et\'e $({\bf D})$ pour les alg\`ebres de Lie simples.
\end{list}
\section{Dimension projective et cohomologie.} \label{p}
On rappelle dans cette section quelques r\'esultats classiques. Soient $X$ une
vari\'et\'e alg\'ebrique affine, de Cohen-Macaulay et ${\Bbb C}[X]$ l'anneau  
des fonctions r\'eguli\`eres sur $X$. On note $U$ un ouvert de $X$, $S$ le  
compl\'ementaire de $U$ dans $X$ et $p$ la codimension de $S$ dans $X$. 
Soient  $P_{\bullet}$ un complexe de ${\Bbb C}[X]$-modules projectifs
de type fini, de longueur finie $l$, et $\varepsilon $ un morphisme 
d'augmentation du complexe $P_{\bullet}$ d'image $R$, d'o\`u un complexe 
augment\'e de modules:
$$ 0 \stackrel{\dd }\longrightarrow P_{l} \stackrel{\dd }\longrightarrow
P_{l-1} \stackrel{\dd }\longrightarrow \cdots
\stackrel{\dd }\longrightarrow P_{0} \stackrel{\varepsilon }
\longrightarrow R \rightarrow 0 \mbox{ .}$$

\begin{prop} \label{pp}
Soient $P$ un module projectif et $R'$ un sous-module de $P$. On suppose que 
les conditions suivantes sont satisfaites:
\begin{list}{}{}
\item {\rm 1)} $p$ est strictement sup\'erieur \`a $l+1$,
\item {\rm 2)} $X$ est normal,
\item {\rm 3)} $S$ contient le support de l'homologie du complexe
augment\'e $P_{\bullet}$,
\item {\rm 4)} $R'$ contient $R$ et $S$ contient le support dans $X$ du
quotient des modules $R'$ et $R$. 
\end{list}
Alors $R'$ est \'egal \`a $R$ et $P_{\bullet}$ est une r\'esolution 
projective de $R$ de longueur $l$.
\end{prop}

\begin{cor} \label{cp}
Soient $d$ un entier naturel et ${\cal C}_{\bullet}$ un complexe d'homologie 
de ${\Bbb C}[X]$-modules de type fini de longueur $l$. On suppose que les
trois conditions suivantes sont r\'ealis\'ees:
\begin{list}{}{}
\item {\rm 1)} $X$ est normal et $S$ contient le support de l'homologie 
du complexe $C_{\bullet}$,
\item {\rm 2)} pour tout entier naturel $i$, $C_{i}$ est un sous-module
d'un module libre,
\item {\rm 3)} pour tout entier strictement positif $i$, la dimension
projective de $C_{i}$ est inf\'erieure \`a $i+d$.
\end{list}
Pour $j$ entier naturel inf\'erieur \`a $l$, on d\'esigne par $Z_{j}$
l'espace des cycles de degr\'e $j$ du complexe $C_{\bullet}$. Si $p$ est
strictement sup\'erieur \`a $2l+d$, alors le complexe $C_{\bullet}$ est
acyclique. En outre, la dimension projective de $Z_{j}$ est inf\'erieure
\`a $2l+d-j-1$.
\end{cor}

Des d\'emon\-strations de ces \'enonc\'es sont donn\'ees dans 
\cite{Ch}(Section 2). 

\begin{lemme}\label{lp}
{\rm i)} Soit
$$0 \rightarrow E_{0} \rightarrow E_{1} \rightarrow E_{2} \rightarrow 0$$
une suite exacte courte de ${\Bbb C}[X]$-modules. On suppose que la dimension
projective des modules $E_{0}$ et $E_{1}$ est inf\'erieure \`a l'entier 
naturel $d$. Alors la dimension projective de $E_{2}$ est inf\'erieure \`a
$d+1$.

{\rm ii)} Soit 
$$ 0 \rightarrow E_{-1}\rightarrow E_{0} \rightarrow \cdots 
\rightarrow E_{l+1} \rightarrow 0 $$
une suite exacte de ${\Bbb C}[X]$-modules. On suppose que $E_{-1}$ est un 
module libre et que pour $i=0,\ldots,l$, la dimension projective de $E_{i}$ 
est inf\'erieure \`a $i$. Alors la dimension projective de $E_{l+1}$ est 
inf\'erieure \`a $l+1$.
\end{lemme}

\begin{proof}
i) Soit $M$ un ${\Bbb C}[X]$-module. Il s'agit de montrer que le groupe
Ext$^{j}(E_{2},M)$ est nul pour $j$ strictement sup\'erieur \`a $d+1$. De la 
suite exacte courte, on d\'eduit la longue suite exacte
$$ \cdots \rightarrow {\rm Ext}^{j}(E_{1},M) \rightarrow 
{\rm Ext}^{j}(E_{0},M) \rightarrow {\rm Ext}^{j+1}(E_{2},M)\rightarrow 
{\rm Ext}^{j+1}(E_{1},M)\rightarrow \cdots  \mbox{ .}$$  
D'apr\`es l'hypoth\`ese, les groupes
${\rm Ext}^{j}(E_{1},M)$, ${\rm Ext}^{j}(E_{0},M)$, ${\rm Ext}^{j+1}(E_{1},M)$
sont nuls pour $j$ strictement sup\'erieur \`a $d$; donc
${\rm Ext}^{j+1}(E_{2},M)$ est nul pour $j$ strictement sup\'erieur \`a
$d$, d'o\`u l'assertion.

ii) Pour $i=0,\ldots,l$, on note $Z_{i}$ l'image de $E_{i-1}$ dans $E_{i}$. On 
montre en raisonnant par r\'ecurrence sur $i$ que la dimension projective de
$Z_{i}$ est inf\'erieure \`a $i$. Pour $i=0$, $Z_{i}$ est isomorphe \`a 
$E_{-1}$. On suppose l'assertion vraie pour $i$. D'apr\`es l'exactitude de la
suite, on a la suite exacte courte
$$ 0 \rightarrow Z_{i} \rightarrow E_{i} \rightarrow Z_{i+1} \rightarrow 0
\mbox{ ;}$$
or par hypoth\`ese, la dimension projective de $E_{i}$ est inf\'erieure \`a
$i$; donc d'apr\`es l'hypoth\`ese de r\'ecurrence et l'assertion (i), la 
dimension projective de $Z_{i+1}$ est inf\'erieure \`a $i+1$. D'apr\`es 
l'exactitude de la suite, on a la suite exacte courte
$$ 0 \rightarrow Z_{l} \rightarrow E_{l} \rightarrow E_{l+1} \rightarrow 0
\mbox{ ;}$$
or par hypoth\`ese, la dimension projective de $E_{l}$ est inf\'erieure \`a
$l$; donc d'apr\`es l'assertion (i), la dimension projective de $E_{l+1}$ est 
inf\'erieure \`a $l+1$.
\end{proof}

\begin{cor} \label{c2p}
Soit ${\cal C}^{\bullet}$ un complexe de cohomologie de ${\Bbb C}[X]$-modules 
de type fini,
$$ 0 \rightarrow C^{-1}\rightarrow C^{0} \rightarrow \cdots 
\rightarrow C^{l} \rightarrow 0 \mbox{ .}$$
On suppose que les trois conditions suivantes sont r\'ealis\'ees:
\begin{list}{}{}
\item {\rm 1)} $X$ est normal et $S$ contient le support de la cohomologie
du complexe $C^{\bullet}$,
\item {\rm 2)} pour tout entier naturel $i$, $C^{i}$ est un sous-module
d'un module libre,
\item {\rm 3)} $C^{-1}$ est un module libre,
\item {\rm 3)} pour tout entier $i$ strictement inf\'erieur \`a $l$, la 
dimension projective de $C^{i}$ est inf\'erieure \`a $i$.
\end{list}
Pour $j$ entier naturel inf\'erieur \`a $l$, on d\'esigne par $Z^{j}$
l'espace des cocycles de degr\'e $j$ du complexe $C^{\bullet}$. Si $p$ est
strictement sup\'erieur \`a $l+1$, alors le complexe $C^{\bullet}$ est
acyclique. En outre, la dimension projective de $Z^{j}$ est inf\'erieure
\`a $j$.
\end{cor}

\begin{proof}
D'apr\`es le lemme \ref{lp}, il s'agit de montrer que le complexe $C^{\bullet}$ est 
acyclique. On montre en raisonnant par r\'ecurrence sur $j$ que pour $j=0,\ldots,l$, 
$Z^{j}$ est l'espace des cobords de degr\'e $j$ du complexe $C^{\bullet}$. 
L'espace $B^{0}$ des cobords de degr\'e $0$ de $C^{\bullet}$ \'etant isomorphe \`a 
$C^{-1}$, $B^{0}$ est \'egal \`a $Z^{0}$ car $C^{0}$ est un module sans torsion et $U$ 
est partout dense dans $X$ comme grand ouvert de la vari\'et\'e de Cohen-Macaulay $X$.
En outre, toute section au dessus de $U$ de la localisation sur $X$ d'un 
${\Bbb C}[X]$-module libre est la restriction \`a $U$ d'un \'el\'ement de ce module car 
$X$ est une vari\'et\'e normale. On suppose que $Z^{j}$ est l'espace des cobords de 
degr\'e $j$ du complexe $C^{\bullet}$. Alors la dimension projective de $Z^{j}$ est 
inf\'erieure \`a $j$; or par hypoth\`ese, la dimension projective de $C^{j}$ est 
inf\'erieure \`a $j$; donc en d\'esignant par $B^{j+1}$ l'espace des cobords de degr\'e 
$j+1$ du complexe $C^{\bullet}$, on a un complexe acyclique
$$ 0 \rightarrow P_{j} \rightarrow P_{j-1}\oplus Q_{j}
\rightarrow \cdots \rightarrow P_{0}\oplus Q_{1} \rightarrow Q_{0}
\rightarrow B^{j+1} \rightarrow 0 \mbox{ ,}$$
o\`u les complexes
$$ 0 \rightarrow P_{j} \rightarrow P_{j-1}
\rightarrow \cdots \rightarrow P_{0} \rightarrow Z^{j} \rightarrow 0 
\mbox{ ,}$$
$$ 0 \rightarrow Q_{j} \rightarrow \cdots \rightarrow Q_{0} \rightarrow 
C^{j} \rightarrow 0 \mbox{ ,}$$
sont des r\'esolutions projectives. Puisque $j+1$ est strictement inf\'erieur
\`a $p$, il r\'esulte de la proposition \ref{pp} que $Z^{j+1}$ est \'egal \`a
$B^{j+1}$ car d'apr\`es la condition (1), $S$ contient le support du quotient 
des modules $Z^{j+1}$ et $B^{j+1}$.
\end{proof}
\section{Sur les id\'eaux.}\label{j}
Soient $X$ une vari\'et\'e alg\'ebrique affine lisse, $Y$ une vari\'et\'e 
alg\'ebrique lisse et $\pi $ un morphisme lisse de $Y$ sur $X$.

\begin{lemme}\label{lj}
Soit $I$ un id\'eal de ${\Bbb C}[X]$ et $J$ son radical. On d\'esigne par
${\cal I}$ et ${\cal J}$ les localisations respectives sur $X$ de $I$ et de
$J$. 

{\rm i)} L'id\'eal $\pi ^{*}({\cal J})$ de $\an Y{}$ est le radical de 
$\pi ^{*}({\cal I})$.

{\rm ii)} Si $\pi ^{*}({\cal I})$ est un id\'eal radiciel, alors $I$ est un 
id\'eal radiciel.
\end{lemme}

\begin{proof}
i) Le morphisme $\pi $ \'etant lisse, les morphismes canoniques de
$\pi ^{*}({\cal I})$ et de $\pi ^{*}({\cal J})$ dans $\an Y{}$ sont des
plongements; donc $\pi ^{*}({\cal I})$ et $\pi ^{*}({\cal J})$ s'identifient \`a 
des id\'eaux de $\an Y{}$. Soient $x$ un point de $X$ et $y$ un point de la 
pr\'eimage de $x$ dans $Y$. On note respectivement ${\cal I}_{x}$, 
${\cal J}_{x}$ les fibres de ${\cal I}$, ${\cal J}$ en $x$ et 
$\pi ^{*}({\cal I})_{y}$, $\pi ^{*}({\cal J})_{y}$ les fibres de 
$\pi ^{*}({\cal I})$, $\pi ^{*}({\cal J})$ en $y$. Soit 
$\poi x1{,\ldots,}{m}{}{}{}$ un syst\`eme de coordonn\'ees de l'anneau local 
$\an Xx$. On identifie l'anneau local $\an Xx$ \`a un sous-anneau de l'anneau 
local $\an Yy$ au moyen du comorphisme de $\pi $. Puisque $\pi $ est un 
morphisme lisse, le syst\`eme de coordonn\'ees $\poi x1{,\ldots,}{m}{}{}{}$ se
compl\`ete en un syst\`eme 
$\poi x1{,\ldots,}{m}{}{}{},\poi y1{,\ldots,}{n}{}{}{}$ de coordonn\'ees de
l'anneau local $\an Yy$. On note ${\goth m}_{x}$ l'id\'eal maximal de
$\an Xx$, $\han Xx$ le compl\'et\'e de $\an Xx$ pour la topologie 
${\goth m}_{x}$-adique, ${\goth m}_{y}$ l'id\'eal maximal de
$\an Yy$, $\han Yy$ le compl\'et\'e de $\an Yy$ pour la topologie 
${\goth m}_{y}$-adique. Puisque $Y$ est lisse en $y$, le morphisme
canonique de l'anneau de s\'eries formelles 
${\Bbb C}[[\poi x1{,\ldots,}{m}{}{}{},\poi y1{,\ldots,}{n}{}{}{}]]$ dans
$\han Yy$ est un isomorphisme. En outre, l'image du sous-anneau 
${\Bbb C}[[\poi x1{,\ldots,}{m}{}{}{}]]$ est \'egale \`a $\han Xx$. 

Puisque $J$ est le radical de $I$, $\pi ^{*}({\cal J})$ est contenu dans le 
radical de $\pi ^{*}({\cal I})$. Vu l'arbitraire du point $x$ de $X$ et du 
point $y$ de la fibre de $\pi $ en $x$, il s'agit de montrer que 
$\pi ^{*}({\cal J})_{y}$ contient le radical de $\pi ^{*}({\cal I})_{y}$. 
Puisque $\pi ^{*}({\cal J})_{y}$ est contenu dans le radical de 
$\pi ^{*}({\cal I})_{y}$, il suffit de montrer que 
$\widehat{\pi ^{*}({\cal J})_{y}}$ contient le radical de
$\widehat{\pi ^{*}({\cal I})_{y}}$ car $\han Yy$ est fid\`element plat sur
$\an Yy$. D'apr\`es ce qui pr\'ec\`ede, on a
$$ \widehat{\pi ^{*}({\cal I})_{y}} = 
\tk {{\Bbb C}[[\poi x1{,\ldots,}{m}{}{}{}]]} 
{{\Bbb C}[[\poi x1{,\ldots,}{m}{}{}{},\poi y1{,\ldots,}{n}{}{}{}]]}
\widehat{{\cal I}_{x}} $$ $$ \mbox{ et }
\widehat{\pi ^{*}({\cal J})_{y}} = 
\tk {{\Bbb C}[[\poi x1{,\ldots,}{m}{}{}{}]]} 
{{\Bbb C}[[\poi x1{,\ldots,}{m}{}{}{},\poi y1{,\ldots,}{n}{}{}{}]]}
\widehat{{\cal J}_{x}} \mbox{ .}$$
Puisque $J$ est le radical de $I$, ${\cal J}_{x}$ est le radical de 
${\cal I}_{x}$. Soient $\poi P1{,\ldots,}{m}{}{}{}$ les id\'eaux premiers 
minimaux de $\an Xx$ qui contiennent ${\cal J}_{x}$. Puisque $\han Xx$ est
une extension plate de $\an Xx$, $\widehat{{\cal J}_{x}}$ est l'intersection
des adh\'erences dans $\han Xx$ des id\'eaux premiers 
$\poi P1{,\ldots,}{m}{}{}{}$; or d'apr\`es 
\cite{Na}(Theorem 36.4 and Theorem 36.5), pour $i=1,\ldots,m$, l'adh\'erence
de $P_{i}$ dans $\han Xx$ est un id\'eal radiciel; donc 
$\widehat{{\cal J}_{x}}$ est le radical de $\widehat{{\cal I}_{x}}$. Soit $a$ 
un \'el\'ement du radical de
$\widehat{\pi ^{*}({\cal I})_{y}}$. Notant $a_{i_{1},\ldots,i_{n}}$ le
coefficient de $y_{1}^{i_{1}}\cdots y_{n}^{i_{n}}$ dans le d\'eveloppement
de $a$ suivant les puissances de $\poi y1{,\ldots,}{n}{}{}{}$, on montre en
raisonnant par r\'ecurrence sur le $n$-uplet $(\poi i1{,\ldots,}{n}{}{}{})$
que $\widehat{{\cal J}_{x}}$ contient $a_{\poi i1{,\ldots,}{n}{}{}{}}$; donc
$\widehat{\pi ^{*}({\cal J})_{y}}$ contient $a$, d'o\`u l'assertion.

ii) On suppose que $\pi ^{*}({\cal I})$ est un id\'eal radiciel. D'apr\`es
l'assertion (i), cela revient \`a dire que les id\'eaux $\pi ^{*}({\cal I})$
et $\pi ^{*}({\cal J})$ sont \'egaux. Il r\'esulte alors de (i) que
$\widehat{{\cal I}_{x}}$ est \'egal \`a $\widehat{{\cal J}_{x}}$ car 
$\han Yy$ est une extension fid\`element plate de $\han Xx$. Par suite,
${\cal I}_{x}$ et ${\cal J}_{x}$ sont \'egaux car $\han Xx$ est une extension
fid\`element plate de $\an Xx$. Vu l'arbitraire de $x$, $I$ est \'egal \`a $J$
car $X$ est affine, d'o\`u l'assertion.
\end{proof}

 \section{Exemples de complexes.} \label{co}
Dans cette section, on donne des exemples de complexes qui seront utilis\'es 
dans la suite de ce m\'emoire. Soient $X$ une vari\'et\'e alg\'ebrique affine,
irr\'eductible et ${\Bbb C}[X]$ l'alg\`ebre des fonctions r\'eguli\`eres sur 
$X$. Dans ce qui suit, on d\'esigne par $V$ un espace vectoriel. On note 
respectivement $\es SV$ et $\ex {}V$ les alg\`ebres sym\'etrique et 
ext\'erieure de $V$. Pour tout entier $i$, $\sy iV$ et $\ex iV$ d\'esignent 
respectivement les sous-espaces de degr\'e $i$ pour les graduations usuelles
des alg\`ebres $\es SV$ et $\ex {}V$. En particulier, pour $i$ strictement 
n\'egatif, $\sy iV$ et $\ex iV$ sont nuls. 

\subsection{} Pour tout entier naturel $d$ inf\'erieur \`a la dimension de $V$,
on d\'esigne par $\ec {Gr}{}{V}{}d$ la grassmannienne de degr\'e $d$ de $V$. 
Pour tout sous-module $L$ de rang $d$ de $\tk {{\Bbb C}}{{\Bbb C}[X]}V$, on 
note $\overline{\ex dL}$ l'ensemble des \'el\'ements de 
$\tk {{\Bbb C}}{{\Bbb C}[X]}\ex dV$ dont le produit par un \'el\'ement non 
nul de ${\Bbb C}[X]$ appartient \`a $\ex dL$. 

\begin{lemme}\label{lco1}
On suppose $X$ lisse et factorielle. Soient $d$ un entier naturel, 
$\alpha $ une application r\'eguli\`ere d'un ouvert non vide $Y$ de $X$
dans $\ec {Gr}{}{V}{}d$ et $L$ l'ensemble des \'el\'ements $\varphi $ de 
$\tk {{\Bbb C}}{{\Bbb C}[X]}V$ tels que $\alpha (x)$ contienne 
$\varphi (x)$ pour tout $x$ dans un ouvert non vide de $U$. 

{\rm i)} Il existe un grand ouvert $\tilde{Y}$ de $X$
contenant $Y$ et une application r\'eguli\`ere $\beta $ de $\tilde{Y}$
dans $\ec {Gr}{}{V}{}{d}$ qui prolonge $\alpha $.

{\rm ii)} L'\'el\'ement $\varphi $ de $\tk {{\Bbb C}}{{\Bbb C}[X]}V$
appartient \`a $L$ si et seulement si $\beta (x)$ contient $\varphi (x)$
pour tout $x$ dans $\tilde{Y}$.

{\rm iii)} Soit ${\cal L}$ la localisation de $L$ sur $X$. Alors
la restriction de ${\cal L}$ \`a $\tilde{Y}$ est un module localement libre
de rang $d$.

{\rm iv)} Le module $\overline{\ex dL}$ est libre de rang $1$.

{\rm v)} Soient $\poi {\varepsilon }1{,\ldots,}{d}{}{}{}$ des \'el\'ements
de $L$. Alors $\poi {\varepsilon }1{,\ldots,}{d}{}{}{}$ est une base du 
${\Bbb C}[X]$-module $L$ si et seulement si 
$\poi x{}{,\ldots,}{}{\varepsilon }{1}{d}$ est une base de $\beta (x)$ pour
tout $x$ dans $\tilde{Y}$.

{\rm vi)} Si $L$ est un module libre, alors les modules $\ex dL$ et
$\overline{\ex dL}$ sont \'egaux.
\end{lemme}

\begin{proof}
i) L'assertion (i) r\'esulte de \cite{Sh}(Ch. VI, Theorem 1).

ii) Soit $\varphi $ dans $\tk {{\Bbb C}}{{\Bbb C}[X]}V$. Si $\alpha (x)$
contient $\varphi (x)$ pour tout $x$ dans un ouvert non vide de $U$, 
$\beta (x)$ contient $\varphi (x)$ pour tout $x$ dans $\tilde{Y}$ car tout
ouvert non vide de $\tilde{Y}$ est partout dense dans $\tilde{Y}$, d'o\`u
l'assertion.

iii) Puisque $\beta $ est une application r\'eguli\`ere, $\tilde{Y}$ est
recouvert par des ouvert affines $Z$ qui ont la propri\'et\'e suivante: il
existe des applications r\'eguli\`eres $\poi {\eta }1{,\ldots,}{d}{}{}{}$ de
$Z$ dans $V$ telles que $\poi x{}{,\ldots,}{}{\eta }{1}{d}$ soit une base de $\beta (x)$ 
pour tout $x$ dans $Z$. Pour $i=1,\ldots,d$, $\eta _{i}$ est une section locale de 
${\cal L}$. D'apr\`es (ii), toute section de ${\cal L}$ au dessus de $Z$ est combinaison 
lin\'eaire \`a coefficients dans ${\Bbb C}[Z]$ des applications 
$\poi {\eta }1{,\ldots,}{d}{}{}{}$; or 
$\poi {\eta }1{,\ldots,}{d}{}{}{}$ sont lin\'eairement
ind\'ependants sur ${\Bbb C}[Z]$; donc la restriction de ${\cal L}$
\`a $Z$ est un $\an Z{}$-module libre de rang $d$ car $Z$ est un ouvert
affine.

iv) D'apr\`es (iii), $L$ est un module de rang $d$; donc 
$\ex dL$ et $\overline{\ex dL}$ sont des modules
de rang $1$. Soit $\varepsilon '$ un \'el\'ement non nul de $\ex dL$.
Puisque ${\Bbb C}[X]$ est un anneau factoriel, $\varepsilon '$ est le
produit d'un \'el\'ement de ${\Bbb C}[X]$ et d'un \'el\'ement de
$\tk {{\Bbb C}}{{\Bbb C}[X]}\ex dV$ dont la vari\'et\'e des
z\'eros est de codimension sup\'erieure \`a $2$ dans $X$. Par
d\'efinition, $\overline{\ex dL}$ contient cet
\'el\'ement. Soit $\varepsilon $ un \'el\'ement de 
$\overline{\ex dL}$ dont la vari\'et\'e des z\'eros dans
$X$ est de codimension sup\'erieure \`a $2$. Pour tout 
$\varphi $ dans $\overline{\ex dL}$, il existe des
\'el\'ements $p$ et $q$ de ${\Bbb C}[X]$ tels que 
$p\varepsilon $ soit \'egal \`a $q\varphi $. Puisque la vari\'et\'e des
z\'eros de $\varepsilon $ est de codimension sup\'erieure \`a $2$, la
vari\'et\'e des z\'eros de $q$ est contenue dans la vari\'et\'e des z\'eros
de $p$; donc $q$ divise $p$ et $\varphi $ appartient au ${\Bbb
C}[X]$-module engendr\'e par $\varepsilon $, d'o\`u l'assertion.

v) On suppose que $\poi {\varepsilon }1{,\ldots,}{d}{}{}{}$ est une base de
$L$. Soient $Z$ un ouvert affine de $\tilde{Y}$ et 
$\poi {\eta }1{,\ldots,}{d}{}{}{}$ des applications r\'eguli\`eres de $Z$
dans $V$ telles que $\poi x{}{,\ldots,}{}{\eta }{1}{d}$ soit une base de
$\beta (x)$ pour tout $x$ dans $Z$. Pour $i=1,\ldots,d$, $\eta _{i}$ est
une section locale de ${\cal L}$; donc $\poi {\eta }1{,\ldots,}{d}{}{}{}$ sont
combinaisons lin\'eaires \`a coefficients dans ${\Bbb C}[Z]$ des
applications $\poi {\varepsilon }1{,\ldots,}{d}{}{}{}$. Il en r\'esulte que
pour tout $x$ dans $Z$, $\poi x{}{,\ldots,}{}{\eta }{1}{d}$ sont
combinaisons lin\'eaires \`a coefficients dans ${\Bbb C}$ de 
$\poi x{}{,\ldots,}{}{\varepsilon }{1}{d}$; donc pour tout $x$ dans
$\tilde{Y}$, $\poi x{}{,\ldots,}{}{\varepsilon }{1}{d}$ est une base de
$\beta (x)$. R\'eciproquement, on suppose que 
$\poi x{}{,\ldots,}{}{\varepsilon }{1}{d}$ est une base de $\beta (x)$ pour
tout $x$ dans $\tilde{Y}$. D'apr\`es l'assertion (ii), pour tout $\varphi $
dans $L$, il existe des fonctions r\'eguli\`eres $\poi a1{,\ldots,}{d}{}{}{}$
sur $\tilde{Y}$ qui satisfont l'\'egalit\'e
$$ \varphi (x) = a_{1}(x) \varepsilon _{1}(x) + \cdots +
a_{d}(x) \varepsilon _{d}(x) \mbox{ ,}$$
pour tout $x$ dans $\tilde{Y}$. Puisque $\tilde{Y}$ est un grand ouvert de
$X$, ces fonctions ont un prolongement r\'egulier \`a 
$X$ car $X$ est une vari\'et\'e lisse; donc  
$\poi {\varepsilon }1{,\ldots,}{d}{}{}{}$ engendrent le module $L$. Puisque
$\poi x{}{,\ldots,}{}{\varepsilon }{1}{d}$ sont lin\'eairement
ind\'ependants pour tout $x$ dans $\tilde{Y}$, 
$\poi {\varepsilon }1{,\ldots,}{d}{}{}{}$ sont lin\'eairement ind\'ependants
sur ${\Bbb C}[X]$; donc $\poi {\varepsilon }1{,\ldots,}{d}{}{}{}$ est une
base de
$L$.

vi) On suppose que $L$ est un module libre. Alors $\ex dL$ est un module
libre de rang $1$. Soient $\varepsilon $ un g\'en\'erateur de
$\overline{\ex dL}$, $\eta $ un g\'en\'erateur de $\ex dL$, $Z$ un ouvert
affine de $\tilde{Y}$ et $\poi {\eta }1{,\ldots,}{d}{}{}{}$ des applications
r\'eguli\`eres de $Z$ dans $V$ telles que $\poi x{}{,\ldots,}{}{\eta }{1}{d}$ soit une 
base de $\beta (x)$ pour tout $x$ dans $Z$. Puisque 
$\poi {\eta }1{\wedge \cdots \wedge}{d}{}{}{}$ est une section au dessus $Z$ de la 
localisation sur $X$ du module $\ex dL$, il est le produit de $\eta $ et d'une fonction 
sur $Z$; donc $\eta $ ne s'annule pas sur $Z$. D'apr\`es (iii), $\tilde{Y}$
est recouvert par des ouverts $Z$ qui ont la propri\'et\'e ci-dessus; donc
$\eta $ ne s'annule pas sur $\tilde{Y}$. Puisque 
$\overline{\ex dL}$ contient $\eta $, $\eta $ est le produit de
$\varepsilon $ et d'un \'el\'ement $p$ de ${\Bbb C}[X]$. Alors la vari\'et\'e 
des z\'eros de $p$ dans $X$ ne rencontre pas $\tilde{Y}$; or
$\tilde{Y}$ est un grand ouvert de $X$; donc $p$ est inversible dans
${\Bbb C}[X]$ et $\ex dL$ contient $\varepsilon $, d'o\`u l'assertion.
\end{proof}

\subsection{} Soient $\lambda $ une application lin\'eaire de $V$ dans
${\Bbb C}[X]$ et $\theta _{\lambda }$ l'application lin\'eaire 
$$\tk {{\Bbb C}}{{\Bbb C}[X]}V \stackrel{\theta _{\lambda }}
\longrightarrow {\Bbb C}[X] \mbox{ , } a\tens v \mapsto a\lambda (v) 
\mbox{ ,}$$ 
o\`u $a$ est dans ${\Bbb C}[X]$ et o\`u $v$ est dans $V$.

\begin{Def} \label{dco2}
On appelle complexe canonique associ\'e \`a $\lambda $ et on le note
$C_{\bullet}(\lambda )$ le complexe d'espace
$\tk {{\Bbb C}}{{\Bbb C}[X]}\ex {}V$ dont la diff\'erentielle est la
${\Bbb C}[X]$-d\'erivation de l'alg\`ebre 
$\tk {{\Bbb C}}{{\Bbb C}[X]}\ex {}V$ qui prolonge $\theta _{\lambda }$. La
graduation naturelle de $\ex {}V$ induit sur $C_{\bullet}(\lambda )$ une
structure de complexe  gradu\'e.

Soient $K_{\lambda }$ et $I_{\lambda }$ le noyau et l'image de 
$\theta _{\lambda }$. Soit $L$ un sous-module de rang $r$ de
$K_{\lambda }$. On appelle espace canonique associ\'e \`a $\lambda $ et \`a $L$,
l'id\'eal $C_{\bullet}(\lambda ,L)$ de $C_{\bullet}(\lambda )$ engendr\'e 
par $\overline{\ex rL}$.
\end{Def}

Pour $\pi $ automorphisme de la vari\'et\'e $X$, on d\'esigne par
$\pi ^{\#}$ l'automorphisme de l'alg\`ebre 
$\tk {{\Bbb C}}{{\Bbb C}[X]}\ex {}V$ qui \`a l'application $\varphi $ de $X$
dans $\ex {}V$ associe $\varphi \rond \pi $ et $\lambda _{\pi }$ l'application
de $V$ dans ${\Bbb C}[X]$ qui \`a l'\'el\'ement $v$ associe la fonction
$\lambda (v)\rond \pi $ sur $X$.  

\begin{lemme}\label{lco2}
Soient $\pi $ un automorphisme de $X$, $\dd $ la diff\'erentielle du
complexe $C_{\bullet}(\lambda )$ et $L$ un sous-module de rang $r$ 
de $K_{\lambda }$.

{\rm i)} Le module $K_{\lambda _{\pi }}$ et l'id\'eal 
$I_{\lambda _{\pi }}$ sont les images respectives de $K_{\lambda }$ et de
$I_{\lambda }$ par $\pi ^{\#}$.

{\rm ii)} La diff\'erentielle du complexe $C_{\bullet}(\lambda _{\pi })$
est l'application $\pi ^{\#}\rond \dd \rond (\pi ^{\#})^{-1}$.

{\rm iii)} L'espace $C_{\bullet}(\lambda ,L)$ est un sous-complexe
de $C_{\bullet}(\lambda )$. En outre, l'image par $\pi ^{\#}$ du complexe 
$C_{\bullet}(\lambda ,L)$ est le sous-complexe
$C_{\bullet}(\lambda _{\pi },\pi ^{\#}(L))$ de 
$C_{\bullet}(\lambda _{\pi })$.
\end{lemme}

\begin{proof}
Les d\'emonstrations des assertions (i) et (ii) sont donn\'ees dans 
\cite{Ch}(2.1). 

iii) Puisque $L$ est contenu dans $K_{\lambda }$, $\ex r{L}$ est un
sous-espace du noyau de $\dd $; or $\tk {{\Bbb C}}{{\Bbb C}[X]}\ex rV$
est un module sans torsion; donc $\overline{\ex rL}$ est un
sous-espace du noyau de $\dd $. Il en r\'esulte que
$C_{\bullet}(\lambda ,L)$ est stable par $\dd $ car $\dd $ est une
d\'erivation de l'alg\`ebre $\tk {{\Bbb C}}{{\Bbb C}[X]}\ex {}V$. Puisque 
$\ex r{\pi ^{\#}(L)}$ est l'image de $\ex rL$ par $\pi ^{\#}$, 
$\overline{\ex r{\pi ^{\#}(L)}}$ et
$C_{\bullet}(\lambda _{\pi },\pi ^{\#}(L))$ sont les images
respectives de $\overline{\ex r{L}}$ et de
$C_{\bullet}(\lambda ,L)$ par $\pi ^{\#}$ car $\pi ^{\#}$ est un
automorphisme de l'alg\`ebre $\tk {{\Bbb C}}{{\Bbb C}[X]}\ex {}V$ et 
$\pi ^{\#}(L)$ est un module de rang $r$. 
\end{proof}

\subsection{} Soient $W$ un espace vectoriel et $\tau $ une application r\'eguli\`ere de 
$X$ dans l'espace des applications lin\'eaires de $V$ dans $W$. 

\begin{Def}\label{dco3}
On appelle complexe canonique associ\'e \`a $\tau $ et on le note 
$C_{\bullet}(\tau )$ le complexe $C_{\bullet}(\lambda )$, d\'efini en
\ref{dco2}, o\`u $\lambda $ est l'application lin\'eaire de $V$ dans 
${\Bbb C}[X\times W^{*}]$ qui \`a $v$ associe la fonction
$(x,w') \mapsto \dv {w'}{\tau (x)(v)}$. 

Si $L$ est un sous-module de $K_{\lambda }$, on note $C_{\bullet}(\tau ,L)$ 
le complexe $C_{\bullet}(\lambda ,L)$ d\'efini en \ref{dco2}.
\end{Def}

On rappelle que $W^{*}$ est le dual de $W$ et que l'alg\`ebre 
${\Bbb C}[X\times W^{*}]$ est canoniquement isomorphe \`a l'alg\`ebre 
$\tk {{\Bbb C}}{{\Bbb C}[X]}\es SW$. 

\begin{lemme}\label{lco3}
On note ${\goth C}_{\tau }$ la vari\'et\'e des z\'eros de $I_{\lambda }$.

{\rm i)} La vari\'et\'e ${\goth C}_{\tau }$ est l'ensemble des \'el\'ements
$(x,w')$ de $X\times W^{*}$ qui satisfont la condition  suivante: $w'$ est
orthogonal \`a l'image de $\tau (x)$. 

{\rm ii)} Le radical de $I_{\lambda }$ est l'ensemble des \'el\'ements
$\varphi $ de $\tk {{\Bbb C}}{{\Bbb C}[X]}\es SW$ qui satisfont la
condition suivante: pour tout $x$ dans $X$, $\varphi (x)$ appartient \`a
l'id\'eal de $\es SW$ engendr\'e par l'image de $\tau (x)$. 

{\rm iii)} Pour tout sous-module $L$ de $K_{\lambda }$, le support de
l'homologie du complexe $C_{\bullet}(\tau ,L)$ est contenu dans 
${\goth C}_{\tau }$.
\end{lemme}

\begin{proof}
Les d\'emonstrations des assertions (i) et (ii) sont donn\'ees dans 
\cite{Ch}(3.2).

iii) On note respectivement ${\cal C}_{\bullet}(\tau )$ et
${\cal C}_{\bullet}(\tau ,L)$ les localisations sur $X\times W^{*}$ des
complexes $C_{\bullet}(\tau )$ et $C_{\bullet}(\tau ,L)$. Soit
$(x,w')$ un point de $X\times W^{*}$ qui n'est pas dans 
${\goth C}_{\tau }$. Puisque la forme lin\'eaire $w'\rond \tau (x)$ n'est
pas nulle, il existe un ouvert affine $U$ de $X\times W^{*}$ qui
contient $(x,w')$ et des applications r\'eguli\`eres 
$\poi {\varphi }1{,\ldots,}{n}{}{}{}$ de $U$ dans $V$ qui satisfont
les conditions suivantes: pour tout $(z,z')$ dans $U$, 
$\dv {z'}{\tau (z)(\varphi _{n}(z,z'))}$ est \'egal
\`a $1$ et $\poi {z,z'}{}{,\ldots,}{}{\varphi }{1}{n-1}$ est une base du
noyau de $z'\rond \tau (z)$. On d\'esigne respectivement par $C_{\bullet}$
et $C_{\bullet}(L)$ les espaces des sections au dessus de $U$ des complexes
${\cal C}_{\bullet}(\tau )$ et ${\cal C}_{\bullet}(\tau ,L)$. Puisque le
bord de $\varphi_{n}$ est \'egal \`a $1$, $c$ est le bord de 
$\varphi _{n}\wedge c$ pour tout cycle $c$ de $C_{\bullet}$; or 
$C_{\bullet}(L)$ contient $\varphi _{n}\wedge c$ s'il contient $c$; donc 
$C_{\bullet}(L)$ est acyclique. Vu l'arbitraire de $(x,w')$,
${\goth C}_{\tau }$ contient le support de l'homologie du complexe 
$C_{\bullet}(\tau ,L)$.
\end{proof}
\subsection{} L'injection canonique de $V$ dans $\ex {}V$ s'\'etend de
mani\`ere unique en une d\'erivation de l'alg\`ebre 
$\tk {{\Bbb C}}{\es SV}\ex {}V$ qui est nulle sur la sous-alg\`ebre 
$\tk {{\Bbb C}}1\ex {}V$. On d\'esigne par $D^{\bullet}(V)$ le complexe
d'espace $\tk {{\Bbb C}}{\es SV}\ex {}V$ dont la diff\'erentielle est la
d\'erivation ainsi d\'efinie. La graduation naturelle de $\ex {}V$ induit
sur $D^{\bullet}(V)$ une structure de complexe gradu\'e. 

\begin{Def}\label{dco4}
Soit $L$ un sous-module de rang $r$ de $\tk {{\Bbb C}}{{\Bbb C}[X]}V$. 
Pour $k$ entier naturel, on pose:
$$ D_{k}^{j}(V,L) = \left \{ \begin{array}{ccc}
\tk {{\Bbb C}[X]}{\sy kL}\ex rL & \mbox{ si } & j=r-1 \\
\tk {{\Bbb C}}{\sy {k-j}V}\ex rL\wedge \ex jV & \mbox{ si } & j\geq r \\
0 & \mbox{ si } & j < r-1 \end{array} \right. \mbox{ ,}$$ 
et on note $D^{\bullet}_{k}(V,L)$ la somme directe des sous-espaces 
$D_{k}^{j}(V,L)$ o\`u \sloppy \hbox{$j=0,1,\ldots$} et $\dd _{k}$ 
l'endomorphisme ${\Bbb C}[X]$-lin\'eaire de $D_{k}^{\bullet}(V,L)$ qui 
prolonge l'injection canonique de $\tk {{\Bbb C}[X]}{\sy kL}\ex {r}L$ dans 
$\tk {{\Bbb C}}{\sy kV}\ex {r}L$ et la restriction \`a \sloppy
\hbox{$\tk {{\Bbb C}}{\sy {k-i}V}\ex {r}L\wedge \ex iV$} 
de la diff\'erentielle de $\tk {{\Bbb C}}{{\Bbb C}[X]}D^{\bullet}(V)$ pour $i$ 
entier naturel.
\end{Def}

\begin{lemme} \label{lco4}
Soit $k$ un entier naturel.  
 
{\rm i)} L'endomorphisme $\dd _{k}$ est une structure de complexe de 
cohomologie sur $D^{\bullet}_{k}(V,L)$. 

{\rm ii)} La cohomologie du complexe $D^{\bullet}(V)$ est \'egale \`a 
${\Bbb C}$. 

{\rm iii)} Pour tout sous-espace $E$ de $V$, $D^{\bullet}_{k}(V,E)$ est un 
complexe acyclique.

{\rm iv)} La cohomologie du complexe $D^{\bullet}_{k}(V,L)$ est un 
${\Bbb C}[X]$-module gradu\'e de type fini.
\end{lemme}

\begin{proof}
i) Pour tout couple $(i,j)$ d'entiers naturels, l'image de \sloppy
$\tk {{\Bbb C}}{\sy {i+1}V}\ex rL\wedge \ex jV$ par la diff\'erentielle du 
complexe $\tk {{\Bbb C}}{{\Bbb C}[X]}D^{\bullet}(V)$ est contenue dans 
$\tk {{\Bbb C}}{\sy iV}\ex rL\wedge \ex {j+1}V$. En outre, 
$\tk {{\Bbb C}[X]}{\sy kL}\ex rL$ est contenu dans le noyau de la restriction de 
$\dd _{k}$ \`a $\tk {{\Bbb C}}{\sy kV}\ex rL$; donc $\dd _{k}$ est une structure de
complexe.

ii) On montre l'assertion en raisonnant par r\'ecurrence sur la
dimension de $V$. Pour $V$ nul, le complexe $D^{\bullet}(V)$ est le
complexe d'espace ${\Bbb C}$ et de diff\'erentielle nulle. On suppose
l'assertion vraie pour tout espace vectoriel de dimension strictement
inf\'erieure \`a celle de $V$. Soient $\dd $ la diff\'erentielle du complexe
$D^{\bullet}(V)$, $W$ un sous-espace de codimension $1$ de $V$ et $v$ un
\'el\'ement de $V$ qui n'est pas dans $W$. Soit $a$ un cocycle homog\`ene de
$D^{\bullet}(V)$. Alors $a$ a un unique d\'eveloppement
$$ a = v^{m}a_{m} + \cdots + a_{0} \mbox{ ,}$$
o\`u pour tout entier naturel $i$, $a_{i}$ est un \'el\'ement de
$\tk {{\Bbb C}}{\es SW}\ex {}V$. Si $a$ est de degr\'e nul, alors on a
$$ mv^{m-1}a_{m}\tens v + \cdots + a_{1}\tens v =0 \mbox{ ;}$$
donc dans ce cas, ${\Bbb C}$ contient $a$. On suppose $a$ de degr\'e 
strictement positif $d$. Alors pour tout entier naturel $i$, on a
$$a_{i} = a'_{i} + a''_{i}\wedge v \mbox{ ,}$$
o\`u $a'_{i}$ et $a''_{i}$ sont respectivement dans 
$\tk {{\Bbb C}}{\es SV}\ex dW$ et $\tk {{\Bbb C}}{\es SV}\ex {d-1}W$. De 
l'\'egalit\'e
$$ 0 = \sum_{i=0}^{m} v^{i}\dd a'_{i} + \sum_{i=1}^{m} iv^{i-1}
a'_{i}\wedge v + \sum_{i=0}^{m} v^{i}\dd a''_{i}\wedge v \mbox{ ,}$$
on d\'eduit que $a'_{0},\ldots,a'_{m}$ sont des cocycles; donc d'apr\`es
l'hypoth\`ese de r\'ecurrence, pour $i=0,\ldots,m$, $a'_{i}$ est le cobord 
d'un \'el\'ement $b_{i}$ de $\tk {{\Bbb C}}{\es SW}\ex {d-1}W$. Par suite,
il vient
$$ a - \dd(\sum_{i=0}^{m} v^{i}b_{i}) = 
v^{m}a''_{m}\wedge v + 
\sum_{i=0}^{m-1} v^{i}((-1)^{d-1}(i+1)b_{i+1}+a''_{i})\wedge v \mbox{ ;}$$
donc $a''_{m}$ et $(-1)^{d-1}(i+1)b_{i+1}+a''_{i}$ sont des cocycles de 
degr\'e $d-1$ pour \sloppy \hbox{$i=0,\ldots,m-1$}. Pour $d=1$, il vient
$$ a = \dd(\sum_{i=0}^{m} v^{i}b_{i} + 
\frac{1}{m+1}v^{m+1}a''_{m} + \sum_{i=0}^{m-1} 
\frac{1}{i+1}v^{i+1}((i+1)b_{i+1}+a''_{i})) \mbox{ ;}$$ 
donc on peut supposer $d$ sup\'erieur \`a $2$. Il r\'esulte alors de 
l'hypoth\`ese de r\'ecurrence que 
$a''_{m}$ est le cobord d'un \'el\'ement $c_{m}$ de 
$\tk {{\Bbb C}}{\es SW}\ex {d-2}W$. De m\^eme, pour $i=0,\ldots,m-1$,
$(-1)^{d-1}(i+1)b_{i+1}+a''_{i}$ est le cobord d'un \'el\'ement $c_{i}$ de 
$\tk {{\Bbb C}}{\es SW}\ex {d-2}V$; donc il vient
$$ a =  \dd(\sum_{i=0}^{m} v^{i}b_{i} + \sum_{i=0}^{m} v^{i}c_{i}\wedge v) 
\mbox{ ,}$$
d'o\`u l'assertion.

iii) On montre en raisonnant par r\'ecurrence sur la codimension $\delta $ de 
$E$ dans $V$ que pour tout entier naturel $k$, le complexe 
$D_{k}^{\bullet}(V,E)$ est acyclique. Le complexe $D_{k}^{\bullet}(V,V)$ est
acyclique par d\'efinition. On suppose l'assertion vraie pour $\delta -1$ et 
$\delta $ strictement positif. Soient $W$ un sous-espace de codimension $1$ 
de $V$ qui contient $E$ et $v$ un \'el\'ement de $V$ qui n'est pas dans $W$. 
Soit $a$ un cocycle de degr\'e $d$ du complexe $D_{k}^{\bullet}(V,E)$. Alors 
$a$ a un unique d\'eveloppement
$$ a = \sum_{i=0}^{m} v^{i}(a'_{i}+a''_{i}\wedge v) \mbox{ ,}$$
o\`u pour $i=0,\ldots,m$, $a'_{i}$ et $a''_{i}$ sont respectivement dans
$D_{k-i}^{d}(W,E)$ et $D_{k-i-1}^{d-1}(W,E)$; or d'apr\`es l'hypoth\`ese de 
r\'ecurrence, le complexe $D_{j}^{\bullet}(W,E)$ est acyclique pour tout 
entier naturel $j$; donc par un raisonnement analogue \`a celui de (ii), 
$a$ est un cobord du complexe $D_{k}^{\bullet}(V,E)$, d'o\`u l'assertion.

iv) La diff\'erentielle du complexe $D_{k}^{\bullet}(V,L)$ est 
${\Bbb C}[X]$-lin\'eaire; donc la cohomologie de $D_{k}^{\bullet}(V,L)$
est un sous-quotient du ${\Bbb C}[X]$-module $D_{k}^{\bullet}(V,L)$.
En outre, la cohomologie de $D_{k}^{\bullet}(V,L)$ a une structure
naturelle d'espace gradu\'e qui est stable par la structure de 
${\Bbb C}[X]$-module. Par d\'efinition, $D^{\bullet}_{k}(V,L)$ est
somme directe d'un nombre fini de modules de type fini; donc la
cohomologie de $D_{k}^{\bullet}(V,L)$ est un \hbox{${\Bbb C}[X]$-module}
de type fini.
\end{proof}

Pour tout $x$ dans $X$, on note $L(x)$ l'image de $L$ par l'application
$\varphi \mapsto \varphi (x)$.

\begin{lemme}\label{l2co4}
Soit $X'$ l'ensemble des \'el\'ements $x$ de $X$ tels que $L(x)$ soit de 
dimension $r$. 

{\rm i)} La partie $X'$ de $X$ est ouverte et non vide. En outre, pour tout 
$x$ dans $X'$, il existe un sous-espace $E$ de $V$ et un ouvert $Y_{x}$ de $X$
qui contient $x$ tels que $E$ soit un suppl\'ementaire de $L(y)$ dans $V$
pour tout $y$ dans $Y_{x}$.

{\rm ii)} Pour tout entier naturel $k$, le support dans $X$ de la cohomologie 
du complexe $D_{k}^{\bullet}(V,L)$ ne rencontre pas $X'$.

{\rm iii)} On suppose que $X$ est une vari\'et\'e normale, que $X'$ est un 
grand ouvert de $X$ et que $L$ est un module libre. Alors le complexe 
$D_{k}^{\bullet}(V,L)$ n'a pas de cohomologie en degr\'e $r$.
\end{lemme}

\begin{proof}
i) Puisque $L$ est de rang $r$, $X$ contient des \'el\'ements $x$ tels que
$L(x)$ soit de dimension $r$. Soit $x$ dans $X'$. Alors il existe des 
\'el\'ements $\poi {\eta }1{,\ldots,}{r}{}{}{}$ de $L$ tels que la famille 
$\poi x{}{,\ldots,}{}{\eta }{1}{r}$ soit lin\'eairement libre. Il en r\'esulte
qu'il existe un sous-espace $E$ de $V$ et un ouvert $Y_{x}$ de $X$ qui
contient $x$ tels que $E$ soit un suppl\'ementaire dans $V$ du
sous-espace engendr\'e par $\poi y{}{,\ldots,}{}{\eta }{1}{r}$ pour tout
$y$ dans $Y_{x}$. Puisque $L$ est de rang $r$, pour tout $y$ dans un
ouvert non vide de $X$, $L(y)$ est de dimension inf\'erieure \`a $r$; or $X$
est irr\'eductible; donc pour tout $y$ dans un ouvert non vide de $Y_{x}$,
$L(y)$ est le sous-espace engendr\'e par $\poi y{}{,\ldots,}{}{\eta }{1}{r}$.
Il en r\'esulte que pour tout $\varphi $ dans $L$, le sous-module de $L$
engendr\'e par $\poi {\eta }1{,\ldots,}{r}{}{}{}$ contient le produit de
$\varphi $ par un \'el\'ement non nul de ${\Bbb C}[X]$. Par suite, il vient
$$ \varphi (y)\wedge \poi y{}{\wedge \cdots \wedge}{}{\eta }{1}{r}
= 0 \mbox{ ,}$$
pour tout $y$ dans $Y_{x}$ et pour tout $\varphi $ dans $L$; donc $L(y)$
est le sous-espace engendr\'e par $\poi y{}{,\ldots,}{}{\eta }{1}{r}$ pour
tout $y$ dans $Y_{x}$. Par suite, $X'$ contient $Y_{x}$. En outre,
$E$ et $Y_{x}$ satisfont la condition de l'assertion. Il en r\'esulte que 
$X'$ est recouvert par des ouverts affines $Y$ qui ont les propri\'et\'es 
suivantes:
\begin{list}{}{}
\item 1) la restriction de ${\cal L}$ \`a $Y$ est libre, 
\item 2) il existe un sous-espace $E$ de $V$ qui est un suppl\'ementaire de
$L(x)$ dans $V$ pour tout $x$ dans $Y$,
\end{list}
en d\'esignant par ${\cal L}$ la localisation de $L$ sur $X$. 

ii) Soient $k$ un entier naturel et $Y$ un ouvert affine de $X'$ qui satisfait
les conditions (1) et (2). On note $L_{Y}$ l'espace des sections de ${\cal L}$
au dessus de $Y$. Alors $D_{k}^{\bullet}(V,L_{Y})$ est l'espace des sections 
au dessus de $Y$ de la localisation sur $X$ de $D_{k}^{\bullet}(V,L)$. Vu 
l'arbitraire de $Y$, il s'agit de montrer que $D_{k}^{\bullet}(V,L_{Y})$ est 
un complexe acyclique.

Soit $x_{0}$ dans $Y$. De la condition (2), on d\'eduit qu'il existe une 
application r\'eguli\`ere de $Y$ dans GL$(V)$ qui \`a l'\'el\'ement $x$ 
associe l'automorphisme lin\'eaire qui prolonge l'identit\'e de $E$ et la 
restriction \`a $L_{Y}(x)$ de la projection lin\'eaire de $V$ sur 
$L_{Y}(x_{0})$ de noyau $E$. Soit $\overline{\tau }$ l'automorphisme de
l'alg\`ebre $\tk {{\Bbb C}}{{\Bbb C}[Y]}\tk {{\Bbb C}}{\es SV}\ex {}V$ qui \`a
l'\'el\'ement $\varphi $ associe l'application $x\mapsto \tau (x)(\varphi (x))$
de $Y$ dans $\tk {{\Bbb C}}{\es SV}\ex {}V$. Les images de $L_{Y}$ et de 
$D_{k}^{\bullet}(V,L_{Y})$ sont respectivement \'egales \`a
$\tk {{\Bbb C}}{{\Bbb C}[Y]}L_{Y}(x_{0})$ et \`a
$\tk {{\Bbb C}}{{\Bbb C}[Y]}D_{k}^{\bullet}(V,L_{Y}(x_{0}))$. En outre, la 
restriction de $\overline{\tau }$ \`a $D_{k}^{\bullet}(V,L_{Y})$ est un 
isomorphisme du complexe $D_{k}^{\bullet}(V,L_{Y})$ sur le complexe
$\tk {{\Bbb C}}{{\Bbb C}[Y]}D_{k}^{\bullet}(V,L_{Y}(x_{0}))$; donc d'apr\`es
l'assertion (iii) du lemme \ref{lco4}, le complexe
$D_{k}^{\bullet}(V,L_{Y})$ est acyclique.

iii) On note $\dd _{k}$ la diff\'erentielle du complexe $D_{k}^{\bullet}(V,L)$.
Soit $a$ un cocycle de degr\'e $r$ du complexe $D_{k}^{\bullet}(V,L)$. 
D'apr\`es l'assertion (ii), la restriction de $a$ \`a $Y$ est l'image d'un 
unique \'el\'ement $\varphi $ de $\tk {{\Bbb C}[Y]}{\sy k{L_{Y}}}
\ex r{L_{Y}}$ par $\dd _{k}$. Vu l'arbitraire de l'ouvert $Y$ de $X'$ 
satisfaisant les conditions (1) et (2) de (i), la restriction de $a$ \`a $X'$ 
est l'image par $\dd _{k}$ d'une section au dessus de $X'$ de la localisation 
sur $X$ de $\tk {{\Bbb C}[X]}{\sy kL}\ex rL$. Puisque $L$ est un module 
libre, $\tk {{\Bbb C}[X]}{\sy kL}\ex rL$ est un module libre et toute 
section au dessus de $X'$ de la localisation sur $X$ de 
$\tk {{\Bbb C}[X]}{\sy kL}\ex rL$ est la restriction \`a $X'$ d'un 
\'el\'ement $b$ de $\tk {{\Bbb C}}{\sy kL}\ex rL$ car $X'$ est un 
grand ouvert de la vari\'et\'e normale $X$; donc $a$ est \'egal \`a $\dd b$.
\end{proof}

\begin{cor}\label{cco4}
Soient $k$ un entier naturel et $L'$ un sous-module de $\tk {{\Bbb C}}{{\Bbb C}[X]}V$. 
On suppose que les conditions suivantes sont satisfaites:
\begin{list}{}{}
\item {\rm 1)} pour tout $x$ dans $X$, $L'(x)$ est de dimension $r$, 
\item {\rm 2)} pour tout $x$ dans un ouvert non vide de $X$, $L(x)$ est \'egal \`a
$L'(x)$,
\item {\rm 3)} les modules $\ex rL$ et $\ex r{L'}$ sont libres.
\end{list}
Alors le complexe $D_{k}^{\bullet}(V,L)$ n'a pas de cohomologie en degr\'e diff\'erent de
$r$.
\end{cor}

\begin{proof}
D'apr\`es la condition (1) et la d\'emonstration de l'assertion (ii) du lemme \ref{lco1},
l'\'el\'ement $\varphi $ de $\tk {{\Bbb C}}{{\Bbb C}[X]}V$ appartient \`a $L'$ si et
seulement si $L'(x)$ contient $\varphi (x)$ pour tout $x$ dans un ouvert non vide de
$X$ car $X$ est une vari\'et\'e affine. En particulier, $L$ est un sous-module de $L'$
d'apr\`es la condition (2). D'apr\`es la condition (1) et l'assertion (ii) du lemme 
\ref{l2co4}, le complexe $D_{k}^{\bullet}(V,L')$ est acyclique. Par d\'efinition, le 
complexe $D_{k}^{\bullet}(V,L)$ n'a pas de cohomologie en degr\'e inf\'erieur \`a $r-1$. 
D'apr\`es la condition (2), $L'$ et $L$ sont des sous-modules de rang $r$; donc d'apr\`es
la condition (3), il existe un \'el\'ement non nul $p$ de ${\Bbb C}[X]$ tel que $\ex rL$ 
soit \'egal \`a $p\ex r{L'}$. Soit $a$ un cocycle de degr\'e $i$, strictement sup\'erieur
\`a $r$, du complexe $D_{k}^{\bullet}(V,L)$. Alors on a $a=pb$ o\`u $b$ est un cocycle de 
degr\'e $i$ du complexe $D_{k}^{\bullet}(V,L')$. Soit $c$ un \'el\'ement de 
$D_{k}^{i-1}(V,L')$ dont le cobord est \'egal \`a $b$. Alors $D_{k}^{i-1}(V,L)$ contient 
$pc$ car $i$ est strictement sup\'erieur \`a $r$; or $a$ est le cobord de $pc$; donc $a$
est un cobord du complexe $D_{k}^{\bullet}(V,L)$.
\end{proof}
\subsection{} Pour tout point $x$ de $X$, on note ${\goth m}_{x}$ l'id\'eal
maximal de $\an Xx$ et $\han Xx$ le compl\'et\'e de $\an Xx$ pour la
topologie ${\goth m}_{x}$-adique. Soit $Y$ une sous-vari\'et\'e
ferm\'ee de codimension $m$ de $X$ qui a les propri\'et\'es suivantes:
pour tout $x$ dans $Y$, l'id\'eal de d\'efinition de $Y$ dans $\an Xx$ est
engendr\'e par $m$ coordonn\'ees $\poi y1{,\ldots,}{m}{}{}{}$, $\an Xx$
contient $\an Yx$ et l'application canonique de l'anneau de s\'eries
formelles $\an Yx[[\poi y1{,\ldots,}{m}{}{}{}]]$ dans $\han Xx$ est un
plongement dont l'image contient $\an Xx$. On note ${\goth m}_{x,Y}$
l'id\'eal de d\'efinition de $Y$ dans $\an Xx$. Pour tout espace vectoriel
$E$, pour tout $m$-uplet $\mu $ d'entiers naturels 
$\poi l1{,\ldots,}{m}{}{}{}$ et pour tout \'el\'ement $\psi $ de 
$\tk {{\Bbb C}}{\an Xx}E$, on d\'esigne par $\psi _{\mu }$ le
coefficient de $y_{1}^{l_{1}}\cdots y_{m}^{l_{m}}$ dans le
d\'eveloppement de $\psi $ suivant les puissances de 
$\poi y1{,\ldots,}{m}{}{}{}$. Alors 
$\psi _{\mu }$ appartient \`a $\tk {{\Bbb C}}{\an Yx}E$ pour tout
$\mu $. 

\begin{lemme}\label{lco5}
Soient $L$ un sous-module de $\tk {{\Bbb C}}{{\Bbb C}[X]}V$ et $L_{Y}$ le module des 
restrictions \`a $Y$ des \'el\'ements de $L$. Soient $l$ un entier naturel et $j$ un 
entier sup\'erieur au rang de $L$. On suppose que les
conditions suivantes sont satisfaites:
\begin{list}{}{}
\item {\rm 1)} le ${\Bbb C}[X]$-module $L$ et le 
${\Bbb C}[Y]$-module $L_{Y}$ sont des modules libres de rang $r$, 
\item {\rm 2)} le complexe $D^{\bullet}_{l}(V,L_{Y})$ n'a pas de cohomologie en degr\'e 
$j$.
\end{list}
Alors le support dans $X$ du $j$-i\`eme groupe de cohomologie du complexe 
$D^{\bullet}_{l}(V,L)$ ne rencontre pas $Y$.
\end{lemme}

\begin{proof}
Soient $M$ un ${\Bbb C}[X]$-module de type fini, $N$ et $N'$ deux
sous-modules de $M$ tels que $N$ contienne $N'$. Dire que le support dans
$X$ du quotient des modules $N$ et $N'$ ne rencontre pas $Y$ revient \`a
dire que pour tout $x$ dans $Y$, les deux modules
$\tk {{\Bbb C}[X]}{\an Xx}N$ et $\tk {{\Bbb C}[X]}{\an Xx}N'$ sont \'egaux.
Dans ce qui suit, on fixe un point $x$ de $Y$ et des g\'en\'erateurs
$\poi y1{,\ldots,}{m}{}{}{}$ de ${\goth m}_{x,Y}$. D\'esignant par
$\poi M0{,}{1}{}{}{},\ldots$ la suite d\'ecroissante de sous-modules de
$\tk {{\Bbb C}[X]}{\an Xx}M$ telle que 
$\tk {{\Bbb C}[X]}{{\goth m}_{x,Y}^{k}}M$ soit \'egale \`a $M_{k}$ pour tout 
$k$, l'\'egalit\'e des modules $\tk {{\Bbb C}[X]}{\an Xx}N$ et 
$\tk {{\Bbb C}[X]}{\an Xx}N'$ est \'equivalente \`a l'inclusion
$$ N\cap M_{k} \subset \tk {{\Bbb C}[Y]}{\an Xx}N' + M_{k+1} \mbox{ ,}$$
pour tout entier naturel $k$. On note $\poi {\eta }1{,\ldots,}{r}{}{}{}$
une base de $L$. Alors d'apr\`es la condition (1), la suite
$\poi {\overline{\eta }}1{,\ldots,}{r}{}{}{}$ des restrictions \`a $Y$
des \'el\'ements de la suite $\poi {\eta }1{,\ldots,}{r}{}{}{}$ est une
base du module $L_{Y}$.

Soient $M$ le module $\tk {{\Bbb C}[X]}{\an Xx}D^{j}_{l}(V,L)$,
$Z^{j}(V,L)$ l'espace des cocycles de degr\'e $j$ du complexe
$\tk {{\Bbb C}[X]}{\an Xx}D^{\bullet}_{l}(V,L)$ et $B^{j}(V, L)$ l'espace des cobords de 
degr\'e $j$ du complexe $\tk {{\Bbb C}[X]}{\an Xx}D^{\bullet}_{l}(V,L)$. Soient $k$ un 
entier naturel et $\psi $ dans l'intersection de $Z^{j}(V,L)$ et de $M_{k}$. D'apr\`es la
remarque pr\'eliminaire, il s'agit de montrer que $\psi $ est la somme d'un \'el\'ement 
de $\tk {{\Bbb C}[X]}{\an Xx}B^{j}(V,L)$ et d'un \'el\'ement de $M_{k+1}$. Soit 
$\kappa $ un $m$-uplet d'entiers naturels de somme $k$. Puisque 
$\poi {\overline{\eta }}1{,\ldots,}{r}{}{}{}$ est une base de $L_{Y}$, $\psi _{\kappa }$ 
est un cocycle de degr\'e $j$ du complexe $D^{\bullet}_{l}(V,L_{Y})$ car $j$ est 
sup\'erieur \`a $r$; donc d'apr\`es la condition (2), $\psi _{\kappa }$ est le cobord 
d'un \'el\'ement $\varphi _{\kappa }$ de 
\sloppy \hbox{$\tk {{\Bbb C}[X]}{\an Xx}D^{j-1}_{l}(V,L_{Y})$}. D\'esignant par 
$\tilde{\varphi }_{\kappa }$ un \'el\'ement de 
$\tk {{\Bbb C}[X]}{\an Xx}D^{j-1}_{l}(V,L)$ qui prolonge $\varphi _{\kappa }$ et posant:
$$ \varphi  = \sum_{\kappa =(\poi k1{,\ldots,}{m}{}{}{})\atop  
{\poi k1{+\cdots +}{m}{}{}{}} = k}  y_{1}^{k_{1}}\cdots y_{m}^{k_{m}}
\tilde{\varphi }_{\kappa } \mbox{ ,}$$
$\psi $ est la somme du cobord de $\varphi $ et d'un \'el\'ement de
$M_{k+1}$, d'o\`u le lemme.
\end{proof}

\section{Sous-module caract\'eristique.}\label{sc} 
Soi $L_{{\goth g}}$ le sous-module des \'el\'ements $\varphi $ du
$\e Sg$-module $\tk {{\Bbb C}}{\e Sg}{\goth g}$ qui satisfont la condition
suivante: pour tout $x$ dans ${\goth g}$, $\varphi (x)$ centralise $x$.
D'apr\`es \cite{Dix}(\S 2), $L_{{\goth g}}$ est un module libre de rang
$\rk {\goth g}$. On note ${\goth g}_{\r}$ l'ensemble des \'el\'ements
r\'eguliers de ${\goth g}$. Par d\'efinition, $x$ est dans ${\goth g}_{\r}$
si et seulement si ${\goth g}(x)$ est de dimension $\j g{}$. On introduit
dans cette section un sous-module $B_{{\goth g}}$ de
$\tk {{\Bbb C}}{\e Sg}\tk {{\Bbb C}}{\e Sg}{\goth g}$, appel\'e sous-module
caract\'eristique pour ${\goth g}$,  et on en donne quelques propri\'et\'es. 
On rappelle que $\k g{}$ d\'esigne la dimension des sous-alg\`ebres de Borel 
de ${\goth g}$.

\subsection{}  Le groupe $\A g$ des automorphismes de ${\goth g}$
op\`ere dans ${\goth g}\times {\goth g}$ par l'action diagonale. Pour tout
\'el\'ement $g$ de GL$_{2}({\Bbb C})$, on note $\kappa _{g}$
l'automorphisme lin\'eaire de ${\goth g}\times {\goth g}$,
$$ (x,y) \mapsto (ax+by,cx+dy) \mbox{ o\`u }
g = \left [ \begin{array}{cc} a & b \\ c & d \end{array} \right ] \mbox{ .}$$
On d\'efinit ainsi une action de GL$_{2}({\Bbb C})$ dans 
${\goth g}\times {\goth g}$. Pour tout $(x,y)$ dans 
${\goth g}\times {\goth g}$, on d\'esigne par ${\goth V}'(x,y)$ la somme
des sous-espaces ${\goth g}(ax+by)$ o\`u $(a,b)$ est dans 
${\Bbb C}^{2}\backslash \{0\}$. Alors pour tout $(x,y)$ dans 
${\goth g}\times {\goth g}$, on a
$$ g^{-1}({\goth V}'(g(x),g(y))) = {\goth V}'(\kappa _{h}(x,y)) =
{\goth V}'(x,y) \mbox{ ,}$$
pour tout $g$ dans $\A g$ et pour tout $h$ dans GL$_{2}({\Bbb C})$.

\begin{lemme}\label{lsc1}
Soit $\Omega _{{\goth g}}$ l'ensemble des \'el\'ements $(x,y)$ de 
${\goth g}\times {\goth g}$ tels que ${\goth g}_{\r}$ contienne 
$ax+by$ pour tout $(a,b)$ dans ${\Bbb C}^{2}\backslash \{0\}$.

{\rm i)} Le sous-ensemble $\Omega _{{\goth g}}$ de 
${\goth g}\times {\goth g}$ est un grand ouvert.

{\rm ii)} Soit $\Omega '_{{\goth g}}$ l'ensemble des \'el\'ements
$(x,y)$ de $\Omega _{{\goth g}}$ tels que ${\goth V}'(x,y)$ soit de
dimension maximale. Alors $\Omega '_{{\goth g}}$ est un ouvert de
${\goth g}\times {\goth g}$.

{\rm iii)} Les sous-ensembles $\Omega _{{\goth g}}$ et 
$\Omega '_{{\goth g}}$ sont invariants pour les actions de $\A g$
et de ${\rm GL}_{2}({\Bbb C})$.
\end{lemme}

\begin{proof}
i) Soient $\tau $ l'application
$$ {\goth g}\times {\goth g}\times ({\Bbb C}^{2}\backslash \{0\})
\rightarrow {\goth g} \mbox{ , } (x,y,a,b) \mapsto ax+by \mbox{ ,}$$
et $X$ l'image r\'eciproque de ${\goth g}\backslash {\goth g}_{\r}$ par
$\tau $. On d\'esigne par $Y$ l'image de $X$ par la projection canonique de
${\goth g}\times {\goth g}\times ({\Bbb C}^{2}\backslash \{0\})$ sur
${\goth g}\times {\goth g}$. La partie $X$ de  
${\goth g}\times {\goth g}\times ({\Bbb C}^{2}\backslash \{0\})$ 
est ferm\'ee car ${\goth g}_{\r}$ est un ouvert de ${\goth g}$. On note
$\tilde{X}$ l'image de $X$ par l'application canonique de 
${\goth g}\times {\goth g}\times ({\Bbb C}^{2}\backslash \{0\})$ sur
${\goth g}\times {\goth g}\times {\Bbb P}^{1}({\Bbb C})$.
Puisque ${\goth g}_{\r}$ est un c\^one, $X$ est l'image r\'eciproque de
$\tilde{X}$; donc $\tilde{X}$ est ferm\'e dans 
${\goth g}\times {\goth g}\times {\Bbb P}^{1}({\Bbb C})$. Il en r\'esulte
que $Y$ est ferm\'e dans ${\goth g}\times {\goth g}$. Par d\'efinition,
$\Omega _{{\goth g}}$ est le compl\'ementaire de $Y$ dans 
${\goth g}\times {\goth g}$; donc $\Omega _{{\goth g}}$ est ouvert dans
${\goth g}\times {\goth g}$. 

On suppose que $Y$ a une composante irr\'eductible $Y_{1}$ de
codimension $1$ dans ${\goth g}\times {\goth g}$. Il s'agit d'aboutir \`a
une contradiction. Puisque $Y_{1}$ est une composante irr\'eductible de
$Y$, $Y_{1}$ est l'image d'une composante irr\'eductible $X_{1}$ de $X$.
Le morphisme $\tau $ est une submersion en tout point; donc $\tau $ est
un morphisme lisse et la codimension de $X$ dans 
${\goth g}\times {\goth g}\times {\Bbb C}^{2}\backslash \{0\}$ est
\'egale \`a la codimension de ${\goth g}\backslash {\goth g}_{\r}$
dans ${\goth g}$. D'apr\`es \cite{V}(Theorem 4.12), la codimension de
${\goth g}\backslash {\goth g}_{\r}$ est \'egale \`a $3$; donc la
dimension de $X_{1}$ est \'egale \`a $2\dim {\goth g} -1$. Il en
r\'esulte que pour tout $(x,y)$ dans un ouvert non vide de $Y_{1}$,
l'ensemble des \'el\'ements $(a,b)$ de ${\Bbb C}^{2}\backslash \{0\}$
tels que $X_{1}$ contienne $(x,y,a,b)$, est fini. Si $X_{1}$ contient
$(x,y,a,b)$, alors $X_{1}$ contient $(x,y,ta,tb)$ pour tout \'el\'ement non
nul $t$ de ${\Bbb C}$ car ${\goth g}_{\r}$ est un c\^one et $X_{1}$ est une
composante irr\'eductible de $X$, d'o\`u la contradiction.

ii) Soit $(x,y)$ un \'el\'ement de $\Omega _{{\goth g}}$ tel que
${\goth V}'(x,y)$ soit de dimension maximale. Puisque ${\goth V}'(x,y)$ est
de dimension finie, il existe un nombre fini \sloppy
\hbox{$(a_{1},b_{1}),\ldots,(a_{k},b_{k})$} d'\'el\'ements de 
${\Bbb C}^{2}\backslash \{0\}$ tels que ${\goth V}'(x,y)$ soit la somme
des sous-espaces 
\hbox{${\goth g}(a_{1}x+b_{1}y),\ldots,{\goth g}(a_{k}x+b_{k}y)$}. Soit 
$\poi {\varepsilon }1{,\ldots,}{\j g{}}{}{}{}$ une base de $L_{{\goth g}}$.
D'apr\`es l'assertion (v) du lemme \ref{lco1}, pour tout \'el\'ement
r\'egulier $v$ de ${\goth g}$, $\poi v{}{,\ldots,}{}{\varepsilon }{1}{\j g{}}$
est une base de ${\goth g}(v)$; donc il existe des suites d'entiers
$\poi i{1,1}{,\ldots,}{l_{1},1}{}{}{},\ldots,
\poi i{1,k}{,\ldots,}{l_{k},k}{}{}{}$ dans $\{1,\ldots,\j g{}\}$ telles que la
suite
$$ \varepsilon _{i_{1,1}}(a_{1}x+b_{1}y) \mbox{ ,\ldots, }
\varepsilon _{i_{l_{1},1}}(a_{1}x+b_{1}y) \mbox{ $,\ldots$ }$$ $$
\mbox{ $\ldots,$ }\varepsilon _{i_{1,k}}(a_{k}x+b_{k}y) \mbox{ ,\ldots, }
\varepsilon _{i_{l_{k},k}}(a_{k}x+b_{k}y) \mbox{ ,}$$
soit une base de ${\goth V}'(x,y)$. Pour tout $(x',y')$ dans un ouvert de
${\goth g}\times {\goth g}$ contenant $(x,y)$, l'ensemble des \'el\'ements
$$ \varepsilon _{i_{1,1}}(a_{1}x'+b_{1}y') \mbox{ ,\ldots, }
\varepsilon _{i_{l_{1},1}}(a_{1}x'+b_{1}y') \mbox{ $,\ldots$ }$$ $$
\mbox{ $\ldots,$ }\varepsilon _{i_{1,k}}(a_{k}x'+b_{k}y') \mbox{ ,\ldots, }
\varepsilon _{i_{l_{k},k}}(a_{k}x'+b_{k}y') \mbox{ ,}$$ 
est une partie libre de ${\goth V}'(x',y')$; donc pour tout $(x',y')$ dans un
ouvert de $\Omega _{{\goth g}}$ contenant $(x,y)$, ${\goth V}'(x,y)$ et
${\goth V}'(x',y')$ ont m\^eme dimension d'apr\`es la maximalit\'e de la
dimension de ${\goth V}'(x,y)$, d'o\`u l'assertion.

iii) Soit $(x,y)$ dans $\Omega _{{\goth g}}$. Puisque ${\goth g}_{\r}$ est
$\A g$-invariant, $\Omega _{{\goth g}}$ contient l'orbite de $(x,y)$ pour
l'action de $\A g$. Le plan $P_{x,y}$ de ${\goth g}$, engendr\'e par $x$ et
$y$, ne d\'ependant que de l'orbite de $(x,y)$ sous l'action de
GL$_{2}({\Bbb C})$, $\Omega _{{\goth g}}$ contient l'orbite de $(x,y)$
pour l'action de GL$_{2}({\Bbb C})$. D'apr\`es les \'egalit\'es
$$ g^{-1}({\goth V}'(g(x),g(y))) = {\goth V}'(\kappa _{h}(x,y)) =
{\goth V}'(x,y) \mbox{ ,}$$
pour tout $g$ dans $\A g$ et pour tout $h$ dans GL$_{2}({\Bbb C})$,
$\Omega '_{{\goth g}}$ contient les orbites de ses \'el\'ements sous les actions
de $\A g$ et GL$_{2}({\Bbb C})$, d'o\`u l'assertion.
\end{proof}

\begin{cor}\label{csc1}
Soit $d$ la dimension de ${\goth V}'(x,y)$ pour $(x,y)$ dans 
$\Omega '_{{\goth g}}$. Alors l'application $(x,y) \mapsto {\goth V}'(x,y)$
de $\Omega '_{{\goth g}}$ dans $\ec {Gr}g{}{}d$ est r\'eguli\`ere.
\end{cor}

\begin{proof}
Soit $(x,y)$ dans $\Omega '_{{\goth g}}$. Selon les notations de
la d\'emonstration de l'assertion (ii) du lemme \ref{lsc1}, pour tout
$(x',y')$ dans un ouvert de $\Omega _{{\goth g}}$ contenant $(x,y)$, la
famille
$$ \varepsilon _{i_{1,1}}(a_{1}x'+b_{1}y') \mbox{ ,\ldots, }
\varepsilon _{i_{l_{1},1}}(a_{1}x'+b_{1}y') \mbox{ $,\ldots$ }$$ $$
\mbox{ $\ldots,$ }\varepsilon _{i_{1,k}}(a_{k}x'+b_{k}y') \mbox{ ,\ldots, }
\varepsilon _{i_{l_{k},k}}(a_{k}x'+b_{k}y') \mbox{ ,}$$ 
est une base de ${\goth V}'(x',y')$. L'application 
$(x',y')\mapsto \varepsilon _{i,j}(a_{j}x'+b_{j}y')$ \'etant
r\'eguli\`ere pour tout $(i,j)$, l'application 
$(x',y') \mapsto {\goth V}'(x',y')$ est r\'eguli\`ere en $(x,y)$, d'o\`u
l'assertion.
\end{proof}\subsection{} Soient $\poi p1{,\ldots,}{\j g{}}{}{}{}$ des \'el\'ements
homog\`enes de $\e Sg$ qui engendrent la sous-alg\`ebre des \'el\'ements
$G$-invariants de $\e Sg$. Pour $i=1,\ldots,\j g{}$, pour $(x,y)$ dans 
${\goth g}\times {\goth g}$ et pour $(a,b)$ dans ${\Bbb C}^{2}$, on a un
unique d\'eveloppement
$$ p_{i}(ax+by) = \sum_{(m,n)\in {\Bbb N}^{2}}
a^{m}b^{n}p_{i,m,n}(x,y) \mbox{ ,}$$
o\`u $p_{i,m,n}$ est un \'el\'ement de $\tk {{\Bbb C}}{\e Sg}\e Sg$ qui est
nul d\`es que $m$ ou $n$ est assez grand. Puisque $p_{i}$ est
$G$-invariant, les \'el\'ements $p_{i,m,n}$ sont $G$-invariants pour
l'action de $G$ dans l'alg\`ebre $\tk {{\Bbb C}}{\e Sg}\e Sg$ qui prolonge
l'action diagonale de $G$ dans ${\goth g}\times {\goth g}$. Pour
$i=1,\ldots,\j g{}$ et pour $(m,n)$ dans ${\Bbb N}^{2}$, on note
$\varepsilon _{i}$ l'\'el\'ement de $\tk {{\Bbb C}}{\e Sg}{\goth g}$
et $\varepsilon _{i,m,n}$ l'\'el\'ement de 
$\tk {{\Bbb C}}{\e Sg}\tk {{\Bbb C}}{\e Sg}{\goth g}$ qui sont d\'efinis par
la condition suivante: pour tout $(x,y)$ dans ${\goth g}\times {\goth g}$,
les formes lin\'eaires sur ${\goth g}$, 
$v\mapsto \dv {\varepsilon _{i}(x)}v$ et 
$v\mapsto \dv {\varepsilon _{i,m,n}(x,y)}v$ sont respectivement les
diff\'erentielles en $x$ des fonctions $p_{i}$ et 
$x\mapsto p_{i,m,n}(x,y)$ sur ${\goth g}$. D'apr\`es \cite{Dix}(\S 2),
$\poi {\varepsilon }1{,\ldots,}{\j g{}}{}{}{}$ est une base de
$L_{{\goth g}}$. On note $\underline{E}$ le sous-ensemble de 
$\tk {{\Bbb C}}{\e Sg}\tk {{\Bbb C}}{\e Sg}{\goth g}$,
$$ \underline{E} = \{\varepsilon _{i,m,n} \mbox{ ; }
m \in {\Bbb N}\backslash \{0\} \mbox{ , } n\in {\Bbb N} \mbox{ , }
m+n = d_{i}\mbox{ , } i=1,\ldots,\j g{}\} \mbox{ ,}$$
o\`u $d_{i}$ est le degr\'e de $p_{i}$ pour $i=1,\ldots,\j g{}$, et
$\underline{e}(x,y)$ l'image de $\underline{E}$ par l'application 
$\varphi \mapsto \varphi (x,y)$ pour tout $(x,y)$ dans 
${\goth g}\times {\goth g}$.

\begin{lemme}\label{lsc2}
{\rm i)} Pour $i=1,\ldots,\j g{}$ et pour tout $(m,n)$ dans ${\Bbb N}^{2}$,
$p_{i,m,n}$ est homog\`ene de degr\'e $m$ en la premi\`ere composante
sur ${\goth g}\times {\goth g}$ et homog\`ene de degr\'e $n$ en la
deuxi\`eme composante sur ${\goth g}\times {\goth g}$. En outre,
$p_{i,m,n}$ est nul si $m+n$ est diff\'erent de $d_{i}$.

{\rm ii)} Pour tout $(x,y)$ dans ${\goth g}\times {\goth g}$ et pour
$i=1,\ldots,\j g{}$, on a
$$ \varepsilon _{i}(ax+by) =  \sum_{m=1}^{d_{i}} 
a^{m-1}b^{d_{i}-m}\varepsilon _{i,m,d_{i}-m}(x,y) \mbox{ ,}$$
pour tout $(a,b)$ dans ${\Bbb C}^{2}$.

{\rm iii)} Pour tout $(x,y)$ dans $\Omega _{{\goth g}}$,
$\underline{e}(x,y)$ engendre le sous-espace ${\goth V}'(x,y)$.

{\rm iv)} Pour $i=1,\ldots,\j g{}$ et pour tout $(m,n)$ dans ${\Bbb N}^{2}$,
$\varepsilon _{i,m,n}(x,y)$ est orthogonal \`a $[x,y]$ pour tout $(x,y)$
dans ${\Bbb N}^{2}$.

{\rm v)} Pour $i=1,\ldots,\j g{}$ et pour tout $(g,m,n)$ dans 
$G\times {\Bbb N}^{2}$, $\varepsilon _{i,m,n}(g(x),g(y))$ est \'egal \`a
$g(\varepsilon _{i,m,n}(x,y))$.
\end{lemme}

\begin{proof}
i) Soient $i=1,\ldots,\j g{}$ et $(x,y)$ dans ${\goth g}\times {\goth g}$.
Pour tout $(a,b,t)$ dans ${\Bbb C}^{3}$, on a
$$ \sum_{(m,n)\in {\Bbb N}^{2}} a^{m}b^{n} p_{i,m,n}(tx,y) =
p_{i}(atx+by) = \sum_{(m,n)\in {\Bbb N}^{2}} a^{m}b^{n}t^{m} p_{i,m,n}(x,y)
\mbox{ ;}$$
donc pour tout $(m,n)$ dans ${\Bbb N}^{2}$, $p_{i,m,n}(tx,y)$ est \'egal \`a
$t^{m}p_{i,m,n}(x,y)$ pour tout $t$ dans ${\Bbb C}$. On montre de m\^eme
que $p_{i,m,n}$ est homog\`ene de degr\'e $n$ en la deuxi\`eme
composante sur ${\goth g}\times {\goth g}$. Puisque $p_{i}$ est
homog\`ene de degr\'e $d_{i}$, le degr\'e total de $p_{i,m,n}$ est \'egal
\`a $d_{i}$ pour tout $(m,n)$; donc $p_{i,m,n}$ est nul si $m+n$ est
diff\'erent de $d_{i}$.

ii) Soient $i=1,\ldots,\j g{}$, $(x,y)$ dans ${\goth g}\times {\goth g}$ et
$v$ un \'el\'ement de ${\goth g}$. Pour tout \'el\'ement $(a,b)$ de 
${\Bbb C}^{2}$ tel que $a$ soit non nul, on a
$$ \dv {\varepsilon _{i}(ax+by)}v =
\frac{\dd }{\dd t}p_{i}(ax+by+tv) \left \vert \right. _{t=0} $$ $$ \mbox{ }
= \sum_{(m,n)\in {\Bbb N}^{2}} 
a^{m}b^{n} \frac{\dd }{\dd t}p_{i,m,n}(x+ta^{-1}v,y) \left \vert \right.
_{t=0}  = \sum_{(m,n)\in {\Bbb N}^{2}} 
a^{m}b^{n} \dv {\varepsilon _{i,m,n}(x,y)}{a^{-1}v}  \mbox{ .}$$
L'assertion r\'esulte alors de l'assertion (i).

iii) Soient $(x,y)$ dans $\Omega _{{\goth g}}$ et $v$ un \'el\'ement de
${\goth g}$. D'apr\`es l'assertion (ii), $v$ est orthogonal \`a 
$\varepsilon _{i,m,n}(x,y)$ pour tout $(m,n)$
dans $({\Bbb N}\backslash \{0\})\times {\Bbb N}$ si et seulement si $v$
est orthogonal \`a $\varepsilon _{i}(ax+by)$ pour tout $(a,b)$
dans ${\Bbb C}^{2}$; donc $v$ est orthogonal \`a $\underline{e}(x,y)$ si et 
seulement si $v$ est orthogonal \`a 
\hbox{$\poi {ax+by}{}{,\ldots,}{}{\varepsilon }{1}{\j g{}}$} pour tout 
$(a,b)$ dans ${\Bbb C}^{2}$. Pour tout \'el\'ement non nul $(a,b)$ de 
${\Bbb C}^{2}$, ${\goth g}_{\r}$ contient $ax+by$; donc d'apr\`es
l'assertion (v) du lemme \ref{lco1},  
$\poi {ax+by}{}{,\ldots,}{}{\varepsilon }{1}{\j g{}}$ est une base de 
${\goth g}(ax+by)$. Il en r\'esulte que $v$ est orhogonal \`a 
${\goth V}'(x,y)$ si et seulement si $v$ est orthogonal \`a 
$\underline{e}(x,y)$, d'o\`u l'assertion.

iv) Soient $i=1,\ldots,\j g{}$ et $(x,y)$ dans ${\goth g}\times {\goth g}$.
D'apr\`es la propri\'et\'e d'invariance de $\dv ..$, pour tout
\'el\'ement $(a,b)$ de ${\Bbb C}^{2}$ tel que $a$ soit non nul, on a
$$ \dv {\varepsilon _{i}(ax+by)}{[x,y]} = a^{-1} 
\dv {\varepsilon _{i}(ax+by)}{[ax+by,y]} $$ $$ \mbox{ }= a^{-1}
\dv {[\varepsilon _{i}(ax+by),ax+by]}{y} = 0 \mbox{ ,}$$
car ${\goth g}(ax+by)$ contient $\varepsilon _{i}(ax+by)$. Il r\'esulte
alors de l'assertion (ii) que $\varepsilon _{i,m,n}(x,y)$ est orthogonal \`a
$[x,y]$ pour tout couple $(m,n)$ d'entiers naturels tels que $m$ soit non nul;
or par d\'efinition, $\varepsilon _{i,0,n}(x,y)$ est nul pour tout $n$;
donc $\varepsilon _{i,m,n}(x,y)$ est orthogonal \`a $[x,y]$ pour tout
$(m,n)$ dans ${\Bbb N}^{2}$.

v) Soient $i=1,\ldots,\j g{}$ et $(m,n)$ dans ${\Bbb N}^{2}$. Puisque
$p_{i}$ est $G$-invariant, $p_{i,m,n}$ est $G$-invariant; donc pour tout
$g$ dans $G$ et pour tout $v$ dans ${\goth g}$, on a
$$ \dv {\varepsilon _{i,m,n}(g(x),g(y))}{v} = 
\frac{\dd }{\dd t}p_{i,m,n}(x+tg^{-1}(v),y) \left \vert \right. _{t=0} =
\dv {\varepsilon _{i,m,n}(x,y)}{g^{-1}(v)} \mbox{ ,}$$
pour tout $(x,y)$ dans ${\goth g}\times {\goth g}$, d'o\`u l'assertion.
\end{proof}

\begin{cor}\label{csc2}
Soit $\Lambda _{{\goth g}}$ l'ensemble des \'el\'ements $(x,y)$ de
${\goth g}\times {\goth g}$ tels que $\underline{e}(x,y)$ soit de rang
maximum.
 
{\rm i)} Le cardinal de $\underline{E}$ est inf\'erieur \`a $\k g{}$.

{\rm ii)} Pour tout $(x,y)$ dans $\Omega _{{\goth g}}$, la dimension de 
${\goth V}'(x,y)$ est inf\'erieure \`a $\k g{}$.

{\rm iii)} La partie $\Lambda _{{\goth g}}$ de ${\goth g}\times {\goth g}$
est un ouvert $G$-invariant qui contient $\Omega '_{{\goth g}}$.
\end{cor}

\begin{proof}
i) Par d\'efinition, $\underline{E}$ a au plus 
$\poi d1{+\cdots +}{\j g{}}{}{}{}$ \'el\'ements distincts; or d'apr\`es 
\cite{Bo}(Ch. V, \S 5, Proposition 3), cette somme est \'egale \`a 
$\k g{}$; donc le cardinal de $\underline{E}$ est inf\'erieur \`a $\k g{}$.

ii) D'apr\`es (i), pour tout $(x,y)$ dans ${\goth g}\times {\goth g}$, le
cardinal de $\underline{e}(x,y)$ est inf\'erieur \`a $\k g{}$; donc d'apr\`es
l'assertion (iii) du lemme \ref{lsc2}, pour tout $(x,y)$ dans 
$\Omega _{{\goth g}}$, la dimension de ${\goth V}'(x,y)$ est inf\'erieure
\`a $\k g{}$.

iii) Par d\'efinition, $\Lambda _{{\goth g}}$ est ouvert dans 
${\goth g}\times {\goth g}$. Cet ouvert est $G$-invariant d'apr\`es
l'assertion (v) du lemme \ref{lsc2}. En outre, d'apr\`es l'assertion (iii) du
lemme \ref{lsc2}, $\Lambda _{{\goth g}}$ contient $\Omega '_{{\goth g}}$
car $\Omega '_{{\goth g}}$ est l'ensemble des \'el\'ements $(x,y)$ de
$\Omega _{{\goth g}}$ tels que ${\goth V}'(x,y)$ soit de dimension
maximale. 
\end{proof}

\subsection{} On suppose que ${\goth g}$ n'est pas une alg\`ebre de Lie
commutative. Soient $(\xi ,\rho ,\eta )$ un ${\goth s}{\goth l}_{2}$-triplet principal
de l'alg\`ebre d\'eriv\'ee de ${\goth g}$. Par d\'efinition, on a
$$[\xi ,\eta ] = \rho \mbox{ , } [\rho ,\xi ] = 2\xi \mbox{ , } 
[\rho ,\eta ] = -2 \eta \mbox{ .}$$
Puisque $\xi $ est un \'el\'ement nilpotent r\'egulier de ${\goth g}$, il existe une 
unique sous-alg\`ebre de Borel ${\goth b}$ de ${\goth g}$ qui contient $\xi $. En outre,
le centralisateur ${\goth h}$ de $\rho $ dans ${\goth g}$ est une sous-alg\`ebre de
Cartan de ${\goth g}$ qui est contenue dans ${\goth b}$. On d\'esigne par $R$ le 
syst\`eme de racines de ${\goth h}$ dans ${\goth g}$ et par $R_{+}$ le syst\`eme de 
racines positives de $R$ d\'efini par ${\goth b}$. Pour tout $\alpha $ dans $R$,
${\goth g}^{\alpha }$ d\'esigne le sous-espace radiciel de racine $\alpha $.

\begin{lemme}\label{lsc3}
{\rm i)} L'ouvert $\Omega _{{\goth g}}$ contient $(\rho ,\xi )$.

{\rm ii)} L'espace ${\goth V}'(\rho ,\xi )$ est \'egal \`a la somme des
sous-espaces \sloppy
\hbox{${\goth h},[\xi ,{\goth h}],(\ad \xi )^{2}({\goth h}),\ldots$}.

{\rm iii)} L'espace ${\goth V}'(\rho ,\xi )$ est \'egal \`a ${\goth b}$.
\end{lemme}

\begin{proof}
i) Pour tout $t$ dans ${\Bbb C}$, $\rho +2t\xi $ est \'egal \`a 
$\exp(-t\ad \xi )(\rho )$; donc $\rho +2t\xi $ est un \'el\'ement semi-simple
r\'egulier de ${\goth g}$. En outre, $\xi $ est r\'egulier; donc 
$\Omega _{{\goth g}}$ contient $(\rho ,\xi )$ car ${\goth g}_{\r}$ est un c\^one.

ii) Soit $V$ la somme des sous-espaces 
\hbox{${\goth h},[\xi ,{\goth h}],(\ad \xi )^{2}({\goth h}),\ldots$}. Puisque 
${\goth g}(\rho )$ est \'egal \`a ${\goth h}$, pour tout $t$ dans ${\Bbb C}$, le 
centralisateur de $\rho +2t\xi $ dans ${\goth g}$ est \'egal
\`a $\exp(-t\ad \xi )({\goth h})$. Si $\varphi $ est une forme lin\'eaire sur
${\goth g}$ qui est nulle sur ${\goth V}'(\rho ,\xi )$, alors pour tout $v$ dans
${\goth h}$, $\varphi $ est nul en $(\ad \xi )^{i}(v)$ pour tout entier
naturel $i$; donc ${\goth V}'(\rho ,\xi )$ contient $V$. R\'eciproquement, $V$
contient ${\goth g}(\rho +2t\xi )$ pour tout $t$ dans ${\Bbb C}$. Pour tout $a$ dans 
${\Bbb C}\backslash \{0\}$, ${\goth g}(a\rho +b\xi )$ est \'egal \`a 
${\goth g}(\rho +a^{-1}b\xi )$; donc $V$ contient ${\goth g}(a\rho +b\xi )$ pour tout $a$
dans ${\Bbb C}\backslash \{0\}$. Puisque $\xi $ est un \'el\'ement nilpotent r\'egulier 
de ${\goth g}$, ${\goth g}(\xi )$ est la limite dans $\ec {Gr}g{}{}{\j g{}}$ de 
${\goth g}(a\rho +\xi )$ quand $a$ tend vers $0$; donc $V$ contient ${\goth g}(\xi )$, 
d'o\`u l'assertion.

iii) Soient ${\goth u}$ le sous-espace des \'el\'ements nilpotents de
${\goth b}$ et ${\goth u}_{-}$ la somme des sous-espaces ${\goth g}^{\alpha }$ o\`u 
$\alpha $ est dans $-R_{+}$. D'apr\`es l'assertion (ii), ${\goth b}$ contient 
${\goth V}'(\rho ,\xi )$. On suppose que ${\goth b}$ contient strictement 
${\goth V}'(\rho ,\xi )$. Il s'agit d'aboutir \`a une contradiction. Soit 
${\goth V}'(\rho ,\xi )^{\perp}$ l'orthogonal de ${\goth V}'(\rho ,\xi )$ dans 
${\goth g}$. Alors ${\goth V}'(\rho ,\xi )^{\perp}$ contient strictement 
${\goth u}$. Puisque $[\rho ,\xi ]$ est \'egal \`a $2\xi $, ${\goth V}'(\rho ,\xi )$ et 
${\goth V}'(\rho ,\xi )^{\perp}$ sont stables par $\ad \rho $ et $\ad \xi $ d'apr\`es 
l'assertion (ii). Le sous-espace ${\goth u}$ \'etant la somme des sous-espaces propres de
$\ad \rho $ pour les valeurs propres strictement positives, la restriction de 
$\ad \rho $ \`a ${\goth V}'(\rho ,\xi )^{\perp}$ a une valeur propre n\'egative. Soient 
$i$ la plus grande valeur propre n\'egative de la restriction de $\ad \rho $ \`a
${\goth V}'(\rho ,\xi )^{\perp}$ et $v$ un vecteur propre non nul de cet
endomorphisme pour la valeur propre $i$. Alors $[\xi ,v]$ est un vecteur
propre de la restriction de $\ad \rho $ \`a ${\goth V}'(\rho ,\xi )^{\perp}$ pour la
valeur propre $i+2$ car $[\rho ,\xi ]$ est \'egal \`a $2\xi $; or $[\xi ,v]$ n'est pas nul
car ${\goth u}$ contient ${\goth g}(\xi )$; donc d'apr\`es la
maximalit\'e de $i$, $i$ est nul et $v$ appartient \`a ${\goth h}$ car
${\goth g}(\rho )$ est \'egal \`a ${\goth h}$. Ceci est absurde car ${\goth h}$
est un sous-espace de ${\goth V}'(\rho ,\xi )$ dont l'orthogonal dans 
${\goth g}$ a une intersection nulle avec ${\goth h}$.
\end{proof}

\begin{cor}\label{csc3}
{\rm i)} Pour tout $(x,y)$ dans $\Omega '_{{\goth g}}$, ${\goth V}'(x,y)$ est
de dimension $\k g{}$.

{\rm ii)} Pour tout $(x,y)$ dans $\Omega '_{{\goth g}}$, le sous-ensemble
$\underline{e}(x,y)$ est une partie libre de ${\goth g}$, de cardinal $\k g{}$.

{\rm iii)} Le sous-ensemble $\underline{E}$ de 
$\tk {{\Bbb C}}{\e Sg}\tk {{\Bbb C}}{\e Sg}{\goth g}$ est une partie libre
de cardinal $\k g{}$.
\end{cor}

\begin{proof}
i) D'apr\`es l'assertion (i) du lemme \ref{lsc3}, $\Omega _{{\goth g}}$
contient $(\rho ,\xi )$; donc d'apr\`es l'assertion (iii) du lemme \ref{lsc3} et le
corollaire \ref{csc2}, $\Omega '_{{\goth g}}$ contient $(\rho ,\xi )$. Il en
r\'esulte que pour tout $(x,y)$ dans $\Omega '_{{\goth g}}$, les
dimensions de ${\goth V}'(x,y)$ et de ${\goth b}$ sont \'egales.

ii) D'apr\`es l'assertion (iii) du lemme \ref{lsc2}, $\underline{e}(x,y)$
engendre ${\goth V}'(x,y)$ pour tout $(x,y)$ dans $\Omega _{{\goth g}}$; or
$\underline{e}(x,y)$ a au plus $\k g{}$ \'el\'ements d'apr\`es
l'assertion (i) du corollaire \ref{csc2}; donc pour tout $(x,y)$ dans
$\Omega '_{{\goth g}}$, $\underline{e}(x,y)$ est une partie libre de
${\goth g}$, de cardinal $\k g{}$, d'apr\`es l'assertion (i).

iii) Pour tout $(x,y)$ dans ${\goth g}\times {\goth g}$,
$\underline{e}(x,y)$ est l'image de $\underline{E}$ par l'application
$\varphi \mapsto \varphi (x,y)$; or d'apr\`es l'assertion (ii) du lemme
\ref{lsc1}, $\Omega '_{{\goth g}}$ est un ouvert de 
${\goth g}\times {\goth g}$; donc d'apr\`es l'assertion (ii), $\underline{E}$
est une partie libre de 
$\tk {{\Bbb C}}{\e Sg}\tk {{\Bbb C}}{\e Sg}{\goth g}$. D'apr\`es l'assertion
(i) du corollaire \ref{csc2}, le cardinal de $\underline{E}$ est inf\'erieur
\`a $\k g{}$; donc d'apr\`es l'assertion (ii), le cardinal de $\underline{E}$
est \'egal \`a $\k g{}$.
\end{proof}

\begin{cor}\label{c2sc3}
Pour ${\goth g}$ alg\`ebre de Lie r\'eductive quelconque, l'ouvert 
$\Lambda _{{\goth g}}$ de ${\goth g}\times {\goth g}$ est un grand ouvert.
\end{cor}

\begin{proof}
Si ${\goth g}$ est une alg\`ebre de Lie commutative, $\Lambda _{{\goth g}}$ est 
\'egal \`a ${\goth g}\times {\goth g}$. On peut donc supposer que ${\goth g}$ n'est pas
une alg\`ebre de Lie commutative. D'apr\`es l'assertion (iii) du corollaire \ref{csc2}, 
$\Lambda _{{\goth g}}$ est un ouvert $G$-invariant de ${\goth g}\times {\goth g}$. On 
suppose que $\Lambda _{{\goth g}}$ n'est pas un grand ouvert de ${\goth g}$. Il s'agit 
d'aboutir \`a une contradiction. Alors le compl\'ementaire de $\Lambda _{{\goth g}}$ dans
${\goth g}\times {\goth g}$ a une composante irr\'eductible $\Sigma $ de codimension $1$ 
dans ${\goth g}\times {\goth g}$. D'apr\`es l'assertion (i) du lemme \ref{lsc1}, 
$\Sigma $ rencontre $\Omega _{{\goth g}}$. En outre, d'apr\`es le corollaire \ref{csc2}, 
l'intersection $\Sigma '$ de $\Sigma $ et de $\Omega _{{\goth g}}$ est une composante 
irr\'eductible du compl\'ementaire de $\Omega '_{{\goth g}}$ dans $\Omega _{{\goth g}}$; 
or d'apr\`es l'assertion (iii) du lemme \ref{lsc2}, pour tout $(x,y)$ dans 
$\Omega _{{\goth g}}$, ${\goth V}'(x,y)$ ne d\'epend que du sous-espace engendr\'e par 
$x$ et $y$; donc $\Sigma '$ et $\Sigma $ sont stables pour l'action de 
GL$_{2}({\Bbb C})$. En outre, $\Sigma $ est $G$-invariant
car $\Lambda _{{\goth g}}$ est $G$-invariant. Soient ${\goth s}$ le sous-espace de 
${\goth g}$ engendr\'e par $\xi $, $\rho $, $\eta $, ${\bf S}$ le sous-groupe ferm\'e 
connexe de $G$ dont l'alg\`ebre de Lie est l'image de ${\goth s}$ par la repr\'esentation
adjointe de ${\goth s}$ et $T$ l'intersection de $\Sigma $ et de 
${\goth s}\times {\goth s}$. Puisque $\Sigma $ est une hypersurface de 
${\goth g}\times {\goth g}$, la codimension de $T$ dans ${\goth s}\times {\goth s}$ est 
inf\'erieure \`a $1$. En outre, $T$ est invariant pour les actions de ${\bf S}$ et de 
GL$_{2}({\Bbb C})$; donc d'apr\`es le lemme \ref{lsu1}, $T$ contient $(\rho ,\xi )$. 
Ceci est absurde car $\Lambda _{{\goth g}}$ contient $(\rho ,\xi )$ d'apr\`es 
l'assertion (iii) du lemme \ref{lsc3}.
\end{proof}

Pour tout $(x,y)$ dans ${\goth g}\times {\goth g}$, on note ${\goth V}(x,y)$ 
le sous-espace de ${\goth g}$ engendr\'e par $\underline{e}(x,y)$. 

\begin{prop}\label{psc}
Soit $B_{{\goth g}}$ le sous-module des \'el\'ements $\varphi $ de 
$\tk {{\Bbb C}}{\e Sg}\tk {{\Bbb C}}{\e Sg}{\goth g}$ tels que 
${\goth V}'(x,y)$ contienne $\varphi (x,y)$ pour tout $(x,y)$ dans un
ouvert non vide de ${\goth g}\times {\goth g}$. 

{\rm i)} L'\'el\'ement $\varphi $ de 
$\tk {{\Bbb C}}{\e Sg}\tk {{\Bbb C}}{\e Sg}{\goth g}$ appartient \`a 
$B_{{\goth g}}$ si et seulement si ${\goth V}'(x,y)$ contient 
$\varphi (x,y)$ pour tout $(x,y)$ dans $\Omega '_{{\goth g}}$.

{\rm ii)} L'\'el\'ement $\varphi $ de 
$\tk {{\Bbb C}}{\e Sg}\tk {{\Bbb C}}{\e Sg}{\goth g}$ appartient \`a
$B_{{\goth g}}$ si et seulement si ${\goth V}(x,y)$ contient 
$\varphi (x,y)$ pour tout $(x,y)$ dans $\Lambda _{{\goth g}}$.

{\rm iii)} Le module $B_{{\goth g}}$ est un module libre de rang \'egal \`a
$\k g{}$. En outre, $\underline{E}$ est une partie libre et g\'en\'eratrice
de $B_{{\goth g}}$.

{\rm iv)} Pour tout $g$ dans $\A g$ et pour tout $\varphi $ dans 
$B_{{\goth g}}$, $B_{{\goth g}}$ contient l'ap\-plication
$(x,y) \mapsto g^{-1}\rond \varphi (g(x),g(y))$.

{\rm v)} Pour tout $g$ dans ${\rm GL}_{2}({\Bbb C})$ et pour tout
$\varphi $ dans  $B_{{\goth g}}$, $B_{{\goth g}}$ contient 
$\varphi \rond \kappa _{g}$.
\end{prop}

\begin{proof}
Si ${\goth g}$ est une alg\`ebre de Lie commutative, $B_{{\goth g}}$ est
\'egal \`a $\tk {{\Bbb C}}{\e Sg}\tk {{\Bbb C}}{\e Sg}{\goth g}$ car 
${\goth V}(x,y)$ est \'egal \`a ${\goth g}$ pour tout $(x,y)$ dans 
${\goth g}\times {\goth g}$. On peut donc supposer que ${\goth g}$ n'est
pas une alg\`ebre de Lie commutative. 

i) Soit $\varphi $ dans 
$\tk {{\Bbb C}}{\e Sg}\tk {{\Bbb C}}{\e Sg}{\goth g}$. D'apr\`es l'assertion
(ii) du lemme \ref{lsc1}, $B_{{\goth g}}$ contient $\varphi $ si et
seulement si ${\goth V}'(x,y)$ contient $\varphi (x,y)$ pour tout $(x,y)$
dans un ouvert non vide de $\Omega '_{{\goth g}}$; donc d'apr\`es le
corollaire \ref{csc1} et par un argument de continuit\'e, 
$B_{{\goth g}}$ contient $\varphi $ si et seulement si ${\goth V}'(x,y)$
contient $\varphi (x,y)$ pour tout $(x,y)$ dans $\Omega '_{{\goth g}}$.

ii) Soit $\varphi $ dans 
$\tk {{\Bbb C}}{\e Sg}\tk {{\Bbb C}}{\e Sg}{\goth g}$. L'application 
$(x,y) \mapsto {\goth V}(x,y)$ est une application r\'eguli\`ere de
$\Lambda _{{\goth g}}$ dans $\ec {Gr}g{}{}{\k g{}}$. En outre, sa
restriction \`a $\Omega '_{{\goth g}}$ est l'application 
$(x,y) \mapsto {\goth V}'(x,y)$ d'apr\`es l'assertion (iii) du lemme
\ref{lsc2}; donc $B_{{\goth g}}$ contient $\varphi $ si et seulement si
${\goth V}(x,y)$ contient $\varphi (x,y)$ pour tout $(x,y)$ dans 
$\Lambda _{{\goth g}}$. 

iii) D'apr\`es l'assertion (ii), $B_{{\goth g}}$ contient
$\underline{E}$ et d'apr\`es l'assertion (iii) du corollaire \ref{csc3}, 
$\underline{E}$ est une partie libre de cardinal $\k g{}$ de 
$\tk {{\Bbb C}}{\e Sg}\tk {{\Bbb C}}{\e Sg}{\goth g}$. Soit 
$\poi {\varphi }1{,\ldots,}{\k g{}}{}{}{}$ les \'el\'ements de
$\underline{E}$. Par d\'efinition, pour tout $(x,y)$ dans 
$\Lambda _{{\goth g}}$, 
$\poi {x,y}{}{,\ldots,}{}{\varphi }{1}{\k g{}}$ est une base de
${\goth V}(x,y)$; donc d'apr\`es l'assertion (ii), pour tout $\varphi $ dans
$B_{{\goth g}}$, il existe des fonctions r\'eguli\`eres 
$\poi {\psi }1{,\ldots,}{\k g{}}{}{}{}$ sur $\Lambda _{{\goth g}}$
qui satisfont l'\'egalit\'e
$$ \varphi (x,y) = \psi _{1}(x,y)\varphi _{1}(x,y) + \cdots +
\psi _{\k g{}}(x,y)\varphi _{\k g{}}(x,y) \mbox{ ,}$$
pour tout $(x,y)$ dans $\Lambda _{{\goth g}}$. D'apr\`es le
corollaire \ref{c2sc3}, $\Lambda _{{\goth g}}$ est un grand ouvert de 
${\goth g}\times {\goth g}$; donc les fonctions
r\'eguli\`eres $\poi {\psi }1{,\ldots,}{\k g{}}{}{}{}$ ont un
prolongement r\'egulier \`a ${\goth g}\times {\goth g}$ car 
${\goth g}\times {\goth g}$ est une vari\'et\'e normale. Il en r\'esulte que
$B_{{\goth g}}$ est engendr\'e par 
$\poi {\varphi }1{,\ldots,}{\k g{}}{}{}{}$, d'o\`u l'assertion.

iv) Pour tout $g$ dans $\A g$ et pour tout $(x,y)$ dans 
$\Omega '_{{\goth g}}$, $\Omega '_{{\goth g}}$ contient $(g(x),g(y))$ et
${\goth V}'(g(x),g(y))$ est \'egal \`a $g({\goth V}'(x,y))$; donc
d'apr\`es l'assertion (i), $B_{{\goth g}}$ contient $\varphi $ si et
seulement si $B_{{\goth g}}$ contient l'application
$(x,y) \mapsto g^{-1}\rond \varphi (g(x),g(y))$.

v) Soit $g$ dans GL$_{2}({\Bbb C})$. Pour tout $(x,y)$ dans 
$\Omega '_{{\goth g}}$, $\Omega '_{{\goth g}}$ contient $\kappa_{g}(x,y)$
et ${\goth V}(\kappa _{g}(x,y))$ est \'egal ${\goth V}(x,y)$ car le
sous-espace engendr\'e par $x$ et $y$ ne d\'epend que de l'orbite de
$(x,y)$ sous l'action de GL$_{2}({\Bbb C})$. L'assertion r\'esulte alors de
l'assertion (i).
\end{proof}

\begin{cor}\label{c3sc3}
L'ouvert $\Lambda _{{\goth g}}$ de ${\goth g}\times {\goth g}$ est 
invariant pour les actions de $\A g$ et de ${\rm GL}_{2}({\Bbb C})$ dans
${\goth g}\times {\goth g}$. 
\end{cor}

\begin{proof}
D'apr\`es l'assertion (iii) de la proposition \ref{psc}, pour tout $(x,y)$ 
dans ${\goth g}\times {\goth g}$, ${\goth V}(x,y)$ est l'image de 
$B_{{\goth g}}$ par l'application $\varphi \mapsto \varphi (x,y)$; or 
$\Lambda _{{\goth g}}$ est l'ensemble des \'el\'ements $(x,y)$ de 
${\goth g}\times {\goth g}$ tels que ${\goth V}(x,y)$ soit de dimension
maximale; donc d'apr\`es l'assertion (iv) de la proposition \ref{psc},
pour tout $g$ dans $\A g$, $\Lambda _{{\goth g}}$ contient $(g(x),g(y))$
s'il contient $(x,y)$. De m\^eme, d'apr\`es l'assertion (v) de la proposition
\ref{psc}, pour tout $h$ dans GL$_{2}({\Bbb C})$, $\Lambda _{{\goth g}}$ 
contient $(g(x),g(y))$ s'il contient $(x,y)$. 
\end{proof}
\subsection{} On rappelle que l'ouvert $\Omega '_{{\goth g}}$ de 
${\goth g}\times {\goth g}$ est introduit au lemme \ref{lsc1}.

\begin{lemme}\label{lsc4}
Soit ${\goth p}$ une sous-alg\`ebre parabolique de ${\goth g}$.

{\rm i)} L'intersection de $\Omega '_{{\goth g}}$ et de 
${\goth p}\times {\goth p}$ est non vide.

{\rm ii)} Pour tout $x$ dans l'intersection de ${\goth g}_{\r}$ et de
${\goth p}$, ${\goth p}$ contient ${\goth g}(x)$.

{\rm iii)} Pour tout $(x,y)$ dans ${\goth p}\times {\goth p}$,
${\goth V}(x,y)$ est un sous-espace de ${\goth p}$.
\end{lemme}

\begin{proof}
i) Puisque ${\goth p}$ est une sous-alg\`ebre parabolique de ${\goth g}$,
${\goth p}$ contient une sous-alg\`ebre de Borel 
${\goth b}$ de ${\goth g}$; or $\Omega '_{{\goth g}}$ rencontre 
${\goth b}\times {\goth b}$ d'apr\`es le lemme \ref{lsc3}; donc 
$\Omega '_{{\goth g}}$ rencontre ${\goth p}\times {\goth p}$.

ii)  Puisque ${\goth p}$ contient ${\goth b}$, l'intersection de ${\goth p}$
et de ${\goth g}_{\r}$ est un ouvert non vide de ${\goth p}$; or
l'application $x\mapsto {\goth g}(x)$ est une application r\'eguli\`ere de 
${\goth g}_{\r}$ dans $\ec {Gr}g{}{}{\j g{}}$; donc il suffit de montrer que
${\goth p}$ contient ${\goth g}(x)$ pour tout $x$ dans un ouvert non vide
de ${\goth p}$. Soient ${\goth l}$ un facteur r\'eductif de ${\goth p}$,
${\goth p}_{\u}$ l'ensemble des \'el\'ements nilpotents du radical de 
${\goth p}$ et $P_{\u}$ le sous-groupe ferm\'e, irr\'eductible de $G$ dont
l'alg\`ebre de Lie est l'image de ${\goth p}_{\u}$ par la repr\'esentation
adjointe de ${\goth g}$. Puisque ${\goth l}$ est une sous-alg\`ebre
r\'eductive dans ${\goth g}$ de rang $\j g{}$, ${\goth l}$ contient 
${\goth g}(x)$ pour tout \'el\'ement semi-simple r\'egulier de ${\goth g}$
qui appartient \`a ${\goth l}$; or ${\goth l}$ rencontre la classe
de conjugaison sous $P_{\u}$ de tout \'el\'ement semi-simple de
${\goth g}$ qui appartient \`a ${\goth p}$; donc ${\goth p}$ contient le
centralisateur dans ${\goth g}$ de tout \'el\'ement semi-simple r\'egulier
de ${\goth g}$ qui appartient \`a ${\goth p}$. L'assertion r\'esulte alors de
ce que l'ensemble des \'el\'ements semi-simples r\'eguliers de 
${\goth g}$ est un ouvert qui rencontre ${\goth p}$.

iii) D'apr\`es l'assertion (ii) et l'assertion (i) de la proposition \ref{psc}, 
${\goth p}$ contient ${\goth V}(x,y)$ pour tout $(x,y)$ dans l'inter\-section de 
$\Omega '_{{\goth g}}$ et de ${\goth p}\times {\goth p}$; or d'apr\`es l'assertion (i), 
l'intersection de $\Omega '_{{\goth g}}$ et de ${\goth p}\times {\goth p}$ est partout
dense dans ${\goth p}\times {\goth p}$; donc ${\goth p}$ contient ${\goth V}(x,y)$ pour 
tout $(x,y)$ dans ${\goth p}\times {\goth p}$, d'o\`u l'assertion.
\end{proof}

Soient ${\goth p}$ une sous-alg\`ebre parabolique de ${\goth g}$, 
${\goth p}_{\u}$ le sous-espace des \'el\'ements nilpotents du radical de 
${\goth p}$ et $\varpi $ l'application canonique de ${\goth p}$ sur le quotient 
${\goth l}$ de ${\goth p}$ et de ${\goth p}_{\u}$.

\begin{lemme}\label{l2sc4}
Soit $U$ l'ensemble des \'el\'ements r\'eguliers de ${\goth g}$ qui
appartiennent \`a ${\goth p}$ et dont l'image par $\varpi $ est \'el\'ement
r\'egulier de ${\goth l}$. 

{\rm i)} Pour tout $x$ dans $U$, ${\goth l}(\varpi (x))$ est l'image de 
${\goth g}(x)$ par $\varpi $ si et seulement si l'intersection de 
${\goth g}(x)$ et de ${\goth p}_{\u}$ est nulle.

{\rm ii)} Soit $U'$ l'ensemble des \'el\'ements $x$ de $U$ tels que 
l'intersection de ${\goth g}(x)$ et de ${\goth p}_{\u}$ soit nulle. Alors
$U'$ est un ouvert non vide de ${\goth p}$.

{\rm iii)} Pour tout $(x,y)$ dans ${\goth p}\times U'$, l'image de ${\goth V}(x,y)$
par $\varpi $ est l'image de $B_{{\goth l}}$ par l'application 
$\varphi \mapsto \varphi (\varpi (x),\varpi (y))$.

{\rm iv)} Soient $x$ dans ${\goth p}$ et $y$ dans $U'$. Alors $\Lambda _{{\goth g}}$ 
contient $(x,y)$ si et seulement si $\Lambda _{{\goth l}}$ contient 
$(\varpi (x),\varpi (y))$ et ${\goth V}(x,y)$ contient ${\goth p}_{\u}$.
\end{lemme}
 
\begin{proof}
i) Soit $x$ dans $U$. D'apr\`es l'assertion (ii) du lemme \ref{lsc4}, 
${\goth p}$ contient ${\goth g}(x)$. Puisque $\varpi $ est un morphisme 
d'alg\`ebres de Lie, l'image de ${\goth g}(x)$ par $\varpi $ est contenue dans
${\goth l}(\varpi (x))$; or ${\goth g}$ et ${\goth l}$ sont des alg\`ebres de 
Lie r\'eductives de m\^eme rang; donc ${\goth l}(\varpi (x))$ est l'image de 
${\goth g}(x)$ par $\varpi $ si et seulement si l'intersection de 
${\goth g}(x)$ et de ${\goth p}_{\u}$ est nulle car les dimensions de 
${\goth g}(x)$ et de ${\goth l}(\varpi (x))$ sont \'egales \`a $\j g{}$.

ii) Si $x$ est un \'el\'ement r\'egulier de ${\goth g}$ qui appartient \`a
une sous-alg\`ebre de Cartan de ${\goth g}$, contenue dans ${\goth p}$,
${\goth g}(x)$ est une sous-alg\`ebre de Cartan de ${\goth g}$; donc son
intersection avec ${\goth p}_{\u}$ est nulle. Il en r\'esulte que $U$ et $U'$
contiennent $x$ car $\varpi (x)$ est un \'el\'ement semi-simple r\'egulier de
${\goth l}$. L'application $x\mapsto {\goth g}(x)$ \'etant une application 
r\'eguli\`ere de $U$ dans $\ec {Gr}p{}{}{\j g{}}$, $U'$ est un ouvert de $U$.

iii) Soit $(x,y)$ dans ${\goth p}\times U'$. On note 
${\goth V}_{{\goth l}}(\varpi (x),\varpi (y))$ l'image de $B_{{\goth l}}$ par 
l'application $\varphi \mapsto \varphi (\varpi (x),\varpi (y))$. Pour tout $z$
dans un ouvert non vide du sous-espace engendr\'e par $x$ et $y$, $U'$ 
contient $z$ car $U'$ contient $y$; donc d'apr\`es l'assertion (i), pour tout $z$ dans un
ouvert non vide du sous-espace engendr\'e par $x$ et $y$, ${\goth l}(\varpi (z))$ est 
contenu dans l'image de ${\goth V}(x,y)$ par $\varpi $ et l'image de 
${\goth g}(z)$ par $\varpi $ est contenue dans 
${\goth V}_{{\goth l}}(\varpi (x),\varpi (y))$. Il en r\'esulte que 
pour tout $\varphi $ dans $L_{{\goth l}}$ et pour tout $(a,b)$ dans 
${\Bbb C}^{2}$, $\varphi (a\varpi (x)+b\varpi (y))$ appartient \`a l'image de 
${\goth V}(x,y)$ par $\varpi $. En outre, pour tout $\varphi $ dans $L_{{\goth g}}$ et 
pour tout $(a,b)$ dans ${\Bbb C}^{2}$, $\varpi \rond \varphi (ax+by)$ appartient 
${\goth V}_{{\goth l}}(\varpi (x),\varpi (y))$, d'o\`u l'assertion. 

iv) Soit $(x,y)$ dans ${\goth p}\times U'$. Puisque $\k g{}$ est la somme des entiers 
$\k l{}$ et $\dim {\goth p}_{\u}$, il r\'esulte de l'assertion (iii) que
$\Lambda _{{\goth g}}$ contient $(x,y)$ si et seulement si ${\goth V}(x,y)$ contient 
${\goth p}_{\u}$ et $\Lambda _{{\goth l}}$ contient $(\varpi (x),\varpi (y))$. 
\end{proof}

On didentifie ${\goth l}$ \`a un facteur r\'eductif de ${\goth p}$.

\begin{cor}\label{csc4}
On note $\Sigma _{{\goth g},{\goth p}}$ l'ensemble des \'el\'ements $(x,y)$ de 
${\goth p}\times U'$ qui n'appartien\-nent pas \`a $\Lambda _{{\goth g}}$ et tels que
la dimension de ${\goth V}_{{\goth l}}(\varpi (x),\varpi (y))$ soit \'egale \`a 
$\k l{}$. Soit $(x,y)$ dans $\Sigma _{{\goth g},{\goth p}}$. Alors il existe
un ouvert affine $V$ de ${\goth p}\times {\goth p}$ et un sous-module $L$ de 
$\tk {{\Bbb C}}{{\Bbb C}[V]}{\goth p}$ qui satisfont les conditions suivantes:
\begin{list}{}{}
\item {\rm 1)} $V$ contient $(x,y)$, 
\item {\rm 2)} le module $L$ est engendr\'e par ${\goth p}_{\u}$ et les applications
\sloppy \hbox{$(x',y')\mapsto \varphi (\varpi (x'),\varpi (y'))$} o\`u $\varphi $ est 
dans $B_{{\goth l}}$,
\item {\rm 3)} pour tout $(x',y')$ dans $V$, l'image de $L$ par l'application 
$\varphi \mapsto \varphi (x',y')$ est la somme de ${\goth p}_{\u}$ et de 
${\goth V}_{{\goth l}}(\varpi (x'),\varpi (y'))$.
\end{list}
En outre, $L$ est un module libre de rang $\k g{}$ et pour tout $(x',y')$ dans un ouvert 
non vide de $V$, ${\goth V}(x',y')$ est l'image de $L$ par l'application 
$\varphi \mapsto \varphi (x',y')$.
\end{cor}

\begin{proof}
Puisque la dimension de ${\goth V}_{{\goth l}}(\varpi (x),\varpi (y))$ est \'egale \`a 
$\k l{}$, il existe un ouvert $V'$ de ${\goth l}\times {\goth l}$ qui satisfait les 
conditions suivantes: $V'$ contient $(\varpi (x),\varpi (y))$ et pour tout $(x',y')$ dans
$V'$, ${\goth V}_{{\goth l}}(x',y')$ est de dimension $\k l{}$. Soit $V$ un ouvert affine
de ${\goth p}\times U'$ qui contient $(x,y)$ et dont l'image par l'application 
$(x',y')\mapsto (\varpi (x'),\varpi (y'))$ est contenue dans $V'$. Alors d'apr\`es 
l'assertion (iii) du lemme \ref{l2sc4}, le module $L$ d\'efini par la condition (2)
satisfait la condition (3). En outre, $L$ est un module libre de rang $\k g{}$ car
$B_{{\goth l}}$ est un module libre de rang $\k l{}$ et $\k g{}$ est la somme des entiers
$\k l{}$ et $\dim {\goth p}_{\u}$. Il r\'esulte de l'assertion (iv) du lemme \ref{l2sc4}
que pour tout $(x',y')$ dans l'intersection de $\Lambda _{{\goth g}}$ et de $V$,
${\goth V}(x',y')$ est l'image de $L$ par l'application $\varphi \mapsto \varphi (x',y')$.
\end{proof}
\subsection{} On suppose ${\goth g}$ simple. On rappelle que $(\xi ,\rho ,\eta )$ est un 
${\goth s}{\goth l}_{2}$-triplet principal de ${\goth g}$ et que ${\goth b}$ est
la sous-alg\`ebre de Borel de ${\goth g}$ qui contient $\xi $. On note ${\goth h}$
le centralisateur de $\rho $ dans ${\goth g}$, ${\goth u}$ l'ensemble des \'el\'ements
nilpotents de ${\goth b}$, ${\goth u}_{+,+}$ la somme des noyaux des endomorphismes
$\ad \rho -i$ o\`u $i$ est entier sup\'erieur \`a $4$, ${\goth b}_{-}$ la somme des 
noyaux des endomorphismes $\ad \rho +i$ o\`u $i$ est entier naturel, ${\bf B}$ 
le normalisateur de ${\goth b}$ dans $G$, ${\bf H}$ le centralisateur de ${\goth h}$ dans
${\bf B}$. Alors ${\goth b}_{-}$ est une sous-alg\`ebre de Borel de ${\goth g}$ qui 
contient ${\goth h}$. 

\begin{lemme}\label{lsc5}
Soit $h$ un \'el\'ement r\'egulier de ${\goth g}$ qui appartient \`a ${\goth h}$.

{\rm i)} L'ouvert $\Lambda _{{\goth g}}$ contient $(h,\xi )$.

{\rm ii)} Pour tout $x$ dans ${\goth u}_{+,+}$, $\Lambda _{{\goth g}}$ contient 
$(h,\xi +x)$.

{\rm iii)} Pour tout $x$ dans ${\goth u}$, $\Lambda _{{\goth g}}$ contient $(h+x,\xi )$.

{\rm iv)} L'ouvert $\Lambda _{{\goth g}}$ contient $\{h\}\times {\bf B}.\xi $.

{\rm v)} Pour tout $x$ dans ${\goth b}_{-}$ et pour tout $y$ dans ${\bf H}.\xi $, 
$\Lambda _{{\goth g}}$ contient $(h,x+y)$.
\end{lemme}

\begin{proof}
i) On pose $\xi _{1}=\xi $. Soient $\poi {\xi }2{,\ldots,}{\j g{}}{}{}{}$ des 
\'el\'ements de ${\goth g}(\xi )$ qui ont les propri\'et\'es suivantes:
\begin{list}{}{}
\item 1) la suite $\poi {\xi }1{,\ldots,}{\j g{}}{}{}{}$ est une base de 
${\goth g}(\xi )$,
\item 2) pour $i=1,\ldots,\j g{}$, $\xi _{i}$ est vecteur propre de $\ad \rho $ 
relativement \`a la valeur propre $a_{i}$,
\item 3) la suite $\poi a1{,\ldots,}{\j g{}}{}{}{}$ est croissante.
\end{list}
Selon ces notations, $a_{1}$ est \'egal \`a $2$. Puisque $\xi $ est \'el\'ement 
r\'egulier de ${\goth g}$, ${\goth g}(\xi )$ est l'image de $L_{{\goth g}}$ par 
l'application $\varphi \mapsto \varphi (\xi )$; donc il existe une base 
$\poi {\varphi }1{,\ldots,}{\j g{}}{}{}{}$ de $L_{{\goth g}}$ telle que
$\poi {\xi }1{,\ldots,}{\j g{}}{}{}{}$ soient respectivement \'egaux \`a
$\poi {\xi }{}{,\ldots,}{}{\varphi }{1}{\j g{}}$. Pour $i=1,\ldots,\j g{}$, l'application
$\nu _{i}$, $t\mapsto \varphi _{i}(h+t\xi )$, est une application polynomiale de 
${\Bbb C}$ dans ${\goth g}$ dont le terme constant centralise $h$ et dont
le terme de plus haut degr\'e est \'egal \`a $\xi _{i}$. On note $v_{i,j}$ le terme de 
degr\'e $j$ de $\nu _{i}$. On a alors les \'egalit\'es
$$ [h,v_{i,0}] = 0 \mbox{ , } [h,v_{i,j}] + [\xi ,v_{i,j-1}] = 0 \mbox{ ,} 
\eqno (\star )$$
pour $j$ entier strictement positif, inf\'erieur au degr\'e de $\nu _{i}$. Puisque $h$
est \'el\'ement r\'egulier de ${\goth g}$, il en r\'esulte que $v_{i,j}$ est vecteur 
propre de $\ad \rho $ pour la valeur propre $2j$ et que $a_{i}/2$ est le degr\'e de 
$\nu _{i}$. En outre, ${\goth b}$ contient $v_{i,j}$ pour $i=1,\ldots,\j g{}$ et 
$j=0,\ldots,a_{i}/2$. D'apr\`es l'assertion (iii) de la proposition \ref{psc}, le
sous-espace ${\goth V}(h,\xi )$ est engendr\'e par les coefficients des applications
polynomiales $\poi {\nu }1{,\ldots,}{\j g{}}{}{}{}$; or la dimension de ${\goth g}$ est 
\'egale \`a la somme des entiers $a_{1}+1,\ldots,a_{\j g{}}+1$; donc le nombre de
ces coefficients est \'egal \`a
$$ \j g{} + \frac{\dim {\goth g}-\j g{}}{2} = \k g{} \mbox{ .}$$ 
Il s'agit alors de montrer que ces coefficients sont lin\'eairement ind\'ependants.

Puisque pour $j$ entier naturel, le noyau de $\ad \rho -j$ contient 
$\poi v{1,j}{,\ldots,}{\j g{},j}{}{}{}$, il suffit de montrer que pour $j$ entier 
naturel, les \'el\'ements non nuls de la suite $\poi v{1,j}{,\ldots,}{\j g{},j}{}{}{}$
sont lin\'eairement ind\'ependants car d'apr\`es les \'egalit\'es ($\star $), $v_{i,j}$ 
n'est pas nul si $j$ est inf\'erieur au degr\'e de $\nu _{i}$. On montre cette assertion
en raisonnant par r\'ecurrence sur $a_{\j g{}}/2-j$. Puisque ${\goth g}$ est simple,
$\nu _{\j g{}}$ est la seule application de degr\'e $a_{\j g{}}/2$; donc l'assertion est 
vraie pour $j$ \'egal \`a $a_{\j g{}}/2$ car $\xi _{\j g{}}$ est \'egal \`a 
$v_{\j g{},j}$. On suppose $j$ strictement positif et l'assertion vraie pour $j+1$. 
Soient $\poi b1{,\ldots,}{\j g{}}{}{}{}$ des \'el\'ements de ${\Bbb C}$ tels que
$$b_{1}v_{1,j}+\cdots +a_{\j g{}}v_{\j g{},j} = 0 \mbox{ .}$$
On note $I$ l'ensemble des entiers $i$ de $\{1,\ldots,\j g{}\}$ tels que $j$ soit 
strictement inf\'erieur au degr\'e de $\nu _{i}$. Alors il vient
$$\sum_{i\in I} b_{i}v_{i,j}\in {\goth g}(\xi ) \mbox{ ,}$$
car ${\goth g}(\xi )$ contient $v_{i,j}$ si $j$ n'est pas dans $I$; or ${\goth u}$ 
contient $\poi v{1,j+1}{,\ldots,}{\j g{},j+1}{}{}{}$; donc d'apr\`es les \'egalit\'es 
($\star $), on a
$$\sum_{i\in I} b_{i}v_{i,j+1} = 0 \mbox{ ,}$$
car $h$ est un \'el\'ement r\'egulier de ${\goth g}$. Il r\'esulte alors de l'hypoth\`ese
de r\'ecurrence que $b_{i}$ est nul si $i$ est dans $I$. Par suite, il vient
$$\sum_{i\not \in I} b_{i}v_{i,j} = 0 \mbox{ ;}$$
or pour $i$ dans $\{1,\ldots,\j g{}\}\backslash I$, $v_{i,j}$ est \'egal \`a $\xi _{i}$ 
s'il n'est pas nul; donc $b_{i}$ est nul si $v_{i,j}$ n'est pas nul, d'o\`u l'assertion
pour $j$. Puisque $h$ est \'el\'ement r\'egulier de ${\goth g}$, la suite 
$\poi h{}{,\ldots,}{}{\varphi }{1}{\j g{}}$ est une base de ${\goth h}$ car la suite 
$\poi {\varphi }1{,\ldots,}{\j g{}}{}{}{}$ est une base de $L_{{\goth g}}$; donc 
l'assertion est vraie pour $j$ nul, d'o\`u l'assertion pour tout $j$.

ii) On note $t\mapsto g(t)$ le sous-groupe \`a un param\`etre de ${\bf H}$ dont 
l'alg\`ebre de Lie est engendr\'ee par $\ad \rho $. Soient $x$ dans ${\goth u}_{+,+}$. 
Puisque $\Lambda _{{\goth g}}$ est un ouvert de ${\goth g}\times {\goth g}$ qui 
contient $(h,\xi )$, $\Lambda _{{\goth g}}$ contient $(h,\xi +t^{-2}g(t)(x))$ 
pour tout \'el\'ement non nul $t$ d'un ouvert de ${\Bbb C}$ qui contient $0$; 
donc il existe un \'el\'ement non nul $t$ de ${\Bbb C}$ tel que 
$\Lambda _{{\goth g}}$ contienne $(h,t^{2}\xi +g(t)(x))$ car 
$\Lambda _{{\goth g}}$ est invariant pour l'action de GL$_{2}({\Bbb C})$. 
Puisque $g(t)(h)$ et $g(t)(\xi +h)$ sont respectivement \'egaux \`a $h$ et \`a
$t^{2}\xi +g(t)(h)$, $\Lambda _{{\goth g}}$ contient $(h,\xi +x)$ car 
$\Lambda _{{\goth g}}$ est invariant pour l'action de $G$.

iii) Soient $x$ dans ${\goth u}$. Puisque $x$ est dans ${\goth u}$, $g(t)(x)$ tend vers 
$0$ quand $t$ tend vers $0$; or $\Lambda _{{\goth g}}$ est un ouvert qui contient 
$(h,\xi )$ d'apr\`es l'assertion (i); donc $\Lambda _{{\goth g}}$ contient 
$(h+g(t)(x),\xi )$ pour tout $t$ dans un ouvert non vide de ${\Bbb C}\backslash \{0\}$. 
Puisque $g(t)(h)$ et $g(t)(\xi )$ sont respectivement \'egaux \`a $h$ et \`a $t^{2}\xi $,
$\Lambda _{{\goth g}}$ contient $(h+x,\xi )$ car $\Lambda _{{\goth g}}$ est invariant 
pour les actions de $G$ et de GL$_{2}({\Bbb C})$.

iv) Soit ${\bf U}$ le sous-groupe connexe ferm\'e de $G$ dont l'alg\`ebre de 
Lie est l'image de ${\goth u}$ par la repr\'esentation adjointe de 
${\goth g}$. Alors l'orbite de ${\bf U}.\xi $ est \'egale \`a 
$\xi +{\goth u}_{+,+}$; or ${\bf B}$ est \'egal \`a ${\bf H}{\bf U}$; donc 
${\bf B}.\xi $ est \'egal \`a ${\bf H}.\xi +{\goth u}_{+,+}$. Soient $x$ dans 
${\goth u}_{+,+}$ et $g$ dans ${\bf H}$. Puisque ${\goth u}_{+,+}$ est stable par 
${\bf H}$, $\Lambda _{{\goth g}}$ contient $(h,\xi +g^{-1}(x))$ d'apr\`es l'assertion 
(ii); or $\Lambda _{{\goth g}}$ est invariant pour l'action de $G$; donc 
$\Lambda _{{\goth g}}$ contient $(h,g(\xi )+x)$ car $g(h)$ est \'egal \`a $h$, d'o\`u 
l'assertion. 

vi) Soient $x$ dans ${\goth b}_{-}$ et $k$ dans ${\bf H}$. Puisque $k^{-1}(x)$ est dans 
${\goth b}_{-}$, \sloppy \hbox{$t^{-2}g(t)\rond k^{-1}(x)$} tend vers $0$ quand $t$ tend 
vers l'infini; donc $\Lambda _{{\goth g}}$ contient \sloppy 
\hbox{$(h,\xi +t^{-2}g(t)\rond k^{-1}(x))$} pour tout $t$ dans un ouvert non vide de 
${\Bbb C}$. Puisque $\Lambda _{{\goth g}}$ est invariant pour l'action de 
GL$_{2}({\Bbb C})$, il contient $(h,g(t)(\xi +k^{-1}(x)))$ pour tout $t$ dans un ouvert 
non vide de ${\Bbb C}$; donc $\Lambda _{{\goth g}}$ contient 
$(h,\xi +k^{-1}(x))$ et $(h,k(\xi )+x)$ car $\Lambda _{{\goth g}}$ est 
${\bf H}$-invariant.
\end{proof}

\subsection{} Soit $B$ la base du syst\`eme de racines positives $R_{+}$. Pour tout 
$\beta $ dans $B$, on d\'esigne par ${\goth u}_{\beta }$ la somme des sous-espaces 
${\goth g}^{\alpha }$ o\`u $\alpha $ est dans $R_{+}\backslash \{\beta \}$.

\begin{lemme}\label{lsc6}
Soient $\beta $ dans $B$ et ${\goth b}_{\beta }$ la somme des sous-espaces
${\goth h}$ et ${\goth u}_{\beta }$. On note ${\bf U}_{\beta }$ le sous-groupe
ferm\'e connexe de $G$ dont l'alg\`ebre de Lie est l'image de 
${\goth u}_{\beta }$ par la repr\'esentation adjointe de ${\goth g}$.

{\rm i)} Le sous-espace ${\goth b}_{\beta }$ est une sous-alg\`ebre de 
${\goth b}$. 

{\rm ii)} L'image de ${\bf U}_{\beta }\times {\goth h}$ par l'application
$(g,x)\mapsto g(x)$ est l'ensemble des \'el\'ements semi-simples de ${\goth g}$
qui appartiennent \`a ${\goth b}_{\beta }$.

{\rm iii)} La sous-alg\`ebre ${\goth b}_{\beta }$ contient les composantes 
nilpotente et semi-simple de ses \'el\'ements.

{\rm iv)} Si $x$ est un \'el\'ement r\'egulier de ${\goth g}$ qui 
appartient \`a ${\goth b}_{\beta }$, alors ${\goth b}_{\beta }$ contient 
${\goth g}(x)$.

{\rm v)} Pour tout couple $(x,y)$ d'\'el\'ements de ${\goth b}_{\beta }$, 
${\goth b}_{\beta }$ contient ${\goth V}(x,y)$.
\end{lemme}

\begin{proof}
i) Puisque $\beta $ est dans $B$, ${\goth u}_{\beta }$ est un id\'eal de 
${\goth b}$; donc ${\goth b}_{\beta }$ est une sous-alg\`ebre de 
${\goth b}$. 

ii) Puisque tout \'el\'ement de ${\goth h}$ est \'el\'ement semi-simple de
${\goth g}$, l'image de ${\bf U}_{\beta }\times {\goth h}$ par l'application
$(g,x)\mapsto g(x)$ est contenue dans l'ensemble des \'el\'ements semi-simples
de ${\goth g}$. Soit $x$ un \'el\'ement semi-simple de ${\goth g}$ qui
appartient \`a ${\goth b}_{\beta }$. Puisque ${\goth b}_{\beta }$ est contenu 
dans ${\goth b}$, il existe un \'el\'ement $g$ de ${\bf U}$ tel que 
${\goth h}$ contienne $g(x)$ car ${\goth h}$ est une sous-alg\`ebre de Cartan 
de ${\goth g}$, contenue dans ${\goth b}$. Le sous-groupe ${\bf U}_{\beta }$
de ${\bf U}$ \'etant un sous-groupe invariant de codimension $1$ de ${\bf U}$,
il existe un \'el\'ement $y$ de ${\goth g}^{\beta }$ tel que 
${\bf U}_{\beta }$ contienne $\exp(\ad y)\rond g$. Il en r\'esulte que 
${\goth b}_{\beta }$ contient $\exp(\ad y)\rond g(x)$; donc 
$\exp (\ad y)\rond g(x)$ est \'egal \`a $g(x)$, d'o\`u l'assertion.  

iii) Soient $x$ dans ${\goth b}_{\beta }$, $x_{\n}$ et $x_{\s}$ les 
composantes nilpotente et semi-simple de $x$. Alors $x_{\s}$ est somme
de trois \'el\'ements $x_{1}$, $x_{2}$, $x_{3}$ qui sont respectivement
dans ${\goth h}$, ${\goth g}^{\beta }$, ${\goth u}_{\beta }$. Puisque
$x$ est dans ${\goth b}_{\beta }$, il existe un \'el\'ement $x_{4}$ de
${\goth u}_{\beta }$ tel que $x_{\n}$ soit \'egal \`a $-x_{2}+x_{4}$. Des
\'egalit\'es:
$$0=[x_{\s},x_{\n}] = -[x_{1},x_{2}]+[x_{2},x_{3}]+[x_{1},x_{4}]+[x_{2},x_{4}] 
+[x_{3},x_{4}] \mbox{ ,}$$
on d\'eduit que $[x_{1},x_{2}]$ est nul car ${\goth u}_{\beta }$ est un 
id\'eal de ${\goth b}$. Soit $g$ un \'el\'ement de ${\bf U}$ tel que 
${\goth h}$ contienne $g(x_{\s})$. Alors il existe un \'el\'ement $t$ de 
${\Bbb C}$ tel que ${\bf U}_{\beta }$ contienne $g\rond \exp (t\ad x_{2})$ et 
on a
$$g(x_{\s}) - g(x_{1}+x_{2}) \in {\goth u}_{\beta } \mbox{ , }
g(x_{1}+x_{2}) = g\rond \exp (t\ad x_{2})(x_{1}+x_{2}) \mbox{ , } $$ $$
g\rond \exp (t\ad x_{2})(x_{1}+x_{2}) - x_{1}-x_{2} \in {\goth u}_{\beta } 
\mbox{ ;}$$
donc $x_{2}$ est nul et ${\goth b}_{\beta }$ contient $x_{\s}$ et $x_{\n}$.

iv) Soient $x$ un \'el\'ement r\'egulier de ${\goth g}$ qui appartient \`a
${\goth b}_{\beta }$, $x_{\n}$ et $x_{\s}$ les composantes nilpotente et 
semi-simple de $x$. D'apr\`es l'assertion (iii), ${\goth b}_{\beta }$ contient
$x_{\s}$ et $x_{\n}$. D'apr\`es l'assertion (ii), on peut supposer que 
${\goth h}$ contient $x_{\s}$. Alors le centre ${\goth z}$ de 
${\goth g}(x_{\s})$ est contenu dans ${\goth h}$ car toute sous-alg\`ebre de 
Cartan de ${\goth g}(x_{\s})$ contient ${\goth z}$. Puisque $x$ est un 
\'el\'ement r\'egulier de ${\goth g}$, $x_{\n}$ est un \'el\'ement r\'egulier 
de l'alg\`ebre d\'eriv\'ee de ${\goth g}(x_{\s})$; or l'intersection de 
${\goth b}_{\beta }$ et de l'alg\`ebre d\'eriv\'ee de ${\goth g}(x_{\s})$ 
contient $x_{\n}$; donc ${\goth u}_{\beta }$ contient l'ensemble des 
\'el\'ements nilpotents de la sous-alg\`ebre de Borel de 
${\goth g}(x_{\s})$ qui contient $x_{\n}$. Il en r\'esulte que ${\goth g}(x)$ 
est contenu dans la somme de ${\goth u}_{\beta }$ et de ${\goth z}$. Par 
suite, ${\goth b}_{\beta }$ contient ${\goth g}(x)$.

v) Puisque ${\goth h}$ est contenu dans ${\goth b}_{\beta }$, l'intersection 
de ${\goth g}_{\r}$ et de ${\goth b}_{\beta }$ est un ouvert partout dense de 
${\goth b}_{\beta }$. D'apr\`es l'assertion (iv), pour tout $\varphi $ dans 
$L_{{\goth g}}$ et pour tout $x$ dans l'intersection de ${\goth g}_{\r}$ et de
${\goth b}_{\beta }$, ${\goth b}_{\beta }$ contient $\varphi (x)$; donc 
${\goth b}_{\beta }$ contient $\varphi (x)$ pour tout $x$ dans 
${\goth b}_{\beta }$. Il en r\'esulte que pour tout couple $(x,y)$ 
d'\'el\'ements de ${\goth b}_{\beta }$ et pour tout $\varphi $ dans 
$L_{{\goth g}}$, ${\goth b}_{\beta }$ contient $\varphi (ax+by)$ pour tout 
$(a,b)$ dans ${\Bbb C}^{2}$, d'o\`u l'assertion d'apr\`es l'assertion (iii) de
la proposition \ref{psc}.
\end{proof}

\section{Complexes canoniques d'une alg\`ebre de Lie r\'eductive.} \label{ca} 
Le complexe canonique $C_{\bullet}({\goth g})$ de l'alg\`ebre de Lie 
r\'eductive ${\goth g}$ est d\'efini en \cite{Ch}(D\'efinition 4.1). L'espace 
sous-jacent \`a $C_{\bullet}({\goth g})$ est l'alg\`ebre
$\tk {{\Bbb C}}{\e Sg}\tk {{\Bbb C}}{\e Sg}\ex {}{{\goth g}}$ et sa 
diff\'erentielle est la $\tk {{\Bbb C}}{\e Sg}\e Sg$-d\'erivation qui \`a
l'\'el\'ement $v$ de ${\goth g}$ associe la fonction
$(x,y) \mapsto \dv v{[x,y]}$ sur ${\goth g}\times {\goth g}$. On d\'esigne par
$I_{{\goth g}}$ l'espace des bords de degr\'e $0$ du complexe
$C_{\bullet}({\goth g})$. Selon les notations de \ref{dco2}, 
$C_{\bullet}({\goth g})$ est le complexe $C_{\bullet}(\lambda _{{\goth g}})$ 
o\`u $\lambda _{{\goth g}}$ est l'application qui \`a l'\'el\'ement $v$ de 
${\goth g}$ associe la fonction $(x,y) \mapsto \dv v{[x,y]}$ sur 
${\goth g}\times {\goth g}$. En outre, $I_{{\goth g}}$ est \'egal \`a
$I_{\lambda _{{\goth g}}}$. On d\'esigne par $J_{{\goth g}}$ le radical de
$I_{{\goth g}}$. On rappelle que $K_{\lambda _{{\goth g}}}$ 
d\'esigne le noyau du morphisme $\theta _{\lambda _{{\goth g}}}$. On utilise 
le sous-module $B_{{\goth g}}$ de 
$\tk {{\Bbb C}}{\e Sg}\tk {{\Bbb C}}{\e Sg}{\goth g}$ introduit dans la
proposition \ref{psc}. D'apr\`es la proposition \ref{psc}, $B_{{\goth g}}$
est un module libre de rang $\k g{}$. En outre, $K_{\lambda _{{\goth g}}}$ 
contient $B_{{\goth g}}$ d'apr\`es l'assertion (iii) de la proposition 
\ref{psc} et l'assertion (iv) du lemme \ref{lsc2}.

\begin{Def} \label{dca}
On appelle complexe canonique de deuxi\`eme esp\`ece de l'alg\`ebre de Lie 
${\goth g}$ le complexe canonique associ\'e \`a $\lambda _{{\goth g}}$ et \`a 
$B_{{\goth g}}$ et on le note $E_{\bullet}({\goth g})$. On appelle complexe 
r\'eduit de l'alg\`ebre de Lie ${\goth g}$ et on le note 
$\overline{E}_{\bullet}({\goth g})$ le sous-complexe de
$E_{\bullet}({\goth g})$ dont l'espace des \'el\'ements de degr\'e $j$ est
\'egal \`a $E_{j}({\goth g})$ pour $j$ diff\'erent de $\k g{}$ et \`a
$J_{{\goth g}}E_{\k g{}}({\goth g})$ pour $j$ \'egal \`a $\k g{}$.
\end{Def}

\subsection{} Soient $Y_{+}({\goth g})$ la r\'eunion des supports 
dans ${\goth g}\times {\goth g}$ des groupes d'homologie en degr\'e
diff\'erent de $\k g{}$ du complexe $E_{\bullet}({\goth g})$, 
$Y_{0}({\goth g})$ le support dans 
${\goth g}\times {\goth g}$ du module $J_{{\goth g}}/I_{{\goth g}}$,
et $X_{0}({\goth g})$ l'image de $Y_{0}({\goth g})$ par la premi\`ere 
projection de ${\goth g}\times {\goth g}$ sur ${\goth g}$. Pour
tout $g$ dans $\A g$ et pour tout $h$ dans GL$_{2}({\Bbb C})$, on note
respectivement $g^{\#}$ et $\kappa _{h}^{\#}$ les automorphismes de 
\hbox{$\tk {{\Bbb C}}{\e Sg} \tk {{\Bbb C}}{\e Sg}\ex {}{{\goth g}}$} qui 
\`a l'application $\varphi $ de  ${\goth g}\times {\goth g}$ dans 
$\ex {}{{\goth g}}$ associent les applications 
$(x,y) \mapsto \varphi(g(x),g(y))$ et 
$(x,y) \mapsto \varphi(\kappa _{h}(x,y))$. On rappelle que l'action de 
GL$_{2}({\Bbb C})$ dans ${\goth g}\times {\goth g}$ est d\'efinie au
lemme \ref{lsc1}. Pour $g$ dans $\A g$, on note $\overline{g}$
l'automorphisme lin\'eaire $(x,y)\mapsto (g(x),g(y))$ de 
${\goth g}\times {\goth g}$.

\begin{lemme} \label{lca1}
Soient $g$ un automorphisme de ${\goth g}$ et $h$ un \'el\'ement de 
${\rm GL}_{2}({\Bbb C})$.

{\rm i)} Les id\'eaux $I_{{\goth g}}$ et $J_{{\goth g}}$ sont stables par
$g^{\#}$, $\kappa _{h}^{\#}$.

{\rm ii)} La partie $Y_{0}({\goth g})$ de ${\goth g}\times {\goth g}$ est 
ferm\'ee et stable par $\overline{g}$.

{\rm iii)} La partie $Y_{0}({\goth g})$ de ${\goth g}\times {\goth g}$ est 
stable par $\kappa _{h}$.
\end{lemme}

\begin{proof}
i) Pour tout $v$ dans ${\goth g}$, les images respectives de 
$\lambda _{{\goth g}}(v)$ par $g^{\#}$ et $\kappa _{h}^{\#}$ sont
respectivement \'egales \`a $\lambda _{{\goth g}}(g^{-1}(v))$ et 
$\det h\lambda _{{\goth g}}(v)$; donc $I_{{\goth g}}$ est stable par
$g^{\#}$ et $\kappa _{h}^{\#}$. Il en r\'esulte que 
$J_{{\goth g}}$ est stable par $g^{\#}$ et $\kappa _{h}^{\#}$. 

Les assertions (ii) et (iii) r\'esultent de l'assertion (i) car $J_{{\goth g}}$
est un \sloppy \hbox{$\tk {{\Bbb C}}{\e Sg}\e Sg$-module} de type fini.
\end{proof}

Soient $V_{1}$ et $V_{2}$ deux espaces vectoriels. Pour tout \'el\'ement
non nul $t$ de ${\Bbb C}$, on note $\rho _{t}$ l'automorphisme lin\'eaire
$(x',x)\mapsto (x',tx)$ de $V_{1}\times V_{2}$.

\begin{lemme} \label{l2ca1}
Soient $T$ une partie ferm\'ee de $V_{1}\times V_{2}$ et
$S$ l'image de $T$ par la projection canonique de $V_{1}\times V_{2}$
sur $V_{1}$. On suppose que pour tout nombre complexe non nul $t$,
$T$ est stable par $\rho _{t}$. Soient $T^{*}$ l'intersection de $T$ avec 
$V_{1}\times (V_{2}\backslash \{0\})$ et
$S^{*}$ l'image de $T^{*}$ par la projection de 
$V_{1}\times V_{2}$ sur $V_{1}$. 

{\rm i)} La partie $S^{*}$ de $V_{1}$ est ferm\'ee.

{\rm ii)} La partie $S$ de $V_{1}$ est ferm\'ee.
\end{lemme}

\begin{proof}
i) Soit ${\Bbb P}(V_{2})$ l'espace projectif de $V_{2}$ et $\varpi $ l'application 
canonique de $V_{2}\backslash \{0\}$ sur ${\Bbb P}(V_{2})$. On note $\tilde{T}$ l'image 
de $T^{*}$ par l'application $(x',x)\mapsto (x',\varpi (x))$ de
$V_{1}\times (V_{2}\backslash \{0\})$ dans 
$V_{1}\times {\Bbb P}(V_{2})$. Puisque $T^{*}$ est invariant
par les automorphismes lin\'eaires $\rho _{t}$, $T^{*}$ est
l'image r\'eciproque de $\tilde{T}$ par cette application; or $T^{*}$ est
ferm\'e dans $V_{1}\times (V_{2}\backslash
\{0\})$; donc $\tilde{T}$ est ferm\'e dans $V_{1}\times 
{\Bbb P}(V_{2})$. En outre, $S^{*}$ est la projection de $\tilde{T}$
sur $V_{1}$; donc $S^{*}$ est ferm\'e dans $V_{1}$ car 
${\Bbb P}(V_{2})$ est une vari\'et\'e projective.

ii) Puisque l'image de ${\Bbb C}^{*}$ par l'application $t\mapsto \rho _{t}$ est un
sous-groupe irr\'eductible du groupe des automorphismes lin\'eaires de
$V_{1}\times V_{2}$ qui laissent stable $T$, chacune des composantes irr\'eductibles de 
$T$ est stable par l'action de ce sous-groupe; donc il suffit de montrer l'assertion dans
le cas o\`u $T$ est irr\'eductible. On suppose alors $T$ irr\'eductible. Si $T^{*}$ est 
non vide, $T$ est l'adh\'erence de $T^{*}$ dans $V_{1}\times V_{2}$; donc $S$ est \'egal
\`a $S^{*}$ car $S^{*}$ est ferm\'e d'apr\`es l'assertion (i). Il en r\'esulte
que $S$ est ferm\'e dans $V_{1}$ dans ce cas. Dans le cas contraire, $T$ est \'egal \`a 
$S\times \{0\}$ et $S$ est ferm\'e dans $V_{1}$ car $T$ est ferm\'e dans 
$V_{1}\times V_{2}$.
\end{proof}

\begin{cor} \label{cca1}
La partie $X_{0}({\goth g})$ de ${\goth g}$ est un c\^one ferm\'e, 
$\A g$-invariant.
\end{cor}

\begin{proof}
D'apr\`es l'assertion (iii) du lemme \ref{lca1}, $Y_{0}({\goth g})$ est 
stable par l'automorphisme lin\'eaire $(x,y)\mapsto (x,ty)$ de 
${\goth g}\times {\goth g}$ pour tout \'el\'ement non nul $t$ de ${\Bbb C}$; 
donc le corollaire r\'esulte du lemme \ref{l2ca1} et des assertions (ii) et
(iii) du lemme \ref{lca1}.
\end{proof}\subsection{} On note ${\goth C}_{{\goth g}}$ la vari\'et\'e commutante
de ${\goth g}$. Par d\'efinition, ${\goth C}_{{\goth g}}$ est la vari\'et\'e
des z\'eros dans ${\goth g}\times {\goth g}$ de l'id\'eal $I_{{\goth g}}$.

\begin{lemme}\label{lca2}
L'id\'eal $I_{{\goth g}}$ est premier si et seulement si 
le complexe r\'eduit $\overline{E}_{\bullet}({\goth g})$ n'a pas
d'homologie en degr\'e $\k g{}$. En outre, le $\k g{}$-i\`eme groupe
d'homologie du complexe $E_{\bullet}({\goth g})$ est isomorphe au
quotient de l'alg\`ebre $\tk {{\Bbb C}}{\e Sg}\e Sg$ par l'id\'eal
$I_{{\goth g}}$.
\end{lemme}

\begin{proof}
D'apr\`es l'assertion (iii) de la proposition \ref{psc} et l'assertion (vi) du
lemme \ref{lco1}, $E_{\k g{}}({\goth g})$ est \'egal \`a 
$\ex {\k g{}}{B_{{\goth g}}}$. Il en r\'esulte que $E_{\k g{}}({\goth g})$ est
un module libre de rang $1$. Si $\varepsilon $ est un g\'en\'erateur de
$E_{\k g{}}({\goth g})$, $I_{{\goth g}}\varepsilon $ est l'espace des bords
de degr\'e $\k g{}$ du complexe $E_{\bullet}({\goth g})$; donc les groupes
d'homologie en degr\'e $\k g{}$ des complexes $E_{\bullet}({\goth g})$ et 
$\overline{E}_{\bullet}({\goth g})$ sont respectivement isomorphes \`a
$\tk {{\Bbb C}}{\e Sg}\e Sg/I_{{\goth g}}$ et $J_{{\goth g}}/I_{{\goth g}}$.
Par suite, $\overline{E}_{\bullet}({\goth g})$ n'a pas d'homologie en
degr\'e $\k g{}$ si et seulement si $I_{{\goth g}}$ est un id\'eal radiciel.
D'apr\`es \cite{Ric}, la vari\'et\'e ${\goth C}_{{\goth g}}$ est
irr\'eductible; or ${\goth C}_{{\goth g}}$ est la vari\'et\'e des z\'eros de
$I_{{\goth g}}$ dans ${\goth g}\times {\goth g}$; donc  
$\overline{E}_{\bullet}({\goth g})$ n'a pas d'homologie en degr\'e 
$\k g{}$ si et seulement si $I_{{\goth g}}$ est un id\'eal premier.
\end{proof}\section{Au voisinage d'un \'el\'ement semi-simple.} \label{s}
On fixe un \'el\'ement semi-simple non central $\rho $ de ${\goth g}$.
Soit $\pi _{\rho }$ l'application $(g,x) \mapsto g(x)$ de 
$G\times {\goth g}(\rho )$ dans ${\goth g}$. On a 
alors le lemme bien connu dont la d\'emonstration est rappel\'ee bri\`evement.

\begin{lemme}\label{ls}
L'alg\`ebre de Lie ${\goth g}$ est somme directe des sous-espaces
${\goth g}(\rho )$ et $[\rho ,{\goth g}]$ qui sont orthogonaux. En outre, il
existe un ouvert affine $W$ de ${\goth g}(\rho )$ qui contient $\rho $ et
qui satisfait les conditions suivantes:
\begin{list}{}{}
\item {\rm 1)} la restriction de $\pi _{\rho }$ \`a $G\times W$ est un
morphisme lisse de $G\times W$ sur un ouvert de ${\goth g}$ qui contient
$\rho $,
\item {\rm 2)} pour tout $x$ dans $W$, la restriction de $\ad x$ \`a 
$[\rho ,{\goth g}]$ est un automorphisme lin\'eaire de $[\rho ,{\goth g}]$.
\end{list}
\end{lemme}

\begin{proof}
D'apr\`es la propri\'et\'e d'invariance de la forme de
Killing, les sous-espaces ${\goth g}(\rho )$ et $[\rho ,{\goth g}]$ sont
orthogonaux. Puisque $\rho $ est semi-simple, ${\goth g}$ est
somme directe des sous-espaces propres de $\ad \rho $ et 
$[\rho ,{\goth g}]$ est somme des sous-espaces propres de $\ad \rho $,
relatifs aux valeurs propres non nulles; donc ${\goth g}$ est somme
directe des sous-espaces ${\goth g}(\rho )$ et $[\rho ,{\goth g}]$. 

Pour tout $x$ dans ${\goth g}(\rho )$, le sous-espace $[\rho ,{\goth g}]$ est stable
par $\ad x$. Puisque la restriction de $\ad \rho $ \`a $[\rho ,{\goth g}]$ est bijective,
il existe un ouvert affine $W$ de ${\goth g}(\rho )$ qui contient $\rho $ et qui 
satisfait la condition (2) du lemme. Alors $\pi _{\rho }$ est une submersion en 
$({\rm id}_{{\goth g}},x)$ pour tout $x$ dans $W$. Puisque l'ensemble des points de 
$G\times {\goth g}(\rho )$ en lesquels $\pi _{\rho }$ est une submersion, est stable pour
la multiplication \`a gauche sur $G$, $\pi _{\rho }$ est une submersion en tout point de
$G\times W$. Par suite, la restriction de $\pi _{\rho }$ \`a $G\times W$ est un morphisme
lisse de $G\times W$ sur un ouvert de ${\goth g}$ qui contient $\rho $. 
\end{proof}

Dans ce qui suit, on fixe un ouvert affine $W$ de ${\goth g}(\rho )$ qui 
contient $\rho $ et qui satisfait les conditions (1) et (2) du lemme. L'image
$U$ de $G\times W$ par $\pi _{\rho }$ est alors un ouvert $G$-invariant qui 
contient $\rho $. On rappelle que $X_{0}({\goth g})$ est la
premi\`ere projection sur ${\goth g}$ du support dans 
${\goth g}\times {\goth g}$ du quotient $J_{{\goth g}}/I_{{\goth g}}$.
On choisit $W$ assez petit pour satisfaire la condition suivante: 
$X_{0}({\goth g}(\rho ))$ ne rencontre pas $W$ s'il ne contient pas $\rho $. Ce choix
est possible car $X_{0}({\goth g}(\rho ))$ est ferm\'e dans ${\goth g}(\rho )$ 
d'apr\`es le corollaire \ref{cca1}. Soient $\lambda $
l'application de ${\goth g}$ dans $\tk {{\Bbb C}}{{\Bbb C}[W]}\e Sg$ qui \`a
$v$ associe la fonction \hbox{$(x,y) \mapsto \dv v{[x,y]}$} et 
$\lambda _{\rho }$ la restriction de $\lambda $ \`a ${\goth g}(\rho )$.
On note $I_{W}$ l'id\'eal de  $\tk {{\Bbb C}}{{\Bbb C}[W]}\e Sg$ engendr\'e
par $\lambda ({\goth g})$.

\begin{lemme}\label{l2s}
Si $X_{0}({\goth g}(\rho ))$ ne contient pas $\rho $, alors l'id\'eal $I_{W}$ 
est radiciel.
\end{lemme}

\begin{proof}
On suppose que $X_{0}({\goth g}(\rho ))$ ne contient pas $\rho $. Soit
$I'_{W}$ l'espace des sections au dessus de $W\times {\goth g}(\rho )$ de
la localisation sur ${\goth g}(\rho )\times {\goth g}(\rho )$ de 
$I_{{\goth g}(\rho )}$. Alors $I'_{W}$ est un id\'eal radiciel. D'apr\`es le
lemme \ref{lco3}, le radical $J_{W}$ de
$I_{W}$ est l'ensemble des \'el\'ements $\varphi $ de 
$\tk {{\Bbb C}}{{\Bbb C}[W]}\e Sg$ qui satisfont la condition suivante:
pour tout $x$ dans $W$, $\varphi (x)$ appartient \`a l'id\'eal de $\e Sg$
engendr\'e par l'image de $\ad x$. Soient $\poi v1{,\ldots,}{m}{}{}{}$ une base
de $[\rho ,{\goth g}]$ et $J$ l'id\'eal 
$\tk {{\Bbb C}}{{\Bbb C}[W]}\e Sg[\rho ,{\goth g}]$ de
$\tk {{\Bbb C}}{{\Bbb C}[W]}\e Sg$. D'apr\`es la condition (2) du lemme 
\ref{ls}, pour tout $w$ dans $[\rho ,{\goth g}]$, il existe des \'el\'ements
$\poi a1{,\ldots,}{m}{}{}{}$ de ${\Bbb C}[W]$ tels que
$$ w = \sum_{i=1}^{m} a_{i}(x)[x,v_{i}] \mbox{ ,}$$
pour tout $x$ dans $W$; donc pour tout $(x,y)$ dans $W\times {\goth g}$, on a
$$\dv wy = \sum_{i=1}^{m} a_{i}(x)\dv {[x,v_{i}]}y = 
-\sum_{i=1}^{m} a_{i}(x) \dv {v_{i}}{[x,y]} \mbox{ .}$$
Il en r\'esulte que $I_{W}$ contient $J$; or $\e Sg$ est somme directe de 
$\es S{{\goth g}(\rho )}$ et de l'id\'eal de $\e Sg$ engendr\'e
par  $[\rho ,{\goth g}]$ car ${\goth g}$ est somme directe de 
${\goth g}(\rho )$  et de $[\rho ,{\goth g}]$; donc $J_{W}$ est somme
directe de $J$ et de l'ensemble $J'_{W}$ des \'el\'ements $\varphi $ de 
$\tk {{\Bbb C}}{{\Bbb C}[W]}\es S{{\goth g}(\rho )}$ qui
satisfont la condition suivante: pour tout $x$ dans $W$, $\varphi (x)$
appartient \`a l'id\'eal de $\es S{{\goth g}(\rho )}$ engendr\'e par l'image
de la restriction de $\ad x$ \`a ${\goth g}(\rho )$ car $[\rho ,{\goth g}]$
est contenu dans l'image de $\ad x$. En outre, $I'_{W}$ est contenu dans 
$I_{W}$ car ${\goth g}(\rho )$ et $[\rho ,{\goth g}]$ sont orthogonaux; or d'apr\`es le 
lemme \ref{lco3}, $J'_{W}$ est le radical de $I'_{W}$; donc $I_{W}$ contient $J'_{W}$ et 
$J_{W}$.
\end{proof}

On note ${\cal I}$ et ${\cal J}$ les restrictions \`a $U\times {\goth g}$ des 
localisations respectives sur ${\goth g}\times {\goth g}$ de $I_{{\goth g}}$ 
et $J_{{\goth g}}$, et on pose:
$$I^{(G)} = \Gamma (G\times W\times {\goth g},
(\pi _{\rho }\times {\rm id})^{*}({\cal I})) \mbox{ , } 
J^{(G)} = \Gamma (G\times W\times {\goth g},
(\pi _{\rho }\times {\rm id})^{*}({\cal J})) \mbox{ .}$$
Soit $\alpha $ l'auto\-morphisme de l'alg\`ebre
$\tk {{\Bbb C}}{{\Bbb C}[G]}\tk {{\Bbb C}}{{\Bbb C}[W]}\e Sg$ qui \`a
l'application $\varphi $ de $G\times W$ dans $\e Sg$ associe 
l'application \hbox{$(g,x) \mapsto g.\varphi (g,x)$} o\`u 
\hbox{$(g,p)\mapsto g.p$} d\'esigne le prolongement canonique de l'action
de $G$ dans ${\goth g}$ \`a l'alg\`ebre $\e Sg$.

\begin{lemme}\label{l3s}
L'id\'eal $I^{(G)}$ est l'image par $\alpha $ de l'id\'eal  
$\tk {{\Bbb C}}{{\Bbb C}[G]}I_{W}$.
\end{lemme}

\begin{proof}
Soit $\poi v1{,\ldots,}{n}{}{}{}$ une base de ${\goth g}$. Pour $i=1,\ldots,n$,
il existe des fonctions r\'eguli\`eres $\poi a{i,1}{,\ldots,}{i,n}{}{}{}$ sur 
$G$ qui sont d\'efinies par l'\'egalit\'e:
$$ g(v_{i}) = a_{i,1}(g)v_{1} + \cdots + a_{i,n}(g)v_{n} \mbox{ .}$$
Pour $i=1,\ldots,n$, l'image de $\lambda (v_{i})$ par $\alpha $ est la 
fonction sur $G\times W\times {\goth g}$,
$$ (g,x,y) \mapsto \dv {v_{i}}{[x,g^{-1}(y)]} \mbox{ .}$$
Il r\'esulte alors  des \'egalit\'es:
$$ \dv {v_{i}}{[x,g^{-1}(y)]} = \dv {g(v_{i})}{[g(x),y]} =
\sum_{j=1}^{n} a_{i,j}(g) \dv {v_{j}}{[g(x),y]} \mbox{ ,}$$
que $I^{(G)}$ contient $\alpha \rond \lambda (v_{i})$. R\'eciproquement,
l'image par $\alpha ^{-1}$ de 
$\lambda _{{\goth g}}(v_{i})\rond (\pi _{\rho }\times {\rm id}_{{\goth g}})$
est la fonction sur $G\times W\times {\goth g}$,
$$(g,x,y) \mapsto \dv {v_{i}}{[g(x),g(y)]} \mbox{ .}$$
Il r\'esulte alors de l'\'egalit\'e:
$$\dv {v_{i}}{[g(x),g(y)]} = \sum_{j=1}^{n} a_{i,j}(g^{-1})
\dv {v_{j}}{[x,y]} \mbox{ ,}$$ 
que $\tk {{\Bbb C}}{{\Bbb C}[G]}I_{W}$ contient la fonction
$\alpha ^{-1}(\lambda _{{\goth g}}(v_{i})\rond 
(\pi _{\rho }\times {\rm id}_{{\goth g}}))$, d'o\`u le lemme.
\end{proof}

\begin{prop} \label{ps}
Si $X_{0}({\goth g}(\rho ))$ ne contient pas $\rho $, alors
$X_{0}({\goth g})$ ne contient pas $\rho $.
\end{prop}

\begin{proof}
On suppose que $X_{0}({\goth g}(\rho ))$ ne contient pas $\rho $.
D'apr\`es le lemme \ref{l2s}, $I_{W}$ est un id\'eal
radiciel; donc d'apr\`es le lemme \ref{l3s}, $I^{(G)}$ est
un id\'eal radiciel. D'apr\`es l'assertion (i) du lemme \ref{lj}, $J^{(G)}$ 
est le radical de $I^{(G)}$ car  $G\times W\times {\goth g}$ est une 
vari\'et\'e affine; donc $I^{(G)}$ est \'egal \`a $J^{(G)}$. Il r\'esulte alors
de l'assertion (ii) du lemme \ref{lj} que ${\cal I}$ est \'egal \`a 
${\cal J}$; donc $\{\rho \}\times {\goth g}$ ne rencontre pas le support
de $J_{{\goth g}}/I_{{\goth g}}$, d'o\`u l'assertion. 
\end{proof}

\section{Sur certaines sous-vari\'et\'es de ${\goth g}\times {\goth g}$.}\label{su}
Dans cette section, on suppose ${\goth g}$ simple, de dimension strictement sup\'erieure 
\`a $3$. On rappelle que $(\xi ,\rho ,\eta )$ est un ${\goth s}{\goth l}_{2}$-triplet 
principal de ${\goth g}$, que ${\goth b}$ est l'unique sous-alg\`ebre de Borel de 
${\goth g}$ qui contient $\xi $ et que ${\goth h}$ est le centralisateur de $\rho $ dans 
${\goth g}$. On d\'esigne par ${\goth N}_{{\goth g}}$ le c\^one nilpotent de ${\goth g}$
et par ${\cal X}_{{\goth g}}$ l'ensemble des \'el\'ements $(x,y)$ de 
${\goth g}\times {\goth g}$ tels que $x$ et $y$ appartiennent \`a une m\^eme 
sous-alg\`ebre de Borel de ${\goth g}$. D'apr\`es le lemme \ref{lq2}, 
${\cal X}_{{\goth g}}$ est ferm\'e dans ${\goth g}\times {\goth g}$.

\subsection{} Soit ${\goth s}$ le sous-espace engendr\'e par 
$\xi $, $\rho $, $\eta $. Alors ${\goth s}$ est une sous-alg\`ebre simple de 
${\goth g}$. On note ${\goth N}_{{\goth s}}$ le c\^one nilpotent de 
${\goth s}$, $T_{5}$ l'intersection de ${\cal X}_{{\goth g}}$ avec 
${\goth s}\times {\goth s}$ et ${\bf S}$ le sous-groupe ferm\'e connexe de $G$ dont 
l'alg\`ebre de Lie est l'image de ${\goth s}$ par la repr\'esentation adjointe
de ${\goth s}$.

\begin{lemme}\label{lsu1}
Soient $T_{3}$ et $T_{4}$ les images respectives de 
${\goth N}_{{\goth s}}\times {\Bbb C}$ et de ${\goth s}\times {\Bbb C}$ par 
l'application $(x,t)\mapsto (x,tx)$ de ${\goth s}\times {\Bbb C}$ dans
${\goth s}\times {\goth s}$. 

{\rm i)} Les parties $T_{3}$, $T_{4}$, $T_{5}$ de ${\goth s}\times {\goth s}$
sont ferm\'ees, irr\'eductibles, invariantes pour les actions de ${\bf S}$ et
de ${\rm GL}_{2}({\Bbb C})$. En outre, leurs dimensions sont respectivement
\'egales \`a $3$, $4$, $5$.

{\rm ii)} Soit $T$ une partie ferm\'ee irr\'eductible de 
${\goth s}\times {\goth s}$ qui est non vide et invariante pour les actions de
${\bf S}$ et de ${\rm GL}_{2}({\Bbb C})$. Alors $T$ est l'une des parties suivantes:
$\{0\}$, $T_{3}$, $T_{4}$, $T_{5}$, ${\goth s}\times {\goth s}$. 
\end{lemme}

\begin{proof}
i) Par d\'efinition, $T_{5}$ est \'egal \`a ${\cal X}_{{\goth s}}$; donc
d'apr\`es l'assertion (i) du lemme \ref{lq2}, $T_{5}$ est une partie ferm\'ee 
irr\'eductible de dimension $5$ de ${\goth s}\times {\goth s}$ qui est 
invariante pour les actions de ${\bf S}$ et de GL$_{2}({\Bbb C})$. Par d\'efinition,
$T_{3}$ est l'ensemble des \'el\'ements $(x,y)$ de 
${\goth N}_{{\goth s}}\times {\goth N}_{{\goth s}}$ tels que $x$ et $y$ soient
colin\'eaires; donc $T_{3}$ est ferm\'e dans ${\goth s}\times {\goth s}$ car
${\goth N}_{{\goth s}}$ est ferm\'e dans ${\goth s}$. En outre, $T_{3}$ est
de dimension $3$ et invariant pour l'action de ${\bf S}$ car ${\goth N}_{{\goth s}}$ est 
de dimension $2$ et ${\bf S}$-invariant. Puisque ${\goth N}_{{\goth s}}$ est un c\^one, 
$T_{3}$ est invariant pour l'action de GL$_{2}({\Bbb C})$. De m\^eme, $T_{4}$ est une 
partie ferm\'ee de dimension $4$ de ${\goth s}\times {\goth s}$ qui est invariante pour 
les actions de ${\bf S}$ et de GL$_{2}({\Bbb C})$. Les parties $T_{3}$ et $T_{4}$ sont 
irr\'eductibles comme images de vari\'et\'es irr\'eductibles. 

ii) La partie $T$ \'etant invariante pour l'action de GL$_{2}({\Bbb C})$, 
$T$ est invariant par l'involution $(x,y)\mapsto (y,x)$. En outre, pour tout 
$(x,y)$ dans $T$ et pour toute combinaison lin\'eaire $v$ de $x$ et de $y$, 
$T$ contient $(x,v)$ et $(v,y)$. Il r\'esulte alors du lemme \ref{l2ca1} que l'image 
$T'$ de $T$ par la premi\`ere projection de ${\goth s}\times {\goth s}$ sur ${\goth s}$ 
est un c\^one ferm\'e ${\bf S}$-invariant de ${\goth s}$ car $T$ est non vide et 
invariant pour l'action de ${\bf S}$. Si $T$ est de dimension nulle, $T$ est alors \'egal
\`a $\{0\}$. On suppose $T$ de dimension strictement positive et inf\'erieure \`a
$3$. Alors $T'$ est de dimension inf\'erieure \`a $2$; or les
${\bf S}$-orbites non nulles dans ${\goth s}$ sont de dimension $2$; donc $T'$ est 
\'egal \`a ${\goth N}_{{\goth s}}$ car la ${\bf S}$-orbite d'un \'el\'ement
semi-simple non nul de ${\goth s}$ n'est pas un c\^one. Il en r\'esulte
que $T$ contient $T_{3}$; or $T$ est irr\'eductible et de dimension $3$; donc
$T$ est \'egal \`a $T_{3}$.

On suppose $T$ de dimension sup\'erieure \`a $4$ et $T'$ contenu dans 
${\goth N}_{{\goth s}}$. Il s'agit d'aboutir \`a une contradiction. Alors $T'$
est \'egal \`a ${\goth N}_{{\goth s}}$ car $T'$ n'est pas nul et ${\bf S}$-invariant.
En outre, $T$ \'etant de dimension sup\'erieure \`a $4$, il existe des 
\'el\'ements $a$, $b$, $c$ de ${\Bbb C}$ qui ont les propri\'et\'es suivantes:
$T$ contient $(\xi ,a\xi +b\eta +c\rho )$ et $b$ et $c$ ne sont pas
tous deux nuls; donc $T$ contient $(\xi ,b\eta +c\rho )$. Puisque 
$b\eta +c\rho $ est semi-simple pour $c$ non nul, $c$ est nul et $T$ 
contient $(\xi ,\eta )$. Ceci est absurde car $\xi +\eta $ est semi-simple;
donc $T'$ n'est pas contenu dans ${\goth N}_{{\goth s}}$. Puisque le c\^one engendr\'e 
par la ${\bf S}$-orbite d'un \'el\'ement semi-simple non nul est de dimension $3$, 
$T'$ est \'egal \`a ${\goth s}$ car $T'$ est ferm\'e dans ${\goth s}$. Par suite, $T$ 
contient $T_{4}$. La partie $T$ \'etant irr\'eductible, $T$ est \'egale \`a $T_{4}$ si
$T$ est de dimension $4$.

On suppose $T$ de dimension sup\'erieure \`a $5$ et $T$ diff\'erent de 
$T_{5}$. Puisque $T_{5}$ est ferm\'e dans ${\goth s}\times {\goth s}$, pour 
tout $x$ dans un ouvert non vide de ${\goth s}$, il existe un \'el\'ement $y$ 
de ${\goth s}$ qui a les deux propri\'et\'es suivantes: $T$ contient $(x,y)$ 
et il n'existe pas de sous-alg\`ebre de Borel de ${\goth s}$ qui contient $x$ 
et $y$. La partie $T$ \'etant invariante pour les actions de ${\bf S}$ et de 
GL$_{2}({\Bbb C})$, il en est de m\^eme pour tout $x$ dans un c\^one ouvert 
${\bf S}$-invariant; donc il existe des \'el\'ements non nuls $a$ et $b$ de 
${\Bbb C}$ tels que $T$ contienne $(\rho ,a\xi +b\eta )$. Puisque $T$ est 
invariant pour l'action du sous-groupe \`a un param\`etre de ${\bf S}$ engendr\'e 
par $\ad \rho $, $T$ contient $(\rho ,at\xi +bt^{-1}\eta )$ pour tout 
\'el\'ement non nul $t$ de ${\Bbb C}$; donc $T$ contient 
$(\rho ,sta\xi +bs\eta )$ pour tout couple $(s,t)$ d'\'el\'ements non nuls de 
${\Bbb C}$. Il en r\'esulte que $T$ contient $(\rho ,u\xi +v\eta )$ pour tout
couple $(u,v)$ d'\'el\'ements de ${\Bbb C}$ car $T$ est ferm\'e. 
Par suite, $T$ contient $\{x\}\times {\goth s}$ pour tout $x$ dans le c\^one 
engendr\'e par l'orbite de $\rho $; donc $T$ est \'egal \`a 
${\goth s}\times {\goth s}$, d'o\`u l'assertion.
\end{proof}

\subsection{} La sous-alg\`ebre ${\goth g}$ est somme directe de 
sous-${\goth s}$-modules simples, deux \`a deux orthogonaux, 
$\poi V1{,\ldots,}{\j g{}}{}{}{}$ tels que $V_{1}$ soit \'egal \`a ${\goth s}$. On 
d\'esigne par $V$ la somme des sous-espaces $\poi V2{,\ldots,}{\j g{}}{}{}{}$, par
$\Delta _{V}$ la sous-vari\'et\'e $T_{4}\oplus V\times V$ de ${\goth g}\times {\goth g}$
et par ${\goth g}_{2}$ l'ensemble des \'el\'ements $(x,y)$ de ${\goth g}\times {\goth g}$
tels que $x$ et $y$ soient colin\'eaires. 

\begin{lemme}\label{lsu2}
Soit $X$ une partie ferm\'ee, irr\'eductible de ${\goth g}\times {\goth g}$ qui est 
invariante pour les actions de $G$ et de ${\rm GL}_{2}({\Bbb C})$ dans 
${\goth g}\times {\goth g}$. 

{\rm i)} Soient $x$ dans ${\goth s}$, $v$ et $w$ dans $V$ tels que $x+v$ et $x+w$ soient 
lin\'eairement ind\'ependants. Si $X$ contient $(x+v,x+w)$, alors pour tout 
$(a,b,c,d,r,t)$ dans ${\Bbb C}^{6}$, $X$ contient $(rx+av+bw,tx+cv+dw)$.

{\rm ii)} L'intersection de $X$ et de ${\goth s}\times {\goth s}$  est invariante pour
les actions de ${\bf S}$ et de ${\rm GL}_{2}({\Bbb C})$.

{\rm iii)} Soit $X_{1}$ l'image de $X$ par la premi\`ere projection de 
${\goth g}\times {\goth g}$ sur ${\goth g}$. Alors l'intersection de $X$ et de 
${\goth g}_{2}$ est l'image de $X_{1}\times {\Bbb C}^{2}$ par l'application 
\sloppy \hbox{$(x,s,t)\mapsto (sx,tx)$}.

{\rm iv)} Soit $Y$ l'intersection de $X$ et de $\Delta _{V}$. Alors $X$ contient
$T_{3}$ si l'intersection de $X$ et de ${\goth g}_{2}$ est strictement contenue dans une
composante irr\'eductible de $Y$. 
\end{lemme}

\begin{proof}
i) On suppose que $X$ contient $(x+v,x+w)$. Soient $\Omega $ l'ensemble des \'el\'ements 
$(a,b,c,d)$ de ${\Bbb C}^{4}$ tels que $(a+b)(c+d)$ soit non nul. On note $P_{v,w}$ le 
sous-espace engendr\'e par $v$ et $w$ et $\delta $ l'application
$$ \Omega  \rightarrow P_{v,w}^{2} \mbox{ , }
(a,b,c,d) \mapsto (\frac{av+bw}{a+b},\frac{cv+dw}{c+d})\mbox{ .}$$
Si $v$ et $w$ ne sont pas colin\'eaires, alors $\delta $ est une submersion en tout point
$(a,b,c,d)$ de $\Omega $ tels que $a+b$ et $c+d$ soient diff\'erents de $1$. Dans le cas
o\`u $v$ et $w$ sont colin\'eaires, il existe un \'el\'ement $z$ de ${\Bbb C}$ qui
est diff\'erent de $1$ et qui satisfait une des deux \'egalit\'es $w=zv$ et $v=zw$ car
$x+v$ et $x+w$ sont lin\'eairement ind\'ependants. Alors dans ce cas, $\delta $ est une 
submersion en tout point $(a,b,c,d)$ de $\Omega $ tels que $abcd$ soit non nul. En 
particulier, l'image de $\delta $ est partout dense dans $P_{v,w}^{2}$. Puisque $X$ est 
invariant pour l'action de GL$_{2}({\Bbb C})$ dans ${\goth g}\times {\goth g}$, $X$ 
contient  
$$(x+\frac{av+bw}{a+b},x+\frac{cv+dw}{c+d})\mbox{ ,}$$
pour tout $(a,b,c,d)$ dans $\Omega $; donc $X$ contient $(x+v',x+w')$ pour tout 
$(v',w')$ dans $P_{v,w}^{2}$ car $X$ est ferm\'e ${\goth g}\times {\goth g}$. Par suite, 
pour tout $(v',w')$ dans $P_{v,w}^{2}$ et pour tout couple $(r,t)$ d'\'el\'ements non 
nuls de ${\Bbb C}$, $X$ contient $(x+rv',x+tw')$ et $(r^{-1}x+v',t^{-1}x+w')$, d'o\`u 
l'assertion car $X$ est ferm\'e dans ${\goth g}\times {\goth g}$.  

ii) Puisque $X$ et ${\goth s}\times {\goth s}$ sont invariants pour les actions de 
${\bf S}$ et de GL$_{2}({\Bbb C})$, il en est de m\^eme de leur intersection.

iii) Soit $(x,y)$ dans $X$. Puisque $X$ est invariant pour l'action de GL$_{2}({\Bbb C})$,
$X$ contient $(x,ty)$ pour tout \'el\'ement non nul $t$ de ${\Bbb C}$; donc $X$ contient
$(x,0)$ car $X$ est ferm\'e dans ${\goth g}\times {\goth g}$. Il en r\'esulte que
l'intersection de $X$ et de ${\goth g}_{2}$ contient l'image de $X_{1}\times {\Bbb C}^{2}$
par l'application $(x,s,t)\mapsto (sx,tx)$. R\'eciproquement, si $X$ contient $(sx,tx)$
alors $X_{1}$ contient $sx$ et $tx$ car $X$ est stable par l'involution 
$(x,y)\mapsto (y,x)$; donc $X_{1}$ contient $x$ si $sx$ et $tx$ ne sont pas tous deux 
nuls, d'o\`u l'assertion. 

iv) On suppose que l'intersection de $X$ et de ${\goth g}_{2}$ est strictement contenue 
dans une composante irr\'eductible de $Y$. Soit $Y'$ une composante irr\'eductible
de $Y$ qui contient strictement l'intersection de $X$ et de ${\goth g}_{2}$. Puisque 
$X$ est invariant pour l'action de $G$, $X_{1}$ est $G$-invariant; donc $X_{1}$ n'est pas
contenu dans $V$ car ${\goth g}$ est simple. Alors $Y'$ contient un \'el\'ement qui n'est
pas dans ${\goth g}_{2}$ et dans $V\times V$. Puisque $\Delta _{V}$ et $X$ sont stables 
pour l'action de GL$_{2}({\Bbb C})$, $Y$ et $Y'$ sont stables pour l'action de 
GL$_{2}({\Bbb C})$; donc il existe un \'el\'ement non nul $x$ de ${\goth s}$ et des 
\'el\'ements $v$ et $w$ de $V$ tels que $(x+v,x+w)$ soit un
\'el\'ement de $X$ qui n'appartient pas \`a ${\goth g}_{2}$. Il r\'esulte alors de 
l'assertion (i) que $X$ contient $(x,x)$; donc d'apr\`es l'assertion (ii) et le lemme 
\ref{lsu1}, $X$ contient $T_{3}$.
\end{proof}

On rappelle que ${\goth g}(\rho )$ est la sous-alg\`ebre de Cartan ${\goth h}$ de 
${\goth g}$.

\begin{lemme}\label{l2su2}
Soient $x$ dans ${\goth g}$, $x_{\n}$ et $x_{\s}$ les composantes nilpotente et 
semi-simple de $x$. On note $C$ le c\^one $G$-invariant ferm\'e engendr\'e par $x$.

{\rm i)} Le c\^one $C$ contient $x_{\n}$ et $x_{\s}$.

{\rm ii)} Si $x$ est semi-simple et non nul, alors $C$ contient $\rho +v$ o\`u $v$ est 
un \'el\'ement non nul de l'intersection de ${\goth h}$ et de $V$.

{\rm iii)} Si $x$ est un \'el\'ement r\'egulier de ${\goth g}$, alors $C$ contient
le c\^one nilpotent de ${\goth g}$.
\end{lemme}

\begin{proof}
i) Soient $h$ et $u$ des \'el\'ements de ${\goth g}(x_{\s})$ tels que $(x_{\n},h,u)$ soit
un ${\goth s}{\goth l}_{2}$-triplet de ${\goth g}$. Puisque $C$ est stable par le 
sous-groupe connexe ferm\'e de $G$ dont l'alg\`ebre de Lie est engendr\'ee par $\ad h$, 
$C$ contient $x_{\s}+tx_{\n}$ pour tout $t$ dans ${\Bbb C}\backslash \{0\}$. Il en 
r\'esulte que $C$ contient $x_{\n}+tx_{\s}$ pour tout $t$ dans 
${\Bbb C}\backslash \{0\}$ car $C$ est un c\^one; donc $C$ contient $x_{\s}$ et 
$x_{\n}$ car $C$ est ferm\'e.

ii) On suppose $x$ semi-simple et non nul. Puisque ${\goth h}$ est une sous-alg\`ebre
de Cartan de ${\goth g}$, l'orbite $G.x$ de $x$ rencontre ${\goth h}$. On peut donc 
supposer $x$ contenu dans ${\goth h}$. Soient $W_{{\goth h}}$ le groupe de Weyl de 
${\goth h}$ et $V_{{\goth h}}$ l'intersection de $V$ et de ${\goth h}$. Puisque 
${\goth g}$ est simple, l'action de $W_{{\goth h}}$ dans ${\goth h}$ est irr\'eductible; 
donc l'orbite de $x$ sous l'action de $W_{{\goth h}}$ n'est pas contenue dans 
$V_{{\goth h}}$ et la droite engendr\'ee par $\rho $. Les sous-espaces 
$\poi V1{,\ldots,}{\j g{}}{}{}{}$ \'etant deux \`a deux orthogonaux, $V_{{\goth h}}$ est 
l'orthogonal de $\rho $ dans ${\goth h}$. Il en r\'esulte que l'orbite $G.x$ 
contient la somme d'un \'el\'ement non nul de $V_{{\goth h}}$ et d'un multiple non nul de
$\rho $; donc $C$ contient la somme de $\rho $ et d'un \'el\'ement non nul de 
$V_{{\goth h}}$.

iii) Soient $\poi p1{,\ldots,}{\j g{}}{}{}{}$ des \'el\'ements homog\`enes et 
$G$-invariants de $\e Sg$, de degr\'e respectif $\poi d1{,\ldots,}{\j g{}}{}{}{}$, qui 
engendrent l'alg\`ebre des \'el\'ements $G$-invariants de $\e Sg$. Alors le c\^one
nilpotent de ${\goth g}$ est la vari\'et\'e des z\'eros communs aux fonctions
$\poi p1{,\ldots,}{\j g{}}{}{}{}$. Pour $i=1,\ldots,\j g{}$, on note $a_{i}$ la 
valeur de $p_{i}$ en $x$. Si $x$ est un \'el\'ement nilpotent r\'egulier de 
${\goth g}$, alors l'orbite $G.x$ est partout dense dans le c\^one nilpotent de 
${\goth g}$. On suppose que $x$ est un \'el\'ement r\'egulier, non nilpotent de 
${\goth g}$. Alors il existe un entier $i$ tel que $a_{i}$ soit non nul; donc
l'orbite $G.x$ n'est pas un c\^one et la codimension de $C$ dans ${\goth g}$ est \'egale
\`a $\j g{}-1$. En outre, $C$ est contenu dans la vari\'et\'e des z\'eros communs 
aux fonctions $a_{i}^{d_{j}}p_{j}^{d_{i}}-a_{j}^{d_{i}}p_{i}^{d_{j}}$ o\`u $j$ est dans 
$\{1,\ldots,\j g{}\}$. Il en r\'esulte que la vari\'et\'e des z\'eros de $p_{i}$ dans $C$
est contenue dans le c\^one nilpotent de ${\goth g}$; or la vari\'et\'e des z\'eros de 
$p_{i}$ dans $C$ est de codimension $\j g{}$ dans ${\goth g}$ car $C$ est de codimension 
$\j g{}-1$ dans ${\goth g}$; donc $C$ contient le c\^one nilpotent de ${\goth g}$ car il 
est irr\'eductible et de codimension $\j g{}$ dans ${\goth g}$. 
\end{proof}

\begin{lemme}\label{l3su2}
Soit $X$ une partie ferm\'ee, irr\'eductible de ${\goth g}\times {\goth g}$ qui est 
invariante pour les actions de $G$ et de ${\rm GL}_{2}({\Bbb C})$, $X_{1}$ l'image de 
$X$ par la premi\`ere projection de ${\goth g}\times {\goth g}$ sur ${\goth g}$, 
$Y$ l'intersection de $X$ et de $\Delta _{V}$. On suppose que $X_{1}$ n'est pas contenu 
dans le c\^one nilpotent de ${\goth g}$ et que l'intersection de $X$ et de 
${\goth g}_{2}$ est strictement contenue dans une composante irr\'eductible de $Y$. Alors
$X$ contient $T_{4}$.
\end{lemme}

\begin{proof}
Soit $Y'$ une composante irr\'eductible de $Y$ qui contient l'inter\-section de $X$ et de
${\goth g}_{2}$. Puisque $X$ est ferm\'e et invariant pour les actions de $G$ et de 
GL$_{2}({\Bbb C})$, $X_{1}$ est un c\^one ferm\'e, $G$-invariant de ${\goth g}$ 
d'apr\`es l'assertion (iii) du lemme \ref{lsu2}. Par hypoth\`ese, $X_{1}$ n'est pas 
contenu dans le c\^one nilpotent de ${\goth g}$; donc d'apr\`es les assertions (i) et 
(ii) du lemme \ref{l2su2}, $X_{1}$ contient des \'el\'ements dont la composante sur 
${\goth s}$ est semi-simple et non nulle. Puisque $Y'$ est irr\'eductible et n'est
pas contenu dans ${\goth g}_{2}$, il existe un \'el\'ement $(x,y)$ de $Y'$ qui n'est
pas dans ${\goth g}_{2}$ et tel que la composante $x'$ de $x$ sur ${\goth s}$ soit 
semi-simple non nulle car l'ensemble des \'el\'ements semi-simples non nuls de 
${\goth s}$ est ouvert. Il r\'esulte alors de l'assertion (i) du lemme \ref{lsu2} que 
$X$ contient $(x',x')$. Par suite, l'intersection de $X$ et de 
${\goth s}\times {\goth s}$ n'est pas contenue dans $T_{3}$; donc d'apr\`es le lemme 
\ref{lsu1} et l'assertion (ii) du lemme \ref{lsu2}, $X$ contient $T_{4}$. 
\end{proof}

\subsection{} On note ${\cal N}_{{\goth g}}$ l'ensemble des \'el\'ements $(x,y)$ de
${\goth g}\times {\goth g}$ tels que le sous-espace engendr\'e par $x$ et $y$ soit
contenu dans le c\^one nilpotent de ${\goth g}$.

\begin{lemme}\label{lsu3}
La partie ${\cal N}_{{\goth g}}$ de ${\goth g}\times {\goth g}$ est invariante pour les 
actions de $G$ et de ${\rm GL}_{2}({\Bbb C})$ dans ${\goth g}\times {\goth g}$. En outre,
${\cal N}_{{\goth g}}$ est ferm\'e dans ${\goth g}\times {\goth g}$ et chacune de ses 
composantes irr\'eductibles est de codimension inf\'erieure \`a $\k g{}+\j g{}$.
\end{lemme}

\begin{proof}
Puisque le c\^one nilpotent de ${\goth g}$ est $G$-invariant pour l'action de $G$, 
${\cal N}_{{\goth g}}$ est invariant pour l'action de $G$ dans 
${\goth g}\times {\goth g}$. Pour $(x,y)$ dans ${\goth g}\times {\goth g}$, le sous-espace
de ${\goth g}$ engendr\'e par $x$ et $y$ ne d\'epend que de l'orbite de $(x,y)$ pour 
l'action de GL$_{2}({\Bbb C})$; donc ${\cal N}_{{\goth g}}$ est invariant pour l'action
de GL$_{2}({\Bbb C})$. Soient $\poi p1{,\ldots,}{\j g{}}{}{}{}$ des \'el\'ements 
homog\`enes de $\e Sg$ qui engendrent la sous-alg\`ebre des \'el\'ements $G$-invariants 
de $\e Sg$. Pour $i=1,\ldots,\j g{}$, pour $(a,b)$ dans ${\Bbb C}^{2}$ et pour $(x,y)$ 
dans ${\goth g}\times {\goth g}$, on a
$$p_{i}(ax+by) = \sum_{(m,n)\in {\Bbb N}^{2}} a^{m}b^{n}p_{i,m,n}(ax+by) \mbox{ ,}$$
o\`u $p_{i,m,n}$ sont des \'el\'ements $G$-invariants de $\tk {{\Bbb C}}{\e Sg}\e Sg$
qui sont nuls si $m+n$ est diff\'erent du degr\'e de $p_{i}$. Puisque la somme  
des degr\'es de $\poi p1{,\ldots,}{\j g{}}{}{}{}$ est \'egale \`a $\k g{}$, il y au plus
$\k g{}+\j g{}$ \'el\'ements non nuls parmi les $p_{i,m,n}$. Le c\^one nilpotent de 
${\goth g}$ \'etant la vari\'et\'e des z\'eros communs aux fonctions 
$\poi p1{,\ldots,}{\j g{}}{}{}{}$, ${\cal N}_{{\goth g}}$ est la vari\'et\'e des z\'eros
communs aux fonctions $p_{i,m,n}$ o\`u $i=1,\ldots,\j g{}$ et o\`u $(m,n)$ est dans
${\Bbb N}^{2}$; donc ${\cal N}_{{\goth g}}$ est une partie ferm\'ee de 
${\goth g}\times {\goth g}$ dont les composantes irr\'eductibles sont de codimension 
inf\'erieure \`a $\k g{}+\j g{}$.
\end{proof}

\begin{prop}\label{psu3}
Soit $X$ une partie ferm\'ee, irr\'eductible de ${\goth g}\times {\goth g}$ qui est
invariante pour les actions de $G$ et ${\rm GL}_{2}({\Bbb C})$.  

{\rm i)} Si $X$ est de codimension inf\'erieure \`a $\dim {\goth g}-4$, alors $X$ 
contient $T_{3}$.

{\rm ii)} Si $X$ est de codimension inf\'erieure \`a $\dim {\goth g}-4$ et si $X$ n'est 
pas contenu dans ${\cal N}_{{\goth g}}$, alors $X$ contient $T_{4}$.  
\end{prop}

\begin{proof}
i) On suppose $X$ de codimension inf\'erieure \`a $\dim {\goth g}-4$. Puisque $X$ et 
$\Delta _{V}$ sont des c\^ones, leur intersection $Y$ est un c\^one. Par hypoth\`ese, la 
dimension de $X$ est sup\'erieure \`a $\dim {\goth g}+4$; donc d'apr\`es un th\'eor\`eme 
de Bezout \cite{Ha}(Ch.I, Theorem 7.2), les dimensions des composantes irr\'eductibles de
$Y$ sont sup\'erieures \`a $\dim {\goth g}+2$ car la dimension de $\Delta _{V}$ est 
\'egale \`a $2\dim {\goth g}-2$. D'apr\`es l'assertion (iii) du lemme \ref{lsu2}, 
l'intersection de $X$ et de ${\goth g}_{2}$ est irr\'eductible car $X$ est 
irr\'eductible; or la dimension de ${\goth g}_{2}$ est \'egale \`a $\dim {\goth g}+1$; 
donc il existe une composante irr\'eductible $Y'$ de $Y$ qui contient strictement 
l'intersection de $X$ et de ${\goth g}_{2}$. Il r\'esulte alors de l'assertion (iv) du 
lemme \ref{lsu2} que $X$ contient $T_{3}$. 

ii) On suppose que la codimension de $X$ est inf\'erieure \`a $\dim {\goth g}-4$ et que
$X$ n'est pas contenu dans ${\cal N}_{{\goth g}}$. Soit $X_{1}$ l'image de $X$ par la 
premi\`ere projection de ${\goth g}\times {\goth g}$ sur ${\goth g}$. D'apr\`es 
l'assertion (i), $X_{1}$ contient $\xi $; or $X_{1}$ est $G$-invariant et ferm\'e dans
${\goth g}$ d'apr\`es l'assertion (iii) du lemme \ref{lsu2}; donc $X_{1}$ contient le
c\^one nilpotent de ${\goth g}$. Si $X_{1}$ est contenu dans le c\^one nilpotent de 
${\goth g}$, alors $X$ est contenu dans ${\cal N}_{{\goth g}}$ car $X$ est invariant pour
l'action de GL$_{2}({\Bbb C})$; donc d'apr\`es l'hypoth\`ese, $X_{1}$ n'est pas contenu
dans le c\^one nilpotent de ${\goth g}$. Il r\'esulte alors du lemme \ref{l3su2} que
$X$ contient $T_{4}$ car d'apr\`es (i), l'intersection de $X$ et de ${\goth g}_{2}$
est strictement contenue dans une composante irr\'eductible de l'intersection de 
$X$ et de $\Delta _{V}$. 
\end{proof}

\begin{Def}\label{dsu3}
On dira que l'alg\`ebre de Lie simple ${\goth g}$ a la propri\'et\'e $({\bf N})$ si 
la codimension de chaque composante irr\'eductible de ${\cal N}_{{\goth g}}$ dans 
${\goth g}\times {\goth g}$ est strictement sup\'erieure \`a $\k g{}-\j g{}$.
\end{Def}

\begin{Rem}\label{rsu3}
Soit $X$ une partie ferm\'ee, irr\'eductible de ${\goth g}\times {\goth g}$ qui est 
invariante pour les actions de $G$ et ${\rm GL}_{2}({\Bbb C})$. On suppose que la 
codimension de $X$ est inf\'erieure \`a $\dim {\goth g}-4$ et que $X$ n'est pas contenu
dans ${\cal N}_{{\goth g}}$. Alors $X$ contient $T_{4}$ d'apr\`es l'assertion (ii) de la 
proposition \ref{psu3} car $\dim {\goth g}-4$ est sup\'erieure \`a $\k g{}-\j g{}+1$; or 
l'image de $T_{4}$ par la premi\`ere projection de ${\goth s}\times {\goth s}$ 
sur ${\goth s}$ contient $\rho $; donc l'image de $X$ par la premi\`ere projection de 
${\goth g}\times {\goth g}$ sur ${\goth g}$ contient des \'el\'ements semi-simples 
r\'eguliers de ${\goth g}$.
\end{Rem}

\begin{lemme}\label{l2su3}
On note ${\cal N}_{{\goth g},\xi }$ l'ensemble des \'el\'ements $x$ de ${\goth g}$ tels
que $(x,\xi )$ appartiennent \`a ${\cal N}_{{\goth g}}$. Alors ${\goth g}$ a la 
propri\'et\'e $({\bf N})$ si et seulement si la dimension de chaque composante 
irr\'eductible de ${\cal N}_{{\goth g},\xi }$ est strictement inf\'erieure \`a
$\k g{}+\j g{}$.
\end{lemme}

\begin{proof}
Soit $X$ une composante irr\'eductible de ${\cal N}_{{\goth g}}$. D'apr\`es le lemme
\ref{lsu3}, la codimension de $X$ est inf\'erieure \`a $\k g{}+\j g{}$; donc d'apr\`es 
la remarque \ref{rsu3}, l'image de $X$ par la deuxi\`eme projection de 
${\goth g}\times {\goth g}$ sur ${\goth g}$ est \'egale au c\^one nilpotent de 
${\goth g}$. Soit $X_{\xi }$ l'ensemble des \'el\'ements $x$ de ${\goth g}$ tels que 
$X$ contient $(x,\xi )$. Puisque $X$ est invariant pour l'action de $G$ dans 
${\goth g}\times {\goth g}$, la dimension de $X$ est la somme des dimensions de 
$X_{\xi }$ et du c\^one nilpotent de ${\goth g}$ car l'orbite de $\xi $ sous $G$ est 
ouverte dans le c\^one nilpotent; donc la codimension de $X$ dans 
${\goth g}\times {\goth g}$ est stictement sup\'erieure \`a $\k g{}-\j g{}$ si et 
seulement si la dimension de $X_{\xi }$ est strictement inf\'erieure \`a $\k g{}+\j g{}$,
d'o\`u le lemme.  
\end{proof}
\section{D\'efinition de la propri\'et\'e $({\bf D})$.}\label{q} 
Dans cette section, on suppose que ${\goth g}$ est une alg\`ebre de Lie simple.
On rappelle que $B_{{\goth g}}$ est le sous-module de 
$\tk {{\Bbb C}}{\e Sg}\tk {{\Bbb C}}{\e Sg}{\goth g}$ introduit par la 
proposition \ref{psc} et que $\Lambda _{{\goth g}}$ est l'ensemble des
\'el\'ements $(x,y)$ de ${\goth g}\times {\goth g}$ tels que l'image ${\goth V}(x,y)$ de 
$B_{{\goth g}}$ par l'application $\varphi \mapsto \varphi (x,y)$ soit
de dimension $\k g{}$. On d\'esigne par $X_{{\goth g}}$ le compl\'ementaire de
$\Lambda _{{\goth g}}$ dans ${\goth g}\times {\goth g}$. Le but de la section est la 
d\'efinition de la propri\'et\'e $({\bf D})$ et la proposition \ref{pq} qui est une 
cons\'equence simple de cette propri\'et\'e. On rappelle que $(\xi ,\rho ,\eta )$ est 
un ${\goth s}{\goth l}_{2}$-triplet principal de ${\goth g}$ et que ${\goth b}$ est 
l'unique sous-alg\`ebre de Borel de ${\goth g}$ qui contient $\xi $. 

\subsection{} Soit ${\goth p}$ une sous-alg\`ebre parabolique de ${\goth g}$ qui contient
${\goth b}$. 
On note ${\goth B}_{{\goth g},{\goth p}}$ l'orbite de ${\goth p}$ pour l'action de $G$ 
dans la grassmannienne $\ec {Gr}g{}{}{\dim {\goth p}}$. Alors 
${\goth B}_{{\goth g},{\goth p}}$ est une vari\'et\'e projective. On  d\'esigne par 
${\goth B}_{{\goth g},{\goth p}}^{(1)}$ la vari\'et\'e 
${\goth B}_{{\goth g},{\goth p}}\times {\goth g}\times {\goth g}$ et par
${\goth I}_{{\goth g},{\goth p}}$ la sous-vari\'et\'e des \'el\'ements $(u,x,y)$ de 
${\goth B}_{{\goth g},{\goth p}}^{(1)}$ tels que $u$ contienne $x$ et $y$.

\begin{lemme}\label{lq1}
{\rm i)} La vari\'et\'e ${\goth I}_{{\goth g},{\goth p}}$ est une 
sous-vari\'et\'e lisse de ${\goth B}_{{\goth g}}^{(1)}$.

{\rm ii)} La vari\'et\'e ${\goth I}_{{\goth g},{\goth p}}$ est ferm\'ee dans 
${\goth B}_{{\goth g}}^{(1)}$.

{\rm iii)} La vari\'et\'e ${\goth I}_{{\goth g},{\goth p}}$ est irr\'eductible et sa 
dimension est \'egale \`a $\dim {\goth g}+\dim {\goth p}$.
\end{lemme}
 
\begin{proof}
On note respectivement $\tau $ et $\pi $ les projections canoniques de 
${\goth B}_{{\goth g},{\goth p}}^{(1)}$ sur ${\goth B}_{{\goth g},{\goth p}}$ et 
${\goth g}\times {\goth g}$. Soit $\tau _{1}$ la restriction de $\tau $ \`a 
${\goth I}_{{\goth g},{\goth p}}$. 

i) Soient ${\goth p}'$ dans ${\goth B}_{{\goth g},{\goth p}}$ et $E$ un 
suppl\'ementaire de ${\goth p}'$ dans ${\goth g}$. On note $U_{E}$ l'ensemble
des \'el\'ements $u$ de ${\goth B}_{{\goth g},{\goth p}}$ dont $E$ est 
suppl\'ementaire dans ${\goth g}$. Alors $U_{E}$ est un ouvert affine de 
${\goth B}_{{\goth g},{\goth p}}^{(1)}$. En outre, il existe une 
application r\'eguli\`ere $\psi $ de $U_{E}$ dans 
Hom$_{{\Bbb C}}({\goth p}',E)$ telle que pour tout $u$ dans $U_{E}$, $u$ est 
l'image de l'application lin\'eaire
$$ {\goth p}' \rightarrow {\goth g} \mbox{ , }
x \mapsto x + \psi (u)(x) \mbox{ .}$$  
Alors l'application
$$ U_{E}\times {\goth p}'\times {\goth p}' \rightarrow 
{\goth I}_{{\goth g},{\goth p}} \mbox{ , }
(u,x,y) \mapsto (u,x + \psi (u)(x),y + \psi (u)(y)) 
\mbox{ ,}$$
est un isomorphisme de $U_{E}\times {\goth p}'\times {\goth p}'$ sur 
$\tau _{1}^{-1}(U_{E})$. 
Puisque les ouverts $U_{E}$ recouvrent ${\goth B}_{{\goth g},{\goth p}}$, les 
ouverts $\tau _{1}^{-1}(U_{E})$ recouvrent ${\goth I}_{{\goth g},{\goth p}}$; or 
$U_{E}$ est une vari\'et\'e lisse car ${\goth B}_{{\goth g},{\goth p}}$ est 
une vari\'et\'e lisse; donc ${\goth I}_{{\goth g},{\goth p}}$ est une 
vari\'et\'e lisse. 

ii) Soit $\chi $ l'application
$$ U_{E}\times {\goth p}'\times {\goth p}'\times E\times E \rightarrow 
{\goth B}_{{\goth g},{\goth p}}^{(1)} \mbox{ , }$$ $$
(u,x,y,z_{1},z_{2}) \mapsto 
(u,x + \psi (u)(x)+z_{1},y + \psi (u)(y)+z_{2}) 
\mbox{ .}$$
Alors $\chi $ est un isomorphisme de 
$U_{E}\times {\goth p}'\times {\goth p}'\times E\times E$ sur 
$\tau ^{-1}(U_{E})$; or d'apr\`es (i),  
$\chi (U_{E}\times {\goth p}'\times {\goth p}')$ est \'egal
\`a $\tau _{1}^{-1}(U_{E})$; donc $\tau _{1}^{-1}(U_{E})$ est ferm\'e dans 
$\tau ^{-1}(U_{E})$. Les ouverts $\tau ^{-1}(U_{E})$ recouvrant
${\goth B}_{{\goth g},{\goth p}}^{(1)}$, ${\goth I}_{{\goth g},{\goth p}}$ est
ferm\'e dans ${\goth B}_{{\goth g},{\goth p}}^{(1)}$.

iii) Soit $X$ l'adh\'erence de $\tau _{1}^{-1}(U_{E})$ dans 
${\goth B}_{{\goth g},{\goth p}}^{(1)}$. Alors ${\goth I}_{{\goth g},{\goth p}}$ contient
$X$ d'apr\`es l'assertion (ii). En outre, $X$ est irr\'eductible car 
$\tau _{1}^{-1}(U_{E})$ est irr\'eductible. D'apr\`es (ii), la dimension de 
$\tau _{1}^{-1}(U_{E})$ est \'egale \`a $\dim {\goth g}+\dim {\goth p}$ car 
la dimension de $U_{E}$ est \'egale \`a $\dim {\goth g}-\dim {\goth p}$; donc il suffit 
de montrer que $X$ est \'egal \`a ${\goth I}_{{\goth g},{\goth p}}$. Pour tout $x$ dans 
$U_{E}$, la dimension de l'intersection de $\tau _{1}^{-1}(U_{E})$ et de 
$\tau _{1}^{-1}(x)$ est \'egale \`a $2\dim {\goth p}'$; donc pour tout $x$ dans $X$, la 
dimension de l'intersection de $X$ et de $\tau _{1}^{-1}(x)$ est sup\'erieure \`a 
$2\dim {\goth p}'$. Il en r\'esulte que $X$ contient l'intersection de 
${\goth I}_{{\goth g},{\goth p}}$ et de $\tau _{1}^{-1}(x)$ pour tout $x$ dans 
$\tau _{1}(X)$. En particulier, $X$ est une sous-vari\'et\'e ferm\'ee de 
${\goth B}_{{\goth g},{\goth p}}^{(1)}$ stable par les automorphismes 
$(u,x,y) \mapsto (u,tx,ty)$ o\`u  $t$ est dans ${\Bbb C}\backslash \{0\}$. Il r\'esulte 
alors du lemme \ref{l2ca1} que $\tau _{1}(X)$ est ferm\'e dans 
${\goth B}_{{\goth g},{\goth p}}$. Puisque ${\goth B}_{{\goth g},{\goth p}}$
est irr\'eductible, $\tau _{1}(X)$ est \'egal \`a ${\goth B}_{{\goth g},{\goth p}}$ et 
$X$ est \'egal \`a ${\goth I}_{{\goth g},{\goth p}}$, d'o\`u l'assertion. 
\end{proof}

\begin{lemme}\label{l2q1}
Soit ${\cal X}_{{\goth g},{\goth p}}$ l'ensemble des \'el\'ements $(x,y)$ de
${\goth g}\times {\goth g}$ tels que $x$ et $y$ appartiennent \`a un 
m\^eme \'el\'ement de ${\goth B}_{{\goth g},{\goth p}}$. 

{\rm i)} La sous-alg\`ebre ${\goth p}$ de ${\goth g}$ est l'unique \'el\'ement
de ${\goth B}_{{\goth g},{\goth p}}$ qui contient $\xi $.

{\rm ii)} La partie ${\cal X}_{{\goth g},{\goth p}}$ de 
${\goth g}\times {\goth g}$ est ferm\'ee, irr\'eductible et de codimension \sloppy
\hbox{$\dim {\goth g} - \dim {\goth p}$}.
\end{lemme}

\begin{proof}
i) Soit $g$ un \'el\'ement de $G$ tel que $g({\goth p})$ contienne $\xi $. 
Puisque ${\goth b}$ est l'unique sous-alg\`ebre de Borel de ${\goth g}$ qui contient 
$\xi $, $g({\goth p})$ contient ${\goth b}$. Soient ${\bf P}$ et ${\bf B}$ les 
normalisateurs respectifs de ${\goth p}$ et de ${\goth b}$ dans $G$. Puisque ${\goth b}$ 
et $g^{-1}({\goth b})$ sont deux sous-alg\`ebres de Borel de ${\goth g}$, contenues dans 
${\goth p}$, il existe un \'el\'ement $p$ de ${\bf P}$ tel que $g^{-1}({\goth b})$ soit
\'egal \`a $p({\goth b})$. Alors ${\bf B}$ contient $gp$ et ${\bf P}$ contient $g$ car 
${\bf P}$ contient ${\bf B}$, d'o\`u l'assertion.

ii) Par d\'efinition, ${\cal X}_{{\goth g},{\goth p}}$ est l'image de 
${\goth I}_{{\goth g},{\goth p}}$ par la projection $\pi $; or 
${\goth B}_{{\goth g},{\goth p}}$ est une vari\'et\'e projective; donc 
d'apr\`es les assertions (ii) et (iii) du lemme \ref{lq1}, 
${\cal X}_{{\goth g},{\goth p}}$ est ferm\'e dans ${\goth g}\times {\goth g}$ 
et irr\'eductible. La fibre en $(\rho ,\xi )$
de la restriction de $\pi $ \`a ${\goth I}_{{\goth g},{\goth p}}$ est 
\'egale \`a $\{{\goth p}\}$ d'apr\`es l'assertion (i). En particulier, elle 
est de dimension nulle; donc les vari\'et\'es 
${\goth I}_{{\goth g},{\goth p}}$ et ${\cal X}_{{\goth g},{\goth p}}$
sont de m\^eme dimension car elles sont irr\'eductibles. La dimension de
${\cal X}_{{\goth g},{\goth p}}$ est alors \'egale \`a 
$\dim {\goth p}+\dim {\goth g}$ d'apr\`es l'assertion (iii) du lemme \ref{lq1}.
\end{proof}

\begin{Def}\label{dq1}
On dira que l'alg\`ebre parabolique ${\goth p}$ est admissible si l'alg\`ebre d\'eriv\'ee
de tout facteur r\'eductif de ${\goth p}$ est simple. 
\end{Def}

\subsection{} On  note ${\goth B}_{{\goth g}}$ la sous-vari\'et\'e des 
sous-alg\`ebres de Borel de ${\goth g}$ et ${\goth I}_{{\goth g}}$ la 
vari\'et\'e d'incidence de ${\goth g}$. Par d\'efinition, 
${\goth I}_{{\goth g}}$ est la sous-vari\'et\'e des \'el\'ements 
$(u,x,y)$ de ${\goth B}_{{\goth g}}\times {\goth g}\times {\goth g}$ 
tels que $u$ contienne $x$ et $y$. On d\'esigne par ${\goth B}_{{\goth g}}^{(1)}$ la 
vari\'et\'e ${\goth B}_{{\goth g}}\times {\goth g}\times {\goth g}$ et par $\pi $ la 
projection canonique de ${\goth B}_{{\goth g}}^{(1)}$ sur 
${\goth g}\times {\goth g}$. D'apr\`es le lemme \ref{lq1}, ${\goth I}_{{\goth g}}$
est une sous-vari\'et\'e lisse, ferm\'ee, irr\'eductible de
${\goth B}_{{\goth g}}^{(1)}$.

\begin{lemme}\label{lq2}
Soit ${\cal X}_{{\goth g}}$ l'ensemble des \'el\'ements $(x,y)$ de
${\goth g}\times {\goth g}$ tels que $x$ et $y$ appartiennent \`a une 
m\^eme sous-alg\`ebre de Borel de ${\goth g}$.

{\rm i)} La partie ${\cal X}_{{\goth g}}$ de ${\goth g}\times {\goth g}$ est
ferm\'ee, irr\'eductible et de codimension $\k g{}-\j g{}$.

{\rm ii)} L'ouvert $\Lambda _{{\goth g}}$ rencontre ${\cal X}_{{\goth g}}$.

{\rm iii)} L'intersection de $X_{{\goth g}}$ et de ${\cal X}_{{\goth g}}$
est de codimension $1$ dans ${\cal X}_{{\goth g}}$. 
\end{lemme}

\begin{proof}
i) L'assertion r\'esulte de l'assertion (ii) du lemme \ref{l2q1} et des
\'egalit\'es:
$$\dim {\goth g} = 2\k g{}-\j g{} \mbox{ , } 
\dim {\goth g}-\dim {\goth b} = \k g{}-\j g{} \mbox{ .}$$ 

ii) D'apr\`es les assertions (i) et (iii) du lemme \ref{lsc3}, 
$\Lambda _{{\goth g}}$ contient $(\rho ,\xi )$; donc ${\cal X}_{{\goth g}}$
rencontre $\Lambda _{{\goth g}}$.

iii) D'apr\`es l'assertion (ii), l'intersection de $X_{{\goth g}}$ et de
${\cal X}_{{\goth g}}$ est strictement contenue dans ${\cal X}_{{\goth g}}$
Soient $\poi {\varphi }1{,\ldots,}{\k g{}}{}{}{}$ une base de $B_{{\goth g}}$ et 
$\omega $ un g\'en\'erateur non nul de $\ex {\k g{}}{{\goth b}}$. D'apr\`es
l'assertion (iii) du lemme \ref{lsc4}, pour tout $(x,y)$ dans 
${\goth b}\times {\goth b}$, ${\goth b}$ contient $\varphi (x,y)$ pour tout
$\varphi $ dans $B_{{\goth g}}$; donc il existe un \'el\'ement non nul $p$
de $\tk {{\Bbb C}}{\es S{{\goth b}^{*}}}\es S{{\goth b}^{*}}$ tel que l'on ait
$$ \poi {x,y}{}{\wedge \cdots \wedge }{}{\varphi }{1}{\k g{}} = 
p(x,y)\omega \mbox{ ,}$$
pour tout $(x,y)$ dans ${\goth b}\times {\goth b}$. Il en r\'esulte que 
l'intersection de ${\goth b}\times {\goth b}$ et de $X_{{\goth g}}$ est une
hypersurface de ${\goth b}\times {\goth b}$. Par suite, l'intersection de 
${\goth I}_{{\goth g}}$ et de $\pi ^{-1}(X_{{\goth g}})$ est une hypersurface
de ${\goth I}_{{\goth g}}$ dont les composantes irr\'eductibles sont toutes
de codimension $1$. Soit $Z$ une composante irr\'eductible de l'intersection de
${\goth I}_{{\goth g}}$ et de $\pi ^{-1}(X_{{\goth g}})$ telle que $\pi (Z)$
soit une composante irr\'eductible de l'intersection de $X_{{\goth g}}$ et de
${\cal X}_{{\goth g}}$ qui contient $(\xi ,\xi )$. Puisque ${\goth b}$ est la seule 
sous-alg\`ebre de Borel de ${\goth g}$ qui contient $\xi $, la fibre en $(\xi ,\xi )$ 
de la restriction de $\pi $ \`a $Z$ est de dimension nulle; donc $Z$ et $\pi (Z)$
ont m\^eme dimension car $Z$ est irr\'eductible. Il en r\'esulte que
$\pi (Z)$ est de codimension $1$ dans ${\cal X}_{{\goth g}}$ car
${\cal X}_{{\goth g}}$ et ${\goth I}_{{\goth g}}$ ont m\^eme dimension. Par 
suite, l'intersection de $X_{{\goth g}}$ et de ${\cal X}_{{\goth g}}$ est de
codimension $1$ dans ${\cal X}_{{\goth g}}$.
\end{proof}

On rappelle que pour ${\goth p}$ sous-alg\`ebre parabolique de ${\goth g}$, la partie 
$\Sigma _{{\goth g},{\goth p}}$ de ${\goth p}\times {\goth p}$ est d\'efinie au 
corollaire \ref{csc4}. On note $\Sigma '_{{\goth g},{\goth p}}$ l'ensemble des 
\'el\'ements de $\Sigma _{{\goth g},{\goth p}}$ dont la deuxi\`eme composante est 
\'el\'ement semi-simple r\'egulier de ${\goth g}$. Pour tout \'el\'ement $x$ de 
${\goth g}$ qui est semi-simple et r\'egulier, on note $G'(x)$ le sous-groupe des
automorphismes lin\'eaires de ${\goth g}\times {\goth g}$ qui est engendr\'e par les 
automorphismes $(v,w)\mapsto (tv,w)$ et $(v,w)\mapsto (g(v),g(w))$ o\`u $t$ et $g$ sont
respectivement dans ${\Bbb C}\backslash \{0\}$ et dans le centralisateur de $x$ dans 
$G$.

\begin{Def}\label{dq}
On d\'esigne par $\Sigma '_{{\goth g}}$ la plus petite partie de 
${\goth g}\times {\goth g}$ qui est invariante pour les actions de $G$ et de 
${\rm GL}_{2}({\Bbb C})$ et qui contient la r\'eunion des 
$\Sigma '_{{\goth g},{\goth p}}$ o\`u ${\goth p}$ est sous-alg\`ebre parabolique 
admissible de ${\goth g}$ au sens de {\rm \ref{dq1}}. Soit $\Sigma _{{\goth g}}$
l'ensemble des \'el\'ements $(x,y)$ de ${\goth g}\times {\goth g}$ qui ont les deux
propri\'et\'es suivantes:
\begin{list}{}{}
\item 1) $y$ est un \'el\'ement semi-simple r\'egulier de ${\goth g}$,
\item 2) l'adh\'erence dans ${\goth g}\times {\goth g}$ de l'orbite de $(x,y)$, sous 
l'action de $G'(x)$, rencontre $\Sigma '_{{\goth g}}$.
\end{list}
On dira que l'alg\`ebre de Lie ${\goth g}$ a la propri\'et\'e $({\bf D})$ si toute partie
ferm\'ee irr\'eductible de $X_{{\goth g}}$, invariante pour les actions de $G$ et de 
${\rm GL}_{2}({\Bbb C})$, de codimension inf\'erieure \`a $\dim {\goth g}-4$, qui ne 
rencontre pas $\Sigma _{{\goth g}}$, est contenue dans la r\'eunion de 
${\cal X}_{{\goth g}}$ et de ${\cal N}_{{\goth g}}$.
\end{Def}

Si ${\goth g}$ est de rang $1$, $X_{{\goth g}}$ est contenu dans ${\cal X}_{{\goth g}}$;
donc d'apr\`es les assertions (i) et (iii) du lemme \ref{lq2}, ${\goth g}$ a la 
propri\'et\'e $({\bf D})$ dans ce cas.

\begin{prop}\label{pq}
On suppose que ${\goth g}$ a la propri\'et\'e $({\bf D})$ au sens de {\rm \ref{dq}} et la
propri\'et\'e $({\bf N})$ au sens de {\rm \ref{dsu3}}. Soit $X$ une partie
ferm\'ee de $X_{{\goth g}}$ qui est irr\'eductible et invariante pour les actions de
$G$ et de ${\rm GL}_{2}({\Bbb C})$. Si $X$ ne rencontre pas $\Sigma _{{\goth g}}$, alors
la codimension de $X$ dans ${\goth g}\times {\goth g}$ est sup\'erieure \`a
$\k g{}-\j g{}+1$.
\end{prop}

\begin{proof}
On suppose que $X$ ne rencontre pas $\Sigma _{{\goth g}}$ et que la codimension de $X$
est inf\'erieure \`a $\dim {\goth g}-4$. Alors d'apr\`es la propri\'et\'e $({\bf D})$,
$X$ est contenu dans la r\'eunion de ${\cal X}_{{\goth g}}$ et de ${\cal N}_{{\goth g}}$;
or d'apr\`es les assertions (i) et (iii) du lemme \ref{lq2}, la codimension dans 
${\goth g}\times {\goth g}$ de toute composante irr\'eductible de l'intersection de 
$X_{{\goth g}}$ et de ${\cal X}_{{\goth g}}$ est sup\'erieure \`a $\k g{}-\j g{}+1$; donc
la codimension de $X$ dans ${\goth g}\times {\goth g}$ est sup\'erieure \`a 
$\k g{}-\j g{}+1$ car la codimension de chaque composante irr\'eductible de 
${\cal N}_{{\goth g}}$ est sup\'erieure \`a $\k g{}-\j g{}+1$, d'apr\`es la 
propri\'et\'e $({\bf N})$.
\end{proof}

\begin{Rem}\label{rq}
Il r\'esultera du th\'eor\`eme \ref{ta} qu'une alg\`ebre de Lie simple a la propri\'et\'e
$({\bf D})$ au sens de \ref{dq}.
\end{Rem}

\subsection{} Soit ${\goth p}$ une sous-alg\`ebre parabolique de ${\goth g}$ qui contient
${\goth b}$. On note ${\goth p}_{\u}$ le sous-espaces des \'el\'ements nilpotents du
radical de ${\goth p}$ et ${\bf P}$ le normalisateur de ${\goth p}$ dans $G$. On rappelle
que ${\cal X}_{{\goth g},{\goth p}}$ d\'esigne l'image de 
$G\times {\goth p}\times {\goth p}$ par l'application $\Psi $, 
$(g,x,y)\mapsto (g(x),g(y))$ de $G\times {\goth g}\times {\goth g}$ dans 
${\goth g}\times {\goth g}$.

\begin{lemme}\label{lq3}
Soit $X$ une partie constructible, irr\'eductible de ${\cal X}_{{\goth g},{\goth p}}$ qui
est invariante pour l'action de $G$ dans ${\goth g}\times {\goth g}$.

{\rm i)} Il existe une composante irr\'eductible $Y$ de l'intersection de $X$ et de 
${\goth p}\times {\goth p}$ telle que $\Psi (G\times Y)$ soit \'egal \`a $X$.

{\rm ii)} L'intersection de $X$ et de ${\goth g}\times \{\xi \}$ est contenue dans $Y$. 

{\rm iii)} Si l'intersection de $X$ et de ${\goth g}\times \{\xi \}$ est non vide, alors
la dimension de $Y$ est \'egale \`a $\dim X-\dim {\goth p}_{\u}$.

{\rm iv)} Si l'image de $X$ par la deuxi\`eme projection de ${\goth g}\times {\goth g}$
sur ${\goth g}$ contient un \'el\'ement semi-simple r\'egulier de ${\goth g}$, alors 
l'image de $Y$ par la deuxi\`eme projection de ${\goth p}\times {\goth p}$ sur 
${\goth p}$ contient un \'el\'ement de ${\goth h}$ qui est \'el\'ement r\'egulier de 
${\goth g}$.
\end{lemme}

\begin{proof}
i) Soient $\poi Y1{,\ldots,}{p}{}{}{}$ les composantes irr\'eductibles de 
\sloppy \hbox{l'intersection} de $X$ et de ${\goth p}\times {\goth p}$. 
Alors $\poi Y1{,\ldots,}{p}{}{}{}$ sont invariants par ${\bf P}$ car ${\bf P}$ est 
irr\'eductible. Puisque ${\goth p}$ est une sous-alg\`ebre parabolique, la vari\'et\'e 
homog\`ene $G/{\bf P}$ est projective; donc $\poi {G\times Y}1{,\ldots,}{p}{\Psi }{}{}$ 
sont ferm\'es dans $X$. En outre, $X$ est la r\'eunion
de $\poi {G\times Y}1{,\ldots,}{p}{\Psi }{}{}$ car $X$ est contenu dans 
${\cal X}_{{\goth g},{\goth p}}$; donc il existe un entier $j$ tel que 
$\Psi (G\times Y_{j})$ soit \'egal \`a $X$ car $X$ est 
irr\'eductible, d'o\`u l'assertion. 

ii) Soit $x$ un \'el\'ement de ${\goth g}$ tel que $X$ contienne $(x,\xi )$. D'apr\`es 
l'assertion (i), il existe un \'el\'ement $g$ de $G$ tel que $Y$ contienne 
$(g(x),g(\xi ))$. Puisque ${\goth p}$ contient $g(\xi )$, il existe un \'el\'ement $h$ de
${\bf P}$ tel que ${\goth b}$ contienne $hg(\xi )$ car ${\goth p}$ contient ${\goth b}$. 
Par suite, il existe un \'el\'ement $k$ de ${\bf B}$ tel que $khg(\xi )$ soit \'egal \`a 
$\xi $ car l'ensemble des \'el\'ements nilpotents r\'eguliers de ${\goth g}$ qui 
appartiennent \`a ${\goth b}$ sont contenus dans une ${\bf B}$-orbite. Puisque $Y$ est 
invariant pour l'action de ${\bf P}$, $Y$ contient $(khg(x),khg(\xi ))$ et $(x,\xi )$ car
${\bf B}$ contient le centralisateur de $\xi $ dans $G$, d'o\`u l'assertion.

iii) On suppose que l'intersection de $X$ et de ${\goth g}\times \{\xi \}$ est non vide. 
Pour $x$ dans ${\goth p}$ et $g$ dans $G$ tels que $X$ contiennent $(x,\xi )$
et $(g(x),g(\xi ))$, ${\bf P}$ contient $g$ d'apr\`es (ii); donc la dimension de 
$X$ est la somme des dimensions de $Y$ et de ${\goth p}_{\u}$ car 
$\dim G-\dim {\bf P}$ est \'egal \`a $\dim {\goth p}_{\u}$.

iv) On suppose que l'image de $X$ par la deuxi\`eme projection de 
${\goth g}\times {\goth g}$ sur ${\goth g}$ contient des \'el\'ements semi-simples 
r\'eguliers de ${\goth g}$. On note $Y'$ l'image de $Y$ par la deuxi\`eme projection de 
${\goth p}\times {\goth p}$ sur ${\goth p}$. Alors $G(Y')$ est l'image de $X$ par la 
deuxi\`eme projection de ${\goth g}\times {\goth g}$ sur ${\goth g}$; donc $Y'$ contient 
des \'el\'ements semi-simples r\'eguliers de ${\goth g}$, d'o\`u l'assertion car $Y'$ est
${\bf P}$-invariant.
\end{proof}

\begin{Def}\label{dq3}
Soit $X$ une partie de ${\goth g}\times {\goth g}$. On dira que $X$ a la propri\'et\'e 
$({\bf Q})$ s'il satisfait les conditions suivantes:
\begin{list}{}{}
\item 1) $X$ est constructible et irr\'eductible,
\item 2) $X$ est invariant pour l'action de $G$ dans ${\goth g}\times {\goth g}$,
\item 3) pour tout $(s,t)$ dans ${\Bbb C}^{2}$, $X$ contient $(sx,ty)$ s'il contient
$(x,y)$,
\item 4) l'image de $X$ par la deuxi\`eme projection de ${\goth g}\times {\goth g}$ sur 
${\goth g}$ contient un \'el\'ement semi-simple r\'egulier de ${\goth g}$.
\end{list}
 
On dira que ${\goth g}$ a la propri\'et\'e $({\bf D}')$ si ${\cal X}_{{\goth g}}$ contient
les parties de $X_{{\goth g}}$ qui ont la propri\'et\'e $({\bf Q})$ et qui ne rencontrent
pas $\Sigma _{{\goth g}}$.
\end{Def}

\begin{prop}\label{pq3}
On suppose que les facteurs simples des centralisateurs des \'el\'ements semi-simples
non nuls de ${\goth g}$ ont la propri\'et\'e $({\bf D}')$ au sens de {\rm \ref{dq3}}.
Soit $X$ une partie de $X_{{\goth g}}$ qui a la propri\'et\'e $({\bf Q})$ au sens de
{\rm \ref{dq3}} et qui ne rencontre pas $\Sigma _{{\goth g}}$. Alors $X$ satisfait l'une 
ou l'autre des deux conditions suivantes:
\begin{list}{}{}
\item {\rm 1)} pour toute sous-alg\`ebre parabolique admissible de ${\goth g}$ au sens
de {\rm \ref{dq1}}, ${\cal X}_{{\goth g},{\goth p}}$ ne contient pas $X$,
\item {\rm 2)} ${\cal X}_{{\goth g}}$ contient $X$.
\end{list}
\end{prop}

\begin{proof}
On suppose que $X$ ne satisfait pas la condition (1). Alors il existe une 
sous-alg\`ebre parabolique admissible ${\goth p}$, qui contient ${\goth b}$, telle que 
${\cal X}_{{\goth g},{\goth p}}$ contienne $X$. D'apr\`es l'assertion (i) du lemme 
\ref{lq3}, il existe une composante irr\'eductible $Y$ de l'intersection de $X$ et de 
${\goth p}\times {\goth p}$ telle que $X$ soit \'egal \`a $\Psi (G\times Y)$. Soient 
${\goth l}$ un facteur r\'eductif de ${\goth p}$, ${\goth z}$ le centre de ${\goth l}$, 
${\goth l}'$ l'alg\`ebre d\'eriv\'ee de ${\goth l}$, $\varpi $ la projection canonique de
${\goth p}$ sur ${\goth l}$. On note $Y'$ l'image de $Y$ par l'application 
$\varpi \times \varpi $ et $Y''$ l'image de $Y'$ par la projection canonique de 
${\goth l}\times {\goth l}$ sur ${\goth l}'\times {\goth l'}$. Puisque ${\goth p}$ est 
une sous-alg\`ebre parabolique admissible de ${\goth g}$, ${\goth l}'$ est une alg\`ebre 
simple. En outre, ${\goth l}$ est le centralisateur dans ${\goth g}$ d'un \'el\'ement 
semi-simple non nul de ${\goth g}$; donc ${\goth l}'$ a la propri\'et\'e $({\bf D}')$. Si
$Y''$ est contenu dans ${\cal X}_{{\goth l}'}$, alors $Y$ est contenu dans 
${\cal X}_{{\goth g}}$ car pour ${\goth b}'$ sous-alg\`ebre de Borel de ${\goth l}'$, 
$\varpi ^{-1}({\goth b}'+{\goth z})$ est une sous-alg\`ebre de Borel de ${\goth g}$; donc
il suffit de montrer que $Y''$ est contenue dans l'adh\'erence d'une partie de 
${\goth l}'\times {\goth l}'$ qui a la propri\'et\'e $({\bf Q})$ et qui ne rencontre pas 
$\Sigma _{{\goth l}'}$.  

Soient $\pi $ la projection canonique de ${\goth l}$ sur ${\goth l}'$ et $V$ l'image 
par $\pi $ de l'ensemble des \'el\'ements semi-simples r\'eguliers de ${\goth g}$ qui
appartiennent \`a ${\goth l}$. Alors $V$ est un c\^one ouvert non vide de ${\goth l}'$, 
invariant par le groupe adjoint de ${\goth l}'$. On note $Y_{1}$ l'intersection de $Y''$ 
et de ${\goth l}'\times V$. D'apr\`es l'assertion (iv) du lemme \ref{lq3}, l'image de 
$Y''$ par la deuxi\`eme projection de ${\goth l}'\times {\goth l}'$ sur ${\goth l}'$ 
rencontre $V$; donc $Y_{1}$ est une partie constructible, irr\'eductible de
${\goth l}'\times \times {\goth l}'$ qui est partout dense dans $Y''$ car $Y''$
est une partie constructible, irr\'eductible. Puisque $X$ a la propri\'et\'e $({\bf Q})$,
$Y''$ est invariant pour l'action diagonale du groupe adjoint de ${\goth l}'$ dans 
${\goth l}'\times {\goth l}'$. En outre, $Y''$ satisfait la condition (2) de la 
d\'efinition \ref{dq3}; donc $Y_{1}$ a la propri\'et\'e $({\bf Q})$ car $V$ est un c\^one
ouvert, invariant pour l'action du groupe adjoint de ${\goth l}'$. Si ${\goth q}$ est une
sous-alg\`ebre parabolique admissible de ${\goth l}'$, alors 
$\varpi ^{-1}({\goth q}+{\goth z})$ est une sous-alg\`ebre parabolique admissible de 
${\goth g}$; donc $\Sigma '_{{\goth l}'}$ est contenu dans $\Sigma '_{{\goth g}}$
d'apr\`es l'assertion (iv) du lemme \ref{lsc4}. Soient $y$ un \'el\'ement semi-simple 
r\'egulier de ${\goth g}$ qui appartient \`a ${\goth l}$ et $x$ un \'el\'ement de 
${\goth l}'$ tels que $Y_{1}$ contienne $(x,\pi (y))$. On note $A$ le sous-groupe des 
automorphismes lin\'eaires de ${\goth l}\times {\goth l}$ qui est engendr\'e par les 
automorphismes $(v,w)\mapsto (tv,w)$ et $(v,w)\mapsto (g(v),g(w))$ o\`u $t$ et $g$ sont 
respectivement dans ${\Bbb C}\backslash \{0\}$ et le centralisateur de $y$ dans le groupe
adjoint de ${\goth l}$. Si $x_{0}$ est un \'el\'ement de ${\goth l}'$ tel que 
$(x_{0},y)$ soit un \'el\'ement adh\'erent \`a l'orbite de $(x,y)$ sous l'action de $A$, 
qui n'appartient pas \`a cette orbite, alors pour tout \'el\'ement nilpotent $v$ du 
radical de ${\goth p}$, $(x_{0},y)$ est adh\'erent \`a l'orbite de $(x+v,y)$ sous 
l'action de $G'(y)$ car $(x,y)$ est adh\'erent \`a l'orbite de $(x+v,y)$ sous l'action de
$G'(y)$; donc dans ce cas, $(x_{0},\pi (y))$ n'appartient pas \`a 
$\Sigma '_{{\goth l}'}$ car $Y$ ne rencontre pas $\Sigma _{{\goth g}}$. Vu l'arbitraire 
du point $(x,\pi (y))$ de $Y_{1}$, $Y_{1}$ ne rencontre pas $\Sigma _{{\goth l}'}$, 
d'o\`u la proposition.
\end{proof}
\section{Th\'eor\`eme d'exactitude et id\'eaux premiers.} \label{ei} 
On rappelle que ${\goth C}_{{\goth g}}$ d\'esigne la sous-vari\'et\'e des
\'el\'ements $(x,y)$ de ${\goth g}\times {\goth g}$ tels que $[x,y]$
soit nul, que $I_{{\goth g}}$ est l'id\'eal de $\tk {{\Bbb C}}{\e Sg}\e Sg$
engendr\'e par les fonctions $(x,y) \mapsto \dv v{[x,y]}$ o\`u $v$ est dans
${\goth g}$, que les complexes $E_{\bullet}({\goth g})$ et
$\overline{E}_{\bullet}({\goth g})$ sont d\'efinis en \ref{dca}.

\begin{Th}\label{tfi}
On suppose que les facteurs simples de ${\goth g}$ ont la propri\'et\'e $({\bf N})$
au sens de {\rm \ref{dsu3}}. Alors pour tout entier naturel $k$, la dimension projective 
du module $\ex k{{\goth g}}\wedge \ex {\k g{}}{B_{{\goth g}}}$ est inf\'erieure \`a 
$k$.
\end{Th}

\begin{proof}
Pour $k$ entier strictement sup\'erieur \`a $\dim {\goth g} - \k g{}$, \sloppy
\hbox{$\ex k{{\goth g}}\wedge \ex {\k g{}}{B_{{\goth g}}}$} est le module nul;
donc il s'agit de montrer que pour \sloppy \hbox{$k=0,\ldots,\dim {\goth g}-\k g{}$}, la
dimension projective de $\ex k{{\goth g}}\wedge \ex {\k g{}}{B_{{\goth g}}}$
est inf\'erieure \`a $k$. On consid\`ere successivement les deux cas suivants:
\begin{list}{}{}
\item 1) ${\goth g}$ est simple,
\item 2) ${\goth g}$ n'est pas simple.
\end{list}

1) On montre en raisonnant par r\'ecurrence sur $k$ que pour \sloppy \hbox{
$k=0,\ldots,\dim {\goth g}-\k g{}$}, la dimension projective de 
$\ex k{{\goth g}}\wedge \ex {\k g{}}{B_{{\goth g}}}$ est inf\'erieure \`a $k$.
Puisque $B_{{\goth g}}$ est un module libre, $\ex {\k g{}}{B_{{\goth g}}}$ est
un module libre. Soit $k$ un entier strictement positif, inf\'erieur \`a
$\dim {\goth g} - \k g{}$ et tel que pour $i=0,\ldots,k-1$, la dimension
projective de $\ex i{{\goth g}}\wedge \ex {\k g{}}{B_{{\goth g}}}$ est 
inf\'erieure \`a $i$. On rappelle que le complexe 
$D_{k}^{\bullet}({\goth g},B_{{\goth g}})$ est d\'efini en \ref{dco4}. D'apr\`es
l'assertion (ii) du lemme \ref{l2co4}, $X_{{\goth g}}$ contient le support
dans ${\goth g}\times {\goth g}$ de la cohomologie du complexe 
$D_{k}^{\bullet}({\goth g},B_{{\goth g}})$. Soit $Y$ une composante irr\'eductible 
du support de la cohomologie du complexe $D_{k}^{\bullet}({\goth g},B_{{\goth g}})$.
Alors $Y$ est invariant pour les actions de $G$ et GL$_{2}({\Bbb C})$ dans 
${\goth g}\times {\goth g}$ d'apr\`es les assertions (iv) et (v) de la proposition 
\ref{psc}. On suppose que $Y$ rencontre la partie $\Sigma '_{{\goth g}}$ de 
${\goth g}\times {\goth g}$ qui est d\'efinie en \ref{dq}. Il s'agit d'aboutir \`a une 
contradiction. Il existe une sous-alg\`ebre parabolique admissible ${\goth p}$ de 
${\goth g}$ au sens de \ref{dq1}, un \'el\'ement semi-simple r\'egulier $y$ de 
${\goth g}$ qui appartient \`a ${\goth p}$ et un \'el\'ement $x$ de ${\goth p}$ tel que 
$(x,y)$ appartienne \`a l'intersection de $Y$ et de $\Sigma _{{\goth g},{\goth p}}$.
Soit $V$ un ouvert affine de ${\goth p}\times {\goth p}$ qui contient $(x,y)$ et 
$L$ un sous-module de $\tk {{\Bbb C}}{{\Bbb C}[V]}{\goth p}$ qui satisfont les conditions
du corollaire \ref{csc4}. D\'esignant par $B_{V}$ le sous-module de 
$\tk {{\Bbb C}}{{\Bbb C}[V]}{\goth p}$ engendr\'e par les restrictions \`a $V$ des 
\'el\'ements de $B_{{\goth g}}$, il r\'esulte du corollaire \ref{cco4} que le
complexe $D_{k}^{\bullet}({\goth g},B_{V})$ n'a pas de cohomologie en degr\'e diff\'erent
de $\k g{}$. Soit $W$ un ouvert affine de ${\goth g}\times {\goth g}$ qui contient 
$(x,y)$ et dont l'intersection avec ${\goth p}\times {\goth p}$ est contenue dans $V$.
On note $B_{W}$ le sous-module de $\tk {{\Bbb C}}{{\Bbb C}[W]}{\goth g}$ engendr\'e
par les restrictions \`a $W$ des \'el\'ements de $B_{{\goth g}}$. Puisque les modules
$B_{V}$ et $B_{W}$ sont libres et de m\^eme rang, il r\'esulte du lemme \ref{lco5}
que le support dans $W$ des groupes de cohomologie en degr\'e diff\'erent de $\k g{}$
du complexe $D_{k}^{\bullet}({\goth g},B_{W})$ ne rencontre pas $V$; or d'apr\`es
l'assertion (iii) du lemme \ref{l2co4}, le complexe $D_{k}^{\bullet}({\goth g},W)$ 
n'a pas de cohomologie en degr\'e $\k g{}$ car $\Lambda _{{\goth g}}$ est un grand ouvert
de ${\goth g}\times {\goth g}$ d'apr\`es le corollaire \ref{csc3}; donc il existe
un ouvert de ${\goth g}\times {\goth g}$ qui contient $(x,y)$ et qui ne rencontre pas
le support de la cohomologie du complexe $D_{k}^{\bullet}({\goth g},B_{{\goth g}})$, 
d'o\`u la contradiction. Puisque $Y$ est ferm\'e et invariant pour les actions de $G$ et
de GL$_{2}({\Bbb C})$ dans ${\goth g}\times {\goth g}$, il en r\'esulte que $Y$ ne 
rencontre pas $\Sigma _{{\goth g}}$. Vu l'arbitraire de $Y$, le support de la cohomologie
du complexe $D_{k}^{\bullet}({\goth g},B_{{\goth g}})$ ne rencontre pas 
$\Sigma _{{\goth g}}$.

Si ${\goth g}$ est de rang $1$, alors $X_{{\goth g}}$ est contenu dans 
${\cal X}_{{\goth g}}$. Il r\'esulte alors de la proposition \ref{pq}, du th\'eor\`eme
\ref{ta} et de la propri\'et\'e $({\bf N})$ pour ${\goth g}$ que la codimension dans 
${\goth g}\times {\goth g}$ du support de la cohomologie du complexe 
$D_{k}^{\bullet}({\goth g},B_{{\goth g}})$ est sup\'erieure \`a $k+1$ car $k$ est 
inf\'erieur \`a $\k g{}-\j g{}$; donc d'apr\`es le corollaire \ref{c2p} et l'hypoth\`ese 
de r\'ecurrence, le complexe $D_{k}^{\bullet}({\goth g},B_{{\goth g}})$ n'a pas de 
cohomologie en degr\'e strictement inf\'erieur \`a $k+\k g{}$. Par d\'efinition, le 
complexe $D_{k}^{\bullet}({\goth g},B_{{\goth g}})$ n'a pas de cohomologie en degr\'e 
$k+\k g{}$ car $k$ est inf\'erieur \`a $\dim {\goth g}-\k g{}$; donc le complexe 
$D_{k}^{\bullet}({\goth g},B_{{\goth g}})$ est acyclique. Il r\'esulte 
alors de l'assertion (ii) du lemme \ref{lp} et de l'hypoth\`ese de r\'ecurrence que la 
dimension projective du module $\ex k{{\goth g}}\wedge \ex {\k g{}}{B_{{\goth g}}}$ est 
inf\'erieure \`a $k$.

2) On note ${\goth z}$ le centre de ${\goth g}$ et 
$\poi {{\goth g}}1{,\ldots,}{m}{}{}{}$ les facteurs simples de ${\goth g}$. 
On montre alors en raisonnant par r\'ecurrence sur $m$ que pour
$k=0,\ldots,\dim {\goth g}-\k g{}$, la dimension projective de 
$\ex k{{\goth g}}\wedge \ex {\k g{}}{B_{{\goth g}}}$ est inf\'erieure \`a $k$.
Pour $m$ nul, $\dim {\goth g}$ est \'egal \`a $\k g{}$. En outre,
$\ex {\k g{}}{B_{{\goth g}}}$ est un module libre car $B_{{\goth g}}$ est un 
module libre. On suppose l'assertion vraie pour $m-1$. Soit ${\goth a}$ la 
somme de ${\goth z}$ et des id\'eaux $\poi {{\goth g}}1{,\ldots,}{m-1}{}{}{}$.
Alors ${\goth g}$ est somme directe de ${\goth a}$ et de ${\goth g}_{m}$; donc
$B_{{\goth g}}$ est le sous-module de 
$\tk {{\Bbb C}}{\e Sg}\tk {{\Bbb C}}{\e Sg}{\goth g}$ engendr\'e par 
$B_{{\goth a}}$ et $B_{{\goth g}_{m}}$, et on a
$$\dim {\goth g} = \dim {\goth a} + \dim {\goth g}_{m} \mbox{ , }
\k g{} = \k a{} + \k gm \mbox{ .}$$
Il en r\'esulte que pour $k=\k g{},\ldots,\dim {\goth g}$, 
$\ex {k-\k g{}}{{\goth g}}\wedge \ex {\k g{}}{B_{{\goth g}}}$ est isomorphe \`a 
$$\bigoplus _{i=\sup \{\k a{},k-\dim {\goth g}_{m}\}}
^{\inf\{\dim {\goth a},k-\k gm \}} 
\tk {{\Bbb C}}{\ex {i-\k a{}}{{\goth a}}\wedge 
\ex {\k a{}}{B_{{\goth a}}}}\ex {k-i-\k gm}{{\goth g}_{m}}\wedge 
\ex {\k gm}{B_{{\goth g}_{m}}} \mbox{ .}$$
Puisque ${\goth g}_{m}$ est simple, la dimension projective de 
$\ex k{{\goth g}_{m}}\wedge \ex {\k gm}{B_{{\goth g}_{m}}}$ est inf\'erieure \`a
$k$ pour tout entier naturel $k$, d'apr\`es (1); donc d'apr\`es
l'hypoth\`ese de r\'ecurrence, la dimension projective du 
$\tk {{\Bbb C}}{\e Sg}\e Sg$-module
$$\tk {{\Bbb C}}{\ex {i-\k a{}}{{\goth a}}\wedge 
\ex {\k a{}}{B_{{\goth a}}}}\ex {k-i-\k gm}{{\goth g}_{m}}\wedge 
\ex {\k gm}{B_{{\goth g}_{m}}} $$
est inf\'erieure \`a l'entier
$$ i- \k a{} + k-i-\k gm = k-\k g{} \mbox{ ,}$$
pour tout entier $k$ sup\'erieur \`a $\k g{}$ et pour tout entier $i$ 
sup\'erieur \`a	$\k a{}$ et inf\'erieur \`a $k-\k gm$; donc pour
$k=0,\ldots,\dim {\goth g}-\k g{}$, la dimension projective de \sloppy \hbox{
$\ex k{{\goth g}}\wedge \ex {\k g{}}{B_{{\goth g}}}$} est inf\'erieure \`a $k$,
d'o\`u le th\'eor\`eme.
\end{proof}

\begin{Th}\label{tei}
On suppose que pour tout \'el\'ement semi-simple $\rho $ de ${\goth g}$, les facteurs 
simples de ${\goth g}(\rho )$ ont la propri\'et\'e $({\bf N})$ au sens de 
{\rm \ref{dsu3}}.

{\rm i)} Le complexe r\'eduit $\overline{E}_{\bullet}({\goth g})$ de
${\goth g}$ est acyclique.

{\rm ii)} L'id\'eal $I_{{\goth g}}$ est premier.

{\rm iii)} Le complexe canonique de deuxi\`eme esp\`ece
$E_{\bullet}({\goth g})$ de ${\goth g}$ n'a pas d'homologie en degr\'e
diff\'erent de $\k g{}$ et son homologie est isomorphe \`a
l'alg\`ebre des fonctions r\'eguli\`eres sur ${\goth C}_{{\goth g}}$.
\end{Th}

\begin{proof}
On d\'esigne par $J_{{\goth g}}$ le radical de $I_{{\goth g}}$. Par
d\'efinition, les complexes $E_{\bullet}({\goth g})$ et 
$\overline{E}_{\bullet}({\goth g})$ ont m\^eme groupe d'homologie en
degr\'e diff\'erent de $\k g{}$; donc d'apr\`es le lemme \ref{lca2}, il
suffit de montrer l'assertion (i) du th\'eor\`eme. Pour $j$ entier strictement
sup\'erieur \`a $\k g{}$, on note $Z_{j}$ l'espace des cycles de degr\'e $j$
du complexe $E_{\bullet}({\goth g})$. On a alors un complexe
$$  0 \rightarrow E_{\dim {\goth g}}({\goth g}) \rightarrow \cdots 
\rightarrow E_{j+1}({\goth g}) \rightarrow Z_{j} \rightarrow 0 
\mbox{ .}$$
D'apr\`es le th\'eor\`eme \ref{tfi}, la dimension projective de
$E_{i}({\goth g})$ est inf\'erieure \`a \sloppy \hbox{$i-\k g{}=i-j+j-\k g{}$}. En outre,
$E_{i}({\goth g})$ est nul si $i-j$ est strictement sup\'erieur \`a 
$\dim {\goth g}-j$. La dimension de 
${\goth C}_{{\goth g}}$ \'etant \'egale \`a $\dim {\goth g} + \j g{}$, sa
codimension dans ${\goth g}\times {\goth g}$ est \'egale 
$\dim {\goth g} - \j g{}$. Des \'egalit\'es
$$ \dim {\goth g} - \j g{} = 2 (\dim {\goth g}-\k g{}) \mbox{ , }
2\dim {\goth g} - \dim {\goth C}_{{\goth g}} = \dim {\goth g} - \j g{}
\mbox{ ,}$$ 
on d\'eduit l'in\'egalit\'e
$$ 2(\dim {\goth g} - j) + j - \k g{} < 2\dim {\goth g} - 
\dim {\goth C}_{{\goth g}}\mbox{ ,}$$	
car $j$ est strictement sup\'erieur \`a $\k g{}$; or d'apr\`es l'assertion
(iii) du lemme \ref{lco3}, ${\goth C}_{{\goth g}}$ contient le support de
l'homologie du  complexe ci-dessus; donc d'apr\`es le corollaire \ref{cp},
le complexe ci-dessus est exact. Il en r\'esulte que le complexe
$E_{\bullet}({\goth g})$ n'a pas d'homologie en degr\'e strictement
sup\'erieur \`a $\k g{}$. Puisque $\overline{E}_{k}({\goth g})$ est nul pour
$k$ strictement inf\'erieur \`a $\k g{}$, le complexe
$\overline{E}_{\bullet}({\goth g})$ n'a pas d'homologie en degr\'e
diff\'erent de $\k g{}$. D'apr\`es le lemme \ref{lca2}, il reste \`a montrer
l'\'egalit\'e des id\'eaux $I_{{\goth g}}$ et $J_{{\goth g}}$. 

On raisonne par r\'ecurrence sur la dimension de ${\goth g}$. Si 
${\goth g}$ est une alg\`ebre de Lie commutative, $I_{{\goth g}}$ est
l'id\'eal nul. On suppose que ${\goth g}$ n'est pas une alg\`ebre de Lie
commutative et que le th\'eor\`eme est vrai pour toute alg\`ebre de Lie
r\'eductive de dimension strictement inf\'erieure \`a celle de ${\goth g}$.
On rappelle que $Y_{0}({\goth g})$ d\'esigne le support dans 
${\goth g}\times {\goth g}$ de $J_{{\goth g}}/I_{{\goth g}}$ et que 
$X_{0}({\goth g})$ est l'image de $Y_{0}({\goth g})$ par la premi\`ere 
projection de ${\goth g}\times {\goth g}$ sur ${\goth g}$. D'apr\`es 
l'hypoth\`ese r\'ecurrence et la proposition \ref{ps}, $X_{0}({\goth g})$ ne 
contient pas d'\'el\'ement semi-simple non central car pour tout \'el\'ement
semi-simple $\rho $ de ${\goth g}$ et tout \'el\'ement semi-simple $\rho '$ de 
${\goth g}(\rho )$, le centralisateur de $\rho '$ dans ${\goth g}(\rho )$ est le
centralisateur dans ${\goth g}$ d'un \'el\'ement semi-simple de ${\goth g}$; or 
$X_{0}({\goth g})$ est ferm\'e et $G$-invariant d'apr\`es le corollaire \ref{cca1}; donc 
$X_{0}({\goth g})$ est contenu dans l'ensemble 
${\goth N}'_{{\goth g}}$ des \'el\'ements de ${\goth g}$ dont l'image par la 
repr\'esentation adjointe est un endomorphisme nilpotent de ${\goth g}$. 
D'apr\`es l'assertion (iii) du lemme \ref{lca1}, $Y_{0}({\goth g})$ est stable 
par l'involution $(x,y) \mapsto (y,x)$; donc $Y_{0}({\goth g})$ est contenu
dans l'intersection de ${\goth N}'_{{\goth g}}\times {\goth N}'_{{\goth g}}$
et de ${\goth C}_{{\goth g}}$ car ${\goth C}_{{\goth g}}$ est la vari\'et\'e
des z\'eros de $I_{{\goth g}}$. D'apr\`es \cite{Pr}, la dimension de l'intersection de 
${\goth N}'_{{\goth g}}\times {\goth N}'_{{\goth g}}$ et de 
${\goth C}_{{\goth g}}$ est \'egale \`a 
$\dim {\goth g}+ \dim {\goth z}_{{\goth g}}$ en d\'esignant
par ${\goth z}_{{\goth g}}$ le centre de ${\goth g}$; donc la codimension
de $Y_{0}({\goth g})$ dans ${\goth g}\times {\goth g}$ est \'egale
\`a $\dim {\goth g} - \dim {\goth z}_{{\goth g}}$. En particulier, elle
strictement sup\'erieure \`a $2(\dim {\goth g} -\k g{})$ car
${\goth g}$ n'est pas une alg\`ebre de Lie commutative; donc d'apr\`es ce
qui pr\'ec\`ede, le support dans ${\goth g}\times {\goth g}$ de la
cohomologie du complexe $\overline{E}_{\bullet}({\goth g})$ est de
codimension strictement sup\'erieure \`a $2(\dim {\goth g} - \k g{})$. Le
complexe $\overline{E}_{\bullet}({\goth g})$ est alors acyclique d'apr\`es
le corollaire \ref{cp} et le th\'eor\`eme \ref{tfi}.
\end{proof}

\section{La propri\'et\'e $({\bf D})$ pour les alg\`ebres de Lie simples.}\label{a}
Dans cette section, on suppose ${\goth g}$ simple de dimension strictement sup\'erieure 
\`a $3$, et on utilise les notations de \ref{q} Le centralisateur ${\goth h}$ de $\rho $ 
dans ${\goth g}$ est une sous-alg\`ebre de Cartan de ${\goth g}$ qui est contenue dans 
${\goth b}$. On note $R$ le syst\`eme de racines de ${\goth h}$ dans ${\goth g}$, 
$R_{+}$ le syst\`eme de racines positives de $R$ d\'efini par ${\goth b}$ et $B$ la base 
de $R_{+}$. Pour tout $\alpha $ dans $R$, ${\goth g}^{\alpha }$ d\'esigne le sous-espace 
radiciel de ${\goth g}$ pour la racine $\alpha $. On d\'esigne respectivement par 
${\bf B}$ et ${\bf H}$ les sous-groupes ferm\'es connexes de $G$ dont les alg\`ebres de 
Lie sont les images de ${\goth b}$ et de ${\goth h}$ par la repr\'esentation adjointe de 
${\goth g}$. On note ${\goth u}$ l'ensemble des \'el\'ements nilpotents de ${\goth b}$, 
${\goth b}_{-}$ la somme des sous-espaces ${\goth h}$ et ${\goth g}^{-\alpha }$ o\`u 
$\alpha $ est dans $R_{+}$, ${\goth u}_{-}$ l'ensemble des \'el\'ements nilpotents de 
${\goth b}_{-}$, ${\bf H}'$ le sous-groupe des automorphismes lin\'eaires de 
${\goth g}\times {\goth g}$ qui est engendr\'e par les automorphismes 
$(v,w)\mapsto (tv,w)$ et $(v,w)\mapsto (g(v),g(w))$ o\`u $t$ et $g$ sont
respectivement dans ${\Bbb C}\backslash \{0\}$ et dans ${\bf H}$.

\subsection{} Soient $\tau $ la projection canonique de ${\goth g}\times {\goth h}$ sur 
${\goth h}$, ${\goth g}'_{\r}$ l'ensemble des \'el\'ements semi-simples r\'eguliers de 
${\goth g}$, ${\goth h}_{\r}$ l'intersection de ${\goth h}$ avec ${\goth g}'_{\r}$ et 
$N_{G}({\goth h})$ le normalisateur de ${\goth h}$ dans $G$. On rappelle que la partie 
$\Sigma _{{\goth g}}$ est d\'efinie en \ref{dq}.

\begin{Def}\label{da1}
On dira qu'une partie $X$ de ${\goth g}\times {\goth h}_{\r}$ a la propri\'et\'e 
$({\bf P})$ si elle satisfait les conditions suivantes:
\begin{list}{}{}
\item {\rm 1)} $X$ est constructible et ne rencontre pas $\Sigma _{{\goth g}}$,
\item {\rm 2)} $X$ est invariant pour l'action de $N_{G}({\goth h})$ dans 
${\goth g}\times {\goth h}_{\r}$, 
\item {\rm 3)} le groupe $N_{G}({\goth h})$ agit transitivement sur les composantes 
irr\'eductibles de $X$,
\item {\rm 4)} pour tout $(s,t)$ dans ${\Bbb C}\times ({\Bbb C}\backslash \{0\})$, $X$ 
contient $(sx,ty)$ s'il contient $(x,y)$.
\end{list}
On d\'esigne par ${\cal P}$ l'ensemble des parties de ${\goth g}\times {\goth h}_{\r}$ 
qui ont la propri\'et\'e $({\bf P})$.
\end{Def}

Pour $\alpha $ dans $B$, on note $R_{\alpha }$ l'ensemble des \'el\'ements
de $R$ dont la coordonn\'ee en $\alpha $, dans la base $B$, est non nulle et
${\goth p}_{\alpha }$ la sous-alg\`ebre parabolique de ${\goth g}$ qui contient 
${\goth b}$ et ${\goth g}^{-\beta }$ pour tout $\beta $ dans $B\backslash \{\alpha \}$. 
Alors la somme ${\goth p}_{\alpha ,\u}$ des sous-espaces ${\goth g}^{\gamma }$, o\`u
$\gamma $ est dans $R_{\alpha }$, est l'ensemble des \'el\'ements nilpotents du radical de
${\goth p}_{\alpha }$ et la sous-alg\`ebre ${\goth l}_{\alpha }$, engendr\'ee par 
${\goth h}$ et les sous-espaces ${\goth g}^{\beta }$ et ${\goth g}^{-\beta }$, o\`u 
$\beta $ est dans $B\backslash \{\alpha \}$, est un facteur r\'eductif de 
${\goth p}_{\alpha }$. En outre, la somme ${\goth p}_{\alpha ,\u,-}$ des sous-espaces
${\goth g}^{-\gamma }$, o\`u $\gamma $ est dans $R_{\alpha }$, est un suppl\'ementaire
de ${\goth p}_{\alpha }$ dans ${\goth g}$ et la somme ${\goth p}_{\alpha ,-}$ de
${\goth l}_{\alpha }$ et de ${\goth p}_{\alpha ,\u,-}$ est une sous-alg\`ebre parabolique
de ${\goth g}$ qui contient la sous-alg\`ebre de Borel ${\goth b}_{-}$.

\begin{Def}\label{d2a1}
On dira qu'une sous-alg\`ebre parabolique ${\goth p}$ de ${\goth g}$ est extr\'emale si
elle est maximale et admissible au sens de \ref{dq1}.
\end{Def}

\begin{Rem}\label{ra1}
Pour $\alpha $ dans $B$, la sous-alg\`ebre parabolique ${\goth p}_{\alpha }$ est 
extr\'emale si et seulement si $\alpha $ est une extr\'emit\'e du diagramme de Dynkin
de support $B$. En outre, toute sous-alg\`ebre parabolique extr\'emale est conjugu\'ee par
l'action de $G$ \`a un ${\goth p}_{\alpha }$, o\`u $\alpha $ est une extr\'emit\'e du
diagramme de Dynkin de support $B$.  
\end{Rem}

Soit $\alpha $ un \'el\'ement de $B$ tel que ${\goth p}_{\alpha }$ soit une 
sous-alg\`ebre parabolique extr\'emale. Pour toute partie $I$ de $R_{\alpha }$, on note 
$F_{\alpha ,I}$ la somme de ${\goth p}_{\alpha }$ et des sous-espaces 
${\goth g}^{-\beta }$ o\`u $\beta $ est dans $I$. On rappelle que la
hauteur d'un \'el\'ement de $R_{+}$ est la somme de ses coordonn\'ees dans la base $B$.

\begin{lemme}\label{la1}
Soient $X$ dans ${\cal P}$ et $Y$ une composante irr\'eductible de $X$. 

{\rm i)} Si $X$ est contenu dans ${\cal X}_{{\goth g}}$, alors $X$ a une composante 
irr\'eductible qui est contenue dans ${\goth b}\times {\goth h}_{\r}$.

{\rm ii)} Il existe une unique partie minimale $I$ de $R_{\alpha }$ telle que
$F_{\alpha ,I}\times {\goth h}_{\r}$ contienne $Y$.

{\rm iii)} Soient $x$, $w$, $v$ des \'el\'ements respectivement dans ${\goth h}_{\r}$,
${\goth p}_{\alpha }$, ${\goth p}_{\alpha ,\u,-}$. Alors $(v,x)$ appartient \`a 
l'adh\'erence dans ${\goth g}\times {\goth h}_{\r}$ de l'orbite de $(v+w,x)$ sous 
l'action de ${\bf H}'$.

{\rm iv)} Soit $\gamma $ un \'el\'ement de plus grande hauteur de $I$. Alors pour tout
$x$ dans un ouvert non vide de $\tau (Y)$, il existe un \'el\'ement $v$ de ${\goth g}$
et un \'el\'ement non nul $w$ de ${\goth g}^{-\gamma }$ qui satisfont les conditions 
suivantes: $Y$ contient $(v,x)$ et $(w,x)$ appartient \`a l'adh\'erence de l'orbite de 
$(v,x)$ sous l'action de ${\bf H}'$.
\end{lemme}

\begin{proof}
i) On suppose $X$ contenu dans ${\cal X}_{{\goth g}}$. Pour tout $g$ dans 
$N_{G}({\goth h})$, on d\'esigne par $Y_{g}$ l'ensemble des \'el\'ements $(x,y)$ de $Y$
tels que ${\goth b}$ contienne $g(x)$ et $g(y)$. Puisque $N_{G}({\goth h})$ op\`ere 
transitivement sur l'ensemble des sous-alg\`ebres de Borel de ${\goth g}$ qui contiennent
un \'el\'ement de ${\goth h}_{\r}$, $Y$ est r\'eunion des parties $Y_{g}$ o\`u $g$ est 
dans $N_{G}({\goth h})$ car $X$ est contenu dans ${\cal X}_{{\goth g}}$. Pour tout $g$ 
dans $N_{G}({\goth h})$ et pour tout $h$ dans ${\bf H}$, $Y_{g}$ est \'egal \`a $Y_{hg}$ 
car ${\goth b}$ contient ${\goth h}$; or $Y$ est irr\'eductible et $Y_{g}$ est ferm\'e 
dans $Y$ pour tout $g$ dans $N_{G}({\goth h})$; donc il existe un \'el\'ement $g$ de 
$N_{G}({\goth h})$ tel que $Y$ soit \'egal \`a $Y_{g}$ car le quotient de 
$N_{G}({\goth h})$ par ${\bf H}$ est fini, d'o\`u l'assertion.  

ii) Puisque ${\goth g}$ est somme directe de ${\goth p}_{\alpha }$ et de 
${\goth p}_{\alpha ,\u,-}$, il existe une partie $I$ de $R_{\alpha }$ telle que
$F_{\alpha ,I}\times {\goth h}_{\r}$ contienne $Y$. Si $I_{1}$ et $I_{2}$ sont deux 
parties de $R_{\alpha }$ telles que $F_{\alpha ,I_{1}}\times {\goth h}_{\r}$ et 
$F_{\alpha ,I_{2}}\times {\goth h}_{\r}$ contiennent $Y$, 
$F_{\alpha ,I}\times {\goth h}_{\r}$ contient $Y$ si $I$ est l'intersection de $I_{1}$ et
de $I_{2}$, d'o\`u l'existence d'une unique partie minimale $I$ de $R_{\alpha }$ qui
satisfait les conditions de l'assertion.

iii) On note $t\mapsto h(t)$ le sous-groupe \`a un param\`etre de $G$ dont l'alg\`ebre de 
Lie est l'image du centre de ${\goth l}_{\alpha }$ par la repr\'esentation adjointe de 
${\goth g}$. On d\'esigne par $w_{1}$ et $w_{2}$ les composantes respectives de $w$ sur 
${\goth l}_{\alpha }$ et ${\goth p}_{\alpha ,\u}$. Puisque $Y$ est invariant par 
${\bf H}$, $Y$ contient $(h(t)(v+w_{1}+w_{2}),h(t)(x))$ pour tout \'el\'ement non nul $t$
de ${\Bbb C}$; donc d'apr\`es la condition (3), $Y$ contient $(v+tw_{1}+t^{2}w_{2},x)$ 
pour tout \'el\'ement non nul $t$ de ${\Bbb C}$, d'o\`u l'assertion.

iv) Soit $E^{*}_{\alpha ,I}$ l'ensemble des \'el\'ements de la somme des 
sous-espaces ${\goth g}^{-\beta }$, o\`u $\beta $ est dans $I$, qui ont une composante 
non nulle sur chacun de ces sous-espaces. D'apr\`es la minimalit\'e de la partie $I$, 
$Y$ \'etant irr\'eductible, pour tout $x$ dans un ouvert non vide $T'$ de $\tau (Y)$, il 
existe un \'el\'ement $v$ de $E^{*}_{\alpha ,I}$ et un \'el\'ement $w$ de 
${\goth p}_{\alpha }$ tel que $Y$ contienne $(v+w,x)$. Soient $x$ dans $T'$ et $v$ dans 
$E^{*}_{\alpha ,I}$ tels qu'il existe un \'el\'ement $w$ de ${\goth p}_{\alpha }$ pour 
lequel $(v+w,x)$ appartient \`a $Y$. D'apr\`es (iii), $(v,x)$ appartient \`a 
l'adh\'erence de l'orbite de $(v+w,x)$ sous l'action de ${\bf H}'$. On note 
$t\mapsto g(t)$ le sous-groupe \`a un param\`etre de $G$ dont l'alg\`ebre 
de Lie est engendr\'ee par $\ad \rho $. Soient $J$ l'ensemble des \'el\'ements de plus
grande hauteur de $I$ et $v_{J}$ la composante de $v$ sur la somme des sous-espaces
${\goth g}^{-\beta }$, o\`u $\beta $ est dans $J$. Alors $Y$ contient $(g(t)(v),x)$ pour 
tout \'el\'ement non  nul $t$ de ${\Bbb C}$; or $t^{\dv {\gamma }{\rho }}g(t)(v)$ tend 
vers $v_{J}$ quand $t$ tend vers $0$; donc $(v_{J},x)$ appartient \`a l'adh\'erence de 
l'orbite de $(v+w,x)$ sous l'action de ${\bf H}'$. Les \'el\'ements de $J$ \'etant 
lin\'eairement ind\'ependants, comme \'el\'ements de m\^eme hauteur, il 
existe un \'el\'ement $k$ de ${\goth h}$ qui a les propri\'et\'es suivantes:
\begin{list}{}{}
\item a) pour $\beta $ dans $B$, $\dv {\beta }k$ est entier, 
\item b) $\dv {\gamma }k$ est un entier strictement positif,
\item c) $k$ appartient \`a l'intersection des noyaux des \'el\'ements de $J$, distincts
de $\gamma $.
\end{list}
Soient $v_{-\gamma }$ la composante de $v_{J}$ sur ${\goth g}^{-\gamma }$ et 
$t\mapsto k(t)$ le sous-groupe \`a un param\`etre de $G$ dont l'alg\`ebre de Lie 
est engendr\'ee par $\ad k$. Puisque $t^{\dv {\gamma }{\rho }}k(t)(v_{J})$ tend vers 
$v_{-\gamma }$ quand $t$ tend vers $0$, $(v_{-\gamma },x)$ appartient \`a 
l'adh\'erence de l'orbite de $(v+w,x)$ sous l'action de ${\bf H}'$, d'o\`u l'assertion. 
\end{proof}
\subsection{} On note ${\cal P}_{0}$ l'ensemble des \'el\'ements de ${\cal P}$ qui sont 
contenus dans $X_{{\goth g}}$. 

\begin{lemme}\label{la2}
Soient $h$ un \'el\'ement de ${\goth h}$ qui est \'el\'ement r\'egulier de ${\goth g}$
et $x$ un \'el\'ement nilpotent de ${\goth g}$ qui appartient \`a ${\goth b}$. On suppose
que $(x,h)$ est un \'el\'ement de $X_{{\goth g}}$ qui n'appartient pas \`a la partie 
$\Sigma '_{{\goth g}}$ de ${\goth g}\times {\goth g}$ d\'efinie en {\rm \ref{dq}}. Alors 
pour tout $\alpha $ dans $B$, la composante de $x$ sur ${\goth g}^{\alpha }$, 
relativement \`a la d\'ecomposition de ${\goth b}$,
$${\goth b} = {\goth h} \oplus \bigoplus _{\beta \in R_{+}} {\goth g}^{\beta } \mbox{ ,}$$
est nulle.
\end{lemme}
 
\begin{proof}
On suppose qu'il existe $\alpha $ dans $B$ tel que la composante de $x$ sur 
${\goth g}^{\alpha }$ ne soit pas nulle. Il s'agit d'aboutir \`a une contradiction. 
Soit $I$ une partie connexe, non vide, du diagramme de Dynkin de support $B$ telle
que la composante de $x$ sur ${\goth g}^{\alpha }$ soit non nulle pour tout $\alpha $
dans $I$. On note ${\goth p}_{I}$ la sous-alg\`ebre de ${\goth g}$ engendr\'ee par 
${\goth b}$ et les sous-espaces ${\goth g}^{-\gamma }$ o\`u $\gamma $ est dans $I$. Alors
${\goth p}_{I}$ est une sous-alg\`ebre parabolique admissible de ${\goth g}$ au sens de 
\ref{dq1} car $I$ est une partie connexe du diagramme de Dynkin de support $B$. En outre,
$h$ est un \'el\'ement r\'egulier du facteur r\'eductif ${\goth l}_{I}$ de 
${\goth p}_{I}$ qui est engendr\'e par ${\goth h}$ et les sous-espaces 
${\goth g}^{\gamma }$, ${\goth g}^{-\gamma }$ o\`u $\gamma $ est dans $I$. D\'esignant 
par $x'$ l'image de $x$ par la projection canonique de ${\goth p}_{I}$ sur 
${\goth l}_{I}$, $x'$ est un \'el\'ement nilpotent r\'egulier de ${\goth l}_{I}$ qui
appartient \`a une sous-alg\`ebre de Borel de ${\goth l}_{I}$ qui contient l'\'el\'ement
semi-simple r\'egulier $h$; donc d'apr\`es l'assertion (iv) du lemme \ref{lsc5},
$(x',h)$ appartient \`a l'ouvert $\Lambda _{{\goth l}_{I}}$ de 
${\goth l}_{I}\times {\goth l}_{I}$. Il en r\'esulte que 
$\Sigma _{{\goth g},{\goth p}_{I}}$ contient $(x,h)$ car $X_{{\goth g}}$ contient
$(x,h)$. Ceci est absurde car $\Sigma _{{\goth g}}$ ne contient pas $(x,h)$.
\end{proof}

\begin{cor}\label{ca2}
Soient $\alpha $ dans $R$, $h$ un \'el\'ement r\'egulier de ${\goth g}$ qui appartient
\`a ${\goth h}$ et $w$ un \'el\'ement non nul de ${\goth g}^{\alpha }$. Alors  
$\Sigma '_{{\goth g}}$ contient $(w,h)$.
\end{cor}

\begin{proof}
Puisque ${\goth g}$ est irr\'eductible, deux racines de m\^eme longueur dans $R$ 
appartiennent \`a une m\^eme orbite du groupe de Weyl de ${\goth h}$ d'apr\`es
\cite{Bo}(Ch. VI, \S 1, Proposition 11); donc il existe un \'el\'ement $g$ dans 
$N_{G}({\goth h})$ et un \'el\'ement $\beta $ de $B$ tels que $g({\goth g}^{\beta })$
soit \'egal \`a ${\goth g}^{\alpha }$. Puisque $\j g{}$ est sup\'erieur \`a $2$,
$X_{{\goth g}}$ contient $(g^{-1}(w),g^{-1}(h))$ d'apr\`es l'assertion (v) du lemme 
\ref{lsc6}; donc d'apr\`es le lemme \ref{la2}, $(g^{-1}(w),g^{-1}(h))$ appartient \`a 
$\Sigma '_{{\goth g}}$. Il en r\'esulte que $\Sigma '_{{\goth g}}$ contient $(w,h)$ car 
$\Sigma '_{{\goth g}}$ est invariant pour l'action de $G$ dans 
${\goth g}\times {\goth g}$.
\end{proof}

On note $B_{\E}$ l'ensemble des exr\'emit\'es du diagramme de Dynkin de support $B$.

\begin{lemme}\label{l2a2}
On suppose que toute alg\`ebre de Lie simple, de rang strictement inf\'erieur \`a 
$\j g{}$, a la propri\'et\'e $({\bf D}')$ au sens de {\rm \ref{dq3}}. 

{\rm i)} Soit ${\goth p}$ une sous-alg\`ebre parabolique extr\'emale de ${\goth g}$. 
Si $X$ est un \'el\'ement de ${\cal P}_{0}$ dont une composante irr\'eductible
est contenue dans ${\goth p}\times {\goth h}_{\r}$, alors ${\cal X}_{{\goth g}}$
contient $X$.

{\rm ii)} Les \'el\'ements de ${\cal P}_{0}$ sont contenus dans ${\cal X}_{{\goth g}}$.
\end{lemme}

\begin{proof}
i) Soi $X$ un \'el\'ement de ${\cal P}_{0}$ dont une composante irr\'eductible est 
contenue dans ${\goth p}\times {\goth h}_{\r}$. On note $\tilde{X}$ l'image de 
$G\times X$ par l'application \sloppy \hbox{$(g,x,y)\mapsto (g(x),g(y))$}. Alors 
$\tilde{X}$ est une partie de $X_{{\goth g}}$ qui a la propri\'et\'e $({\bf Q})$ au sens 
de \ref{dq3}. En outre, ${\cal X}_{{\goth g},{\goth p}}$ contient $\tilde{X}$; donc 
d'apr\`es l'hypoth\`ese du lemme et la proposition \ref{pq3}, ${\cal X}_{{\goth g}}$ 
contient $\tilde{X}$ et $X$ car les facteurs simples des centralisateurs des \'el\'ements 
semi-simples non nuls de ${\goth g}$ sont de rang strictement inf\'erieur \`a $\j g{}$.

ii) Pour $X$ dans ${\cal P}$, on note $\nu (X)$ le plus petit entier tel qu'il existe
un \'el\'ement $\alpha $ de $B_{\E}$, une partie $I$ de $R_{\alpha }$ de cardinal 
$\nu (X)$ pour laquelle $F_{\alpha ,I}\times {\goth h}_{\r}$ contient une composante 
irr\'eductible de $X$ dont l'intersection avec ${\goth b}\times {\goth h}_{\r}$ est de 
dimension maximale. L'existence de l'entier $\nu (X)$ r\'esulte de l'assertion (ii) du 
lemme \ref{la1}. On suppose que ${\cal P}_{0}$ contient des \'el\'ements qui ne sont pas
contenus dans ${\cal X}_{{\goth g}}$. Il s'agit d'aboutir \`a une contradiction. 
D'apr\`es la condition (2) de la d\'efinition \ref{da1}, $N_{G}({\goth h})$ agit 
transitivement sur les composantes irr\'eductibles de $X$; donc $\nu (X)$ est nul si et 
seulement s'il existe une sous-alg\`ebre parabolique extr\'emale ${\goth p}$ de 
${\goth g}$ telle que ${\cal X}_{{\goth g},{\goth p}}$ contienne $X$. 

Soit $\nu $ le plus petit entier tel qu'il existe un \'el\'ement $X$ de ${\cal P}_{0}$ 
qui n'est pas contenu dans ${\cal X}_{{\goth g}}$ et tel que $\nu (X)$ soit \'egal \`a 
$\nu $. D'apr\`es l'assertion (i), $\nu $ est strictement positif. Par d\'efinition, il 
existe un \'el\'ement $\alpha $ de $B_{\E}$, une partie $I$ de $R_{\alpha }$, une 
composante irr\'eductible $Y$ de $X$ qui est contenue dans 
$F_{\alpha ,I}\times {\goth h}_{\r}$ et dont l'intersection avec 
${\goth b}\times {\goth h}_{\r}$ est de dimension maximale. Soit $\gamma $ un \'el\'ement
de $I$ de plus grande hauteur. D'apr\`es l'assertion (iv) du lemme \ref{la1}, il existe
un \'el\'ement $(v,x)$ de $Y$ et un \'el\'ement non nul $w$ de ${\goth g}^{-\gamma }$
tels que $(w,x)$ appartienne \`a l'adh\'erence de l'orbite de $(v,x)$ sous l'action de 
${\bf H}'$. Ceci est absurde d'apr\`es le corollaire \ref{ca2} car $Y$ ne rencontre pas 
$\Sigma _{{\goth g}}$, d'o\`u l'assertion.
\end{proof}

\begin{cor}\label{c2a2}
L'alg\`ebre de Lie ${\goth g}$ a la propri\'et\'e $({\bf D}')$ au sens de 
${\rm \ref{dq3}}$.
\end{cor}

\begin{proof}
On rappelle que $\Psi $ d\'esigne l'application \sloppy 
\hbox{$(g,x,y)\mapsto (g(x),g(y))$} de $G\times {\goth g}\times {\goth g}$ dans 
${\goth g}\times {\goth g}$. On raisonne par r\'ecurrence sur le rang de ${\goth g}$. Si 
${\goth g}$ est de rang $1$, alors ${\goth g}$ a la propri\'et\'e $({\bf D}')$ car 
$X_{{\goth g}}$ est contenu dans ${\cal X}_{{\goth g}}$. On suppose que toute alg\`ebre 
de Lie simple de rang strictement inf\'erieur \`a $\j g{}$ a la propri\'et\'e 
$({\bf D}')$. Soit $X$ une partie de $X_{{\goth g}}$ qui a la propri\'et\'e $({\bf Q})$ 
au sens de \ref{dq3} et qui ne rencontre pas $\Sigma _{{\goth g}}$. D'apr\`es 
l'assertion (ii) du lemme \ref{l2a2}, il suffit de montrer qu'il existe un \'el\'ement 
$Y$ de ${\cal P}_{0}$ tel que $X$ soit contenu dans l'adh\'erence dans 
${\goth g}\times {\goth g}$ de $\Psi (G\times Y)$.

On rappelle que ${\goth g}'_{\r}$ d\'esigne l'ensemble des \'el\'ements semi-simples
r\'eguliers de ${\goth g}$ et que ${\goth h}_{\r}$ d\'esigne l'intersection de ${\goth h}$
et de ${\goth g}'_{\r}$. Soient $\tilde{X}$ l'intersection de $X$ et de 
${\goth g}\times {\goth h}_{\r}$, $\poi Y1{,\ldots,}{m}{}{}{}$ les composantes 
irr\'eductibles de $\tilde{X}$. Pour $i=1,\ldots,m$, on note 
$\overline{\Psi (G\times Y_{i})}$ l'adh\'erence de $\Psi (G\times Y_{i})$ dans 
${\goth g}\times {\goth g}$. Puisque ${\goth g}'_{\r}$ est un ouvert de ${\goth g}$ qui 
rencontre l'image de $X$ par la deuxi\`eme projection de ${\goth g}\times {\goth g}$ sur 
${\goth g}$, l'intersection $X'$ de $X$ et de ${\goth g}\times {\goth g}'_{\r}$ est 
partout dense dans $X$ car $X$ est irr\'eductible; or $X'$ est \'egal \`a 
$\Psi (G\times \tilde{X})$ car l'orbite sous $G$ de tout \'el\'ement semi-simple 
r\'egulier de ${\goth g}$ rencontre ${\goth h}_{\r}$; donc $\Psi (G\times \tilde{X})$ est
partout dense dans $X$ et $X$ est contenu dans la r\'eunion des parties 
$\overline{\Psi (G\times Y_{1})},\ldots,\overline{\Psi (G\times Y_{m})}$. Par suite, il 
existe un entier $i$ tel que $\overline{\Psi (G\times Y_{i})}$ contienne $X$ car $X$ est 
irr\'eductible. On note $Y$ la r\'eunion des $g.Y_{i}$ o\`u $g$ est dans 
$N_{G}({\goth h})$. Puisque $X$ est invariant pour l'action de ${\bf H}$ dans
${\goth g}\times {\goth g}$, $\tilde{X}$ et $Y_{i}$ sont invariants pour l'action de 
${\bf H}$ dans ${\goth g}\times {\goth h}_{\r}$. Il en r\'esulte que $Y$ est une partie
constructible de ${\goth g}\times {\goth h}_{\r}$ comme r\'eunion finie de parties 
constructibles. Puisque $Y_{i}$ est une composante irr\'eductible de $\tilde{X}$, pour 
tout $(s,t)$ dans ${\Bbb C}\times ({\Bbb C}\backslash \{0\})$, $Y_{i}$ contient 
$(sx,ty)$ s'il contient $(x,y)$; donc pour tout $(s,t)$ dans 
${\Bbb C}\times ({\Bbb C}\backslash \{0\})$, $Y$ contient $(sx,ty)$ s'il contient 
$(x,y)$. Par suite, $Y$ est un \'el\'ement de ${\cal P}$ car $X$ ne rencontre pas 
$\Sigma _{{\goth g}}$. En outre, $X$ est contenu dans l'adh\'erence 
de $\Psi (G\times Y)$ et ${\cal P}_{0}$ contient $Y$ car $X$ est contenu dans 
$X_{{\goth g}}$, d'o\`u le corollaire.
\end{proof}

\begin{Th}\label{ta}
Soit ${\goth g}$ une alg\`ebre de Lie simple. Alors ${\goth g}$ a la propri\'et\'e 
$({\bf D})$ au sens de {\rm \ref{dq}}.
\end{Th}

\begin{proof}
Soit $X$ une partie ferm\'ee, irr\'eductible de $X_{{\goth g}}$, invariante pour les
actions de $G$ et de ${\rm GL}_{2}({\Bbb C})$, de codimension inf\'erieure \`a 
$\dim {\goth g}-4$, qui ne rencontre pas $\Sigma _{{\goth g}}$. Alors d'apr\`es la 
proposition \ref{psu3}, $X$ a la propri\'et\'e $({\bf Q})$ si $X$ n'est pas contenu dans 
${\cal N}_{{\goth g}}$; donc d'apr\`es le corollaire \ref{c2a2}, $X$ est contenu dans la 
r\'eunion de ${\cal X}_{{\goth g}}$ et de ${\cal N}_{{\goth g}}$, d'o\`u le th\'eor\`eme.
\end{proof}

\references

\lastpage

\end{document}